\documentclass[12pt,reqno]{amsart} 
\usepackage{amsthm,amsmath,amssymb,amsfonts,amscd,latexsym,verbatim}
\usepackage{graphicx,bm,psfrag,color} 

\usepackage[scaled=0.92]{helvet}
\usepackage[T1]{fontenc}
\usepackage{textcomp}
 
\newtheorem{Theorem}{Theorem}[section]
\newtheorem{Lemma}[Theorem]{Lemma}
\newtheorem{Proposition}[Theorem]{Proposition}

\newtheorem{Definition}[Theorem]{Definition}

\newcommand{\cbrac}{)\hspace{-1.5mm}(} 
\newcommand{\const}{\Box} 
\renewcommand{\P}{\mathbb P}        
\newcommand{\ra}{\rightarrow}
\newcommand{\lra}{\longrightarrow}
\newcommand{\la}{\leftarrow}

\renewcommand{\le}{\leqslant} 
\renewcommand{\ge}{\geqslant} 
\newcommand{\complex}{{\mathbb C}} 
\newcommand{\Sym}{\text{Sym}}

\newcommand{\cX}{\mathcal X}
\newcommand{\tcY}{{\widetilde{\mathcal Y}}^\circ}
\newcommand{\cY}{\mathcal Y}
\newcommand{\SYS}{\mathfrak S}

\newcommand{\demo}{\noindent {\sc Proof.}\;}

\newcommand{\NN}{\mathbf N}
\newcommand{\ZZ}{\mathbf Z}
\newcommand{\ID}{\text{ID}} 
\def\bbone{{\mathchoice {\rm 1\mskip-4mu l} {\rm 1\mskip-4mu l}
{\rm 1\mskip-4.5mu l} {\rm 1\mskip-5mu l}}}
\begin{document} 
\title[Quadratic involutions]{Quadratic involutions \\ on \\ binary forms} 
\maketitle 
\centerline{Abdelmalek Abdesselam and Jaydeep Chipalkatti} 
\pagestyle{myheadings}
\markboth{{\sc Quadratic Involutions}}{{\sc Abdesselam and Chipalkatti}}

\bigskip \bigskip 

\parbox{12cm} {\small 
{\sc Abstract:} 
There is a classical geometric construction which uses a binary quadratic form to define 
an involution on the space of binary $d$-ics. We give a complete characterisation of a general 
class of such involutions which are definable using compound transvectant formulae. We also 
study the associated varieties of forms which are preserved by such involutions. 
Along the way we prove a recoupling formula for transvectants, which is used to 
deduce a system of equations satisfied by the coefficients in these involutions. } 

\bigskip 
\thispagestyle{empty} 

\parbox{12cm}{\small 
} 

\bigskip  

\setcounter{tocdepth}{1} 
\tableofcontents

\section{Introduction} 
Given a smooth conic $C$ in the projective plane $\P^2$, a point in $\P^2 \setminus C$  
will define an involution (i.e., a degree $2$ automorphism) on $C$. 
Several familiar objects in the invariant theory of binary forms 
(such as the quartic catalecticant or the Hermite invariant) can be defined as sets of divisors 
which are fixed by such an involution. In this 
paper we study a wide class of such involutions and the corresponding fixed loci. 

We begin with an elementary introduction to the subject. The main results are described in 
\S\ref{summary.results} after the required notation is available. Many of the proofs involve 
elaborate calculations using the graphical or symbolic method, 
as in~\cite{spinnet, Advances}. 
In such cases, as far as possible we have stated the formula which is to be of immediate use, 
and relegated its derivation to a later section. Our object is to ensure that the reader 
who is not familiar with this 
method should be able to follow the argument without loss of continuity. 

We refer the reader to~\cite{Glenn, GY,Salmon} for classical introductions to the invariant 
theory of binary forms, and~\cite{Dolgachev, Olver, Procesi, Sturmfels} for more 
modern accounts. 

\subsection{Representations of $SL_2$} 
Throughout, the base field will be $\complex$ (complex numbers). Let $V$ denote a 
two-dimensional $\complex$-vector space. For a nonnegative 
integer $m$, let $S_m = \Sym^m \, V$ denote the $m$-th symmetric power. 
If $x = \{x_1,x_2\}$ is a basis of $V$, then $S_m$ can be identified with the space of 
binary forms of order $m$ in the variables $x$. 
The $\{S_m : m \ge 0\}$ are a complete set of finite-dimensional irreducible 
representations of $SL(V)$ (see e.g.,~\cite[\S  I.9]{Knapp}). 

Following a notation introduced by Cayley, we will write the binary form 
$\sum\limits_{i=0}^d \, a_i \, \binom{d}{i} \, x_1^{d-i} \, x_2^i$ as 
$(a_0,\dots,a_d \cbrac x_1,x_2)^d$. 

\subsection{Transvectants} 
Given integers $m,n \ge 0$ and $0 \le r \le \min (m,n)$, we have a transvectant morphism 
(see~\cite{AIF}) 
\[ S_m \otimes S_n \longrightarrow S_{m+n-2r}. \] 
If $A,B$ are binary forms of orders $m,n$ respectively, the image of $A \otimes B$ via 
this morphism is called their $r$-th transvectant; it will be denoted by $(A,B)_r$. We have 
an explicit formula 
\[ 
(A,B)_r = \frac{(m-r)! \, (n-r)!}{m! \, n!} 
\sum\limits_{i=0}^r \, (-1)^i \binom{r}{i} \, 
\frac{\partial^r A}{\partial x_1 ^{r-i} \, \partial x_2^i} 
\frac{\partial^r B}{\partial x_2^{r-i} \, \partial x_1^i} \, . 
\] 
In the notation of symbolic calculus, if $A = \alpha_x^m$ and $B = \beta_x^n$, then 
$(A,B)_r = (\alpha \, \beta)^r \, \alpha_x^{m-r} \, \beta_x^{n-r}$.  
\subsection{} 
Consider the Veronese embedding 
\[ v: \P V \lra \P S_2, \quad [\ell] \lra [\ell^2];  \] 
whose image $C$ is a smooth conic in $\P^2$. If $(a_0,a_1,a_2 \cbrac x_1,x_2)^2$ is 
identified with $[a_0,a_1,a_2] \in \P^2$, then $C$ is defined by the 
equation $a_1^2 = a_0 \, a_2$. 

Fix a point $q \in \P^2 \setminus C$. Given $t \in C$, draw the line 
$\overline{q \, t}$, and let $\sigma_q(t)$ denote the other point where the line intersects $C$. 
This defines an order $2$ automorphism 
\[ \sigma_q: C \lra C, \qquad t \lra \sigma_q(t). \] 
The two intersection points of $C$ with the polar line of $q$ are the 
fixed points of this automorphism. 

Assume that $q$ corresponds to $Q = (q_0,q_1,q_2 \cbrac x_1,x_2)^2 \in S_2$. 
Define 
\[ \Delta_Q = - 2 \, (Q,Q)_2 = 4 \, (q_1^2 - q_0 \, q_2). \] 
(The normalisation is chosen such that $\Delta_{x_1x_2}=1$.) We have $\Delta_Q \neq 0$, 
since $[Q] \notin C$. We will merely write $\Delta$ for $\Delta_Q$ if no confusion is likely. 
\begin{Lemma} \sl 
If the point $t \in C$ corresponds via $v$ to $\ell \in V$, then $\sigma_q(t)$ corresponds to 
$(Q,\ell)_1$. 
\end{Lemma} 
\demo It is enough to show that the three forms 
$Q, \ell^2, {(Q,\ell)_1}^2$ are linearly dependent, and hence the corresponding 
points in $\P^2$ are collinear. By $SL_2$-equivariance, we may take 
$\ell = x_1$. Then $(Q,x_1)_1 = - (q_1 \, x_1 + q_2 \, x_2)$, and 
it is immediate that 
\begin{equation} 
{(Q,x_1)_1}^2 - q_2 \, Q = (q_1^2 - q_0 \, q_2) \, x_1^2, 
\label{formula.sigmaq} \end{equation} 
which proves the claim. \qed 

\subsection{} \label{section.defn.sigmaQF} 
If $\ell \in V$, then a simple calculation shows the identity 
\begin{equation} 
(Q,(Q,\ell)_1)_1 = \frac{1}{4} \, \Delta \, \ell. 
\label{formula.QQl} \end{equation} 
Let $F \in S_d$ be a nonzero binary $d$-ic, which factors into linear forms as 
$F = \prod\limits_{i=1}^d \, \ell_i$. Define 
\begin{equation} \sigma_Q(F) = 2^d \, \prod\limits_{i=1}^d \, (Q,\ell_i)_1. 
\label{sigmaQF.l} \end{equation} 
By formula~(\ref{formula.QQl}), 
\begin{equation} 
\sigma_Q^2(F) = \Delta^d \, F. \label{sigma2} \end{equation} 
Now $F$ corresponds to the divisor $A = \sum \, v([\ell_i])$ 
on $C$, and $\sigma_Q(F)$ to its image $\sigma_q(A) = \sum \, \sigma_q([\ell_i^2])$. 
The divisor $A$ is said to be {\sl in involution} with respect to $q$, if $\sigma_q(A) = A$. 
The following diagrams illustrate typical divisors in involution for orders $6$ and $7$. 

\smallskip 

\[ 
\parbox{12cm}{\psfrag{q}{$q$}
\psfrag{C}{$ C$}
\includegraphics[width=12cm]{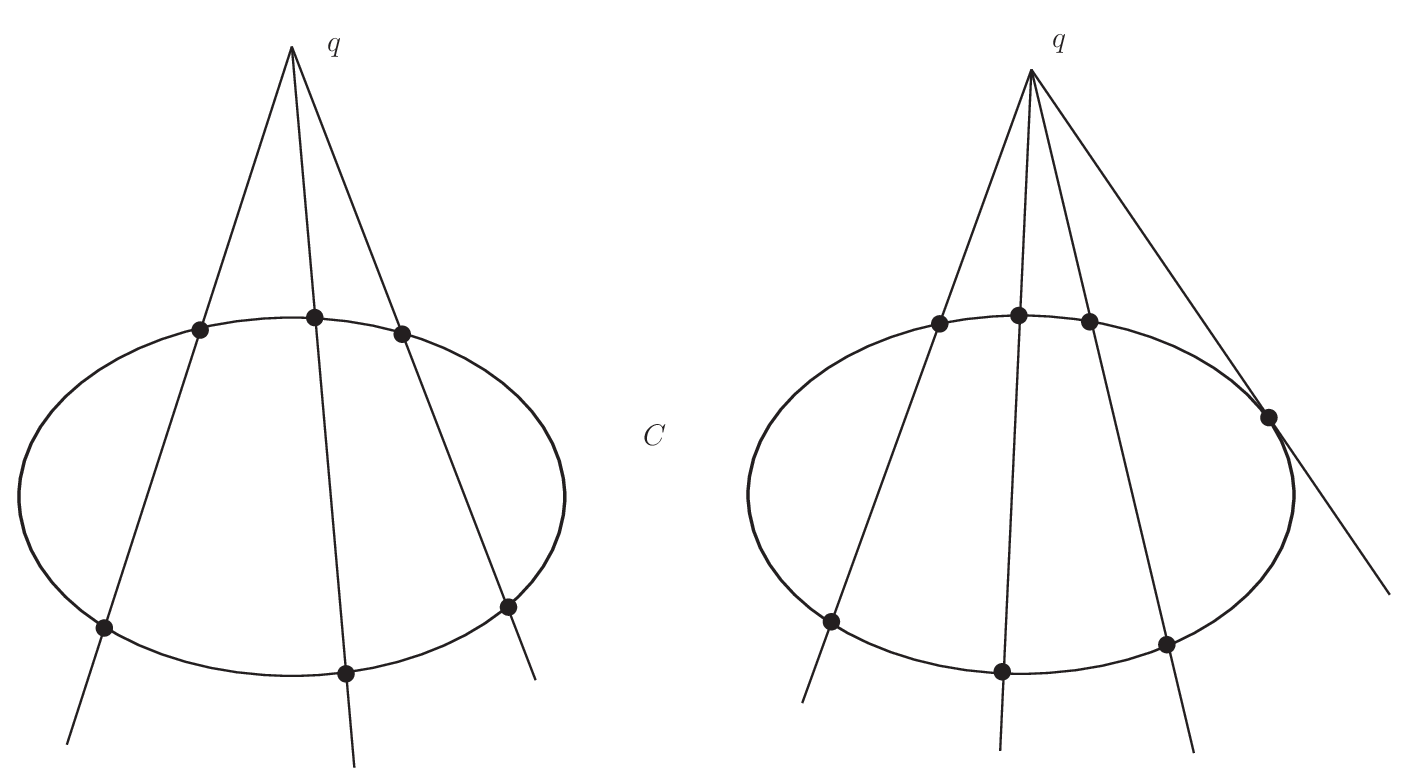}}
\] \label{conics.inv} 

\bigskip 

\subsection{} Define 
\begin{equation} 
\cX^\circ_d = \{A \in \P S_d: \text{$A$ is in involution with respect to some $q \notin C$} \},   
\label{defn.cX} \end{equation} 
and let $\cX_d \subseteq \P S_d$ denote its Zariski closure. It is an $SL_2$-invariant 
irreducible projective subvariety of $\P^d$. 

If $d \le 4$, then given a general set of $d$ points on $C$, one can always complete 
the diagram to find a $q$, hence $\cX_d = \P^d$. Assume $d \ge 5$. 
If $d = 2n$, then a typical point in $\cX_d^\circ$ can be constructed as follows: choose 
an arbitrary point $q$ away from $C$, a degree $n$ divisor $B$, and let $A = B + \sigma_q(B)$. 
If $d = 2n+1$, then choose $q,B$ as above, together with 
any one point $t \in C$ such that $\overline{qt}$ is tangent to $C$, and let 
$A = B + \sigma_q(B) + t$. In either case, a parameter count shows that $\dim \cX_d = n+2$. 

The hypersurfaces $\cX_5$ and $\cX_6$ are respectively of degrees $18$ and $15$. 
Their defining equations are given by well-known skew-invariants of 
binary quintics and sextics (see~\cite[\S 260]{Salmon} and \cite{JC_Hermite}). 
No such equations seem to be known for higher values of $d$. 

\subsection{Canonical form} Let us write $\const$ for an arbitrary scalar which 
need not be precisely specified. Suppose that the divisor $A$ is in involution with 
respect to $q \in \P^2 \setminus C$. Then $A$ is a sum of pairs of the form $t + \sigma_q(t)$, 
together with some points $w$ such that $w = \sigma_q(w)$. 

By a change of variables, we may assume that 
$q = [Q]$ for $Q = x_1 \, x_2$. If $\ell = l_1 \, x_1 + l_2 \, x_2$, then 
$(Q, \ell)_1 = -\frac{1}{2} \, (l_1 \, x_1 - l_2 \, x_2)$. 
If $\ell = \const \, (Q,\ell)_1$, then $\ell = \const \, x_1$ or $\const \, x_2$. If 
$\ell \neq \const \, (Q,\ell)_1$, then $\ell \, (Q,\ell)_1$ is a form in $x_1^2$ and $x_2^2$. 
Alternately, $\const \, x_1^2 + \const \, x_2^2$ can be factored as  $\ell \, (Q,\ell)_1$. 
We have proved the following: 
\begin{Proposition} \sl 
A divisor $A$ is in involution, if and only if, up to a linear change of variables, 
it corresponds to a binary form which can be written as 
\[ x_1^r \, x_2^s \, \times \text{a form in $x_1^2$ and $x_2^2$},  \] 
for some $r,s$. 
\label{proposition.canonical1} \end{Proposition} \qed 

\noindent Later we will prove a similar result for a more general class of involutions.

\section{The system $\SYS(d)$}    \label{system.sigmad} 
\subsection{} We begin by generalising the map $\sigma_Q$, and subsequently the 
notion of an involution. The first step is to rewrite the expression (\ref{sigmaQF.l}) 
for $\sigma_Q(F)$ in terms of transvectants of only $Q$ and $F$, without involving the 
factors of $F$. This is done in \S\ref{proof.geometric.formula}; here we only 
state the result. 

\smallskip 

Let $n = \lfloor d/2\rfloor$. Then we have an expansion 
\begin{equation} 
\sigma_Q(F) = \sum\limits_{i=0}^{n} \; g_i \, \Delta^i \, (Q^{d-2i},F)_{d-2i}, 
\label{geom.involution} \end{equation} 
where 
\[ g_i = \frac{2^{d-2i} \cdot d! \, (d-i)! \, (2d-4i+1)!}{i! \, (d-2i)!^2 \, (2d-2i+1)!}. 
\] 

\smallskip 

Since the construction of $\sigma_Q(F)$ is covariant in $Q$ and $F$, a 
result of Gordan (see~\cite[\S103]{GY}) implies that it should be expressible in 
terms of compound transvectants of the two arguments. 
Thus, it follows~{\it a priori}, that an identity such as 
(\ref{geom.involution}) should exist for \emph{some} rational numbers $g_i$. 
\subsection{} \label{section.Sd} 
All of this suggests the following construction. 
Let $z = (z_0,\dots,z_n)$ be a sequence of complex numbers, and consider the function 
\[ \sigma_{Q,z}: S_d \lra S_d, \quad 
F \lra \sum\limits_{i=0}^n \, z_i \, \Delta^i \, (Q^{d-2i},F)_{d-2i}. \] 
One should like to write a set of equations in $z_i$  
which encodes the condition that $\sigma_{Q,z}$ be involutive. Now, 
\[ \begin{aligned} 
\sigma_{Q,z}^2(F) = & \sum\limits_{i=0}^n \, z_i \, \Delta^i \, 
\sum\limits_{j=0}^n \, z_j \, \Delta^j \, (Q^{d-2j},(Q^{d-2i},F)_{d-2i})_{d-2j} \\ 
= & \sum\limits_{0 \le i, j \le n} \, z_i \, z_j \, \Delta^{i+j} \, 
\underbrace{(Q^{d-2j},(Q^{d-2i},F)_{d-2i})_{d-2j}}_{(\star)} . 
\end{aligned} \] 
One can rewrite  $(\star)$ by expanding the compound transvectant 
expression $(Q^\bullet,(Q^\bullet,F))$ into a sum of terms of the 
form $\Delta^\bullet (Q^\bullet,F)$. A general formula along these lines is proved 
in~\S\ref{section.omega}, where the reader will find the definition of the 
rational numbers $\omega$ which are needed below. Here we only need to state the result to be 
used. There is an expansion 
\begin{equation} (\star) = \sum\limits_{t} \, 
\alpha_{i,j}^{(t)} \, \Delta^{d-i-j-\frac{t}{2}} \, (Q^t,F)_t, 
\label{expansion.star} \end{equation} 
where the sum is quantified over all even $t$ in the range 
\begin{equation} 
2 |i-j| \le t \le \min \, \{d,2(d-i-j)\}, \label{range.t} \end{equation} 
and we have written 
$\alpha_{i,j}^{(t)} = \omega(d-2j,d-2i;d-2i,d-2j; t)$ for brevity. 
If $t$ does not lie in the range (\ref{range.t}), then define $\alpha_{i,j}^{(t)}$ to be zero. 
In particular, 
$\alpha_{i,j}^{(0)} = 0$, unless $i=j$. Thus 
\[ \sigma_{Q,z}^2(F) = \sum\limits_{\stackrel{t=0}{t \; \text{even}}}^d \; 
\left\{ 
\Delta^{d-\frac{t}{2}} \, (\sum\limits_{0 \le i,j \le n} \, \alpha_{i,j}^{(t)} \, z_i z_j) \, 
(Q^t,F)_t \right\}. \] 
Now we require that the coefficient of $\Delta^{d-\frac{t}{2}}$ be $1$ for $t=0$, and 
vanish for $t \neq 0$, which would force 
\begin{equation} \sigma_{Q,z}^2(F) = \Delta^d \, F. 
\label{identity.inv} \end{equation} 
This gives the following system of $n$ homogeneous 
quadratic equations 
\begin{equation} 
\sum\limits_{0 \le i,j \le n} \, \alpha_{i,j}^{(t)} \, z_i z_j = 0, \quad 
(t=2,4,\dots,2n),  
\label{system.inv} \end{equation} 
together with the condition 
\begin{equation} 
\sum\limits_{i=0}^n \, \alpha_{i,i}^{(0)} \, z_i^2 =1. 
\label{condition.inv.one} \end{equation} 
The combined set~(\ref{system.inv}) and (\ref{condition.inv.one}) will be denoted by 
$\SYS(d)$. For instance, the system $\SYS(6)$ is comprised of 
\[ \left. 
\begin{array}{l} 
- \frac{25}{20328} \, z_0^2+\frac{5}{3234} \, z_0 \, z_1-\frac{1}{2058} \, z_1^2
+\frac{22}{735} \, z_1 \, z_2+\frac{11}{210} \, z_2^2+2 \, z_2 \, z_3 \\ \\ 
\frac{5}{1331} \, z_0^2-\frac{15}{847} \, z_0 \, z_1+
\frac{5}{121} \, z_0 \, z_2-\frac{69}{5390} \, z_1^2-\frac{2}{77} \, z_1 \,z_2+ 
2 \, z_1 \, z_3+\frac{2}{5} \, z_2^2 \\ \\ 
- \frac{5}{2541} \, z_0^2+\frac{4}{165} \, z_0 \, z_1-\frac{7}{33} \, z_0 \, z_2+
2 \, z_0 \, z_3-\frac{1}{35} \, z_1^2+\frac{2}{15} \, z_1 \, z_2
\end{array} \right\} =0, \] 
together with 
\[ \frac{1}{6468} \, z_0^2 + \frac{11}{22050} \, z_1^2 + \frac{1}{75} \, z_2^2 + z_3^2 =1. 
\] 

\subsection{} 
A sequence $z=(z_0,\dots,z_n)$ will be called an involutor if it satisfies $\SYS(d)$. In 
particular, $g = (g_0,\dots,g_n)$ will be called the {\sl geometric} involutor. 
It is clear that if $z$ is an involutor, then 
so is $- z = (-z_0,\dots,-z_n)$. If $d$ is even, then $(0,\dots,0, \pm 1)$ will be called the 
{\sl improper} involutors. (In these cases $\sigma_{Q,z}$ merely multiples $F$ by a 
scalar, i.e., it is the identity map at the level of divisiors on $C$.) 

\subsection{A summary of results} \label{summary.results} 
Since (\ref{system.inv}) is a system of $n$ homogeneous quadratic 
equations in $n+1$ variables, B{\'e}zout's theorem implies that the number 
of homogeneous solutions, if finite, should be at most $2^n$. 
After imposing condition~(\ref{condition.inv.one}), one expects at most $2^{n+1}$ affine 
solutions. We programmed the system $\SYS(d)$ in {\sc Maple}, and found that 
for the first few values of $d$, there are always precisely $2^{n+1}$ involutors. This 
prompted us to ask whether it should be possible to write down all the solutions to $\SYS(d)$. 
We will carry this out in \S\ref{signseq.involutors} below. 

The key idea is to introduce a graphically motivated basis $\{{\mathcal O}_i\}$ for the space 
$\text{Hom}_{SL_2}(S_d(S_2) \otimes S_d, S_d)$. It turns out that an involutor arises 
naturally from a {\sl sign sequence} (see Definition~\ref{defofs}), and moreover there is 
an explicit bijection between sign sequences and involutors. This is formally 
stated in Theorem~\ref{MainTheorem}. 

Given an involutor $z$, one can define a subvariety $\cY(z) \subseteq \P S_d$ in 
analogy with~\S\ref{section.defn.sigmaQF}. We initiate a study of this varieties 
in \S\ref{section.varieties}--\ref{section.vartheta.locus}. 
In particular, Theorem~\ref{can.form} will show that a point in $\cY(z)$ admits 
a canonical form which can be read off from the corresponding sign sequence. It is evident 
from the examples in \S\ref{neg.z}--\ref{d.even.ex} and \S\ref{section.vartheta.locus} 
that these varieties display a wide range of geometries depending on the sign sequence, and they 
provide ample matter for further study. 

We prove a general recoupling formula for transvectants in Theorem~\ref{Trecoupling}. 
In brief, it rewrites a compount transvectant of the form $(A,(B,C))$ as a sum 
of terms (with coefficients) of the form $((A,B),C)$. Then we specialise to the 
case where $A$ and $B$ are powers of a quadratic form $Q$. The resulting coefficients 
(denoted by $\omega$) are precisely the ones needed to build $\SYS(d)$. 

\section{Sign sequences and Involutors} \label{signseq.involutors} 
\subsection{} In this section we will describe a complete classification of all involutors. 
The following definition and formulae are justified by Theorem~\ref{MainTheorem}. 
Recall that $n = \lfloor d/2 \rfloor$. 
\begin{Definition}\label{defofs} \rm A sign sequence for $d$ is of the form 
$s = (s_0,s_1,\dots,s_d)$, where $s_i = \pm 1$, and $s_{d-i} = (-1)^d \, s_{i}$. 
\end{Definition} 
It is clear that the segment $(s_0,\dots,s_n)$ can be made up arbitrarily, and then it 
determines the rest. Thus there are $2^{n+1}$ sign sequences for $d$. 
We will write 
$\pm$ for $\pm 1$ if no confusion is likely. Define $\gamma$ to be the 
alternating sign sequence $(-,+,\dots,-,+)$ if $d$ is odd, and $(+,-,+,\dots,-,+)$ if $d$ is even. 

\subsection{} \label{z.formula.sec}
Given a sign sequence 
$s$, and an index $i$ such that $0 \le i \le n$, let 
\[ E_{1,i} = \frac{d! \, (2d-4i+1)!}{2^{2i-1} \, (d-2i)!^2},  \] 
and 
\[ E_{2,i} = \sum \left\{ s_\ell \, m_\ell \, (-1)^q \, 
\frac{(d-2e)! \, (d-i-e)!}{(2d-2i-2e+1)! \, (i-e)! \, p! \, q! \, (\ell-p)! \, (d-\ell-q)!} \right\}, \] 
where the sum is quantified over all integer quadruples $(e,\ell,p,q)$ such that 
\[ 0 \le e \le i, \quad 0 \le \ell \le n, \quad 0 \le p \le \ell, \quad 0 \le q \le d-\ell, \quad 
p+q=d-2e. \] 
Here $m_\ell$ is defined to be $\frac{1}{2}$ if $d = 2 \, \ell$, and $1$ otherwise. Define 
\begin{equation} z_i(s) = E_{1,i} \, E_{2,i}. 
\label{formula.zis} \end{equation} 
\begin{Theorem} \label{MainTheorem}\sl 
With notation as above, $z(s) = (z_0(s),\dots,z_n(s))$ is an involutor. Moreover, every involutor 
arises in this way from a unique sign sequence. 
\end{Theorem} 
The proof will be given in~\S\ref{main.proof.sec}.
The following proposition lists some properties of this 
correspondence. 
\begin{Proposition}\label{special.prop} \sl With notation as above, 
\begin{enumerate} 
\item 
If $z=z(s)$, and $s'$ is the sign sequence such that $s'_i = - s_i$, then 
$z(s') = -z$. 
\item 
The geometric involutor corresponds to $\gamma$. 
\item 
When $d$ is even, the improper involutors correspond to the sequences 
$(+,+,\dots,+)$ and $(-,-,\dots,-)$. 
\end{enumerate} 
\end{Proposition} 
Part (1) is obvious from the definition of $E_{2,i}$. The rest will be proved in \S\ref{special.sec}. 

\smallskip 

For instance, if $d=4$, then $(+,-,-,-,+)$ corresponds to the involutor $(4,48/7,-1/5)$, and 
$(+,-,+,-,+)$ corresponds to the geometric involutor $(16,24/7,1/5)$. 

\section{Varieties of forms in involution} \label{section.varieties} 
\subsection{} 
Consider the product $\P S_2 \times \P S_d$ with projections $\pi_1,\pi_2$ onto 
$\P S_2$ and $\P S_d$ respectively. Given an involutor $z$, define the locus 
$\tcY(z) \subseteq \P S_2 \times \P S_d$ as the set of pairs 
$\langle Q,F \rangle$ satisfying $\Delta_Q \neq 0$, and 
\begin{equation} 
\begin{array}{cl} 
\sigma_{Q,z}(F) = \Delta^{d/2} \, F, & \text{if $d$ is even,} \\ 
\left[\sigma_{Q,z}(F)\right]^2 = \Delta^d \, F^2, & \text{if $d$ is odd.}
\end{array} \label{equation.xtilda} \end{equation} 
If (\ref{equation.xtilda}) holds, then we will say 
that $[Q] \in \P^2$ is a centre of involution for $[F]$ with respect to $z$. For instance, 
in the diagrams on page \pageref{conics.inv}, the point $q$ is such a centre 
with respect to $\gamma$ for the forms corresponding to the divisors shown. However, 
if $z \neq \gamma$, then it is not clear to us whether there is any hidden `geometry' in 
the relationship defined by (\ref{equation.xtilda}). 

Define $\cY^\circ(z) \subseteq \P S_d$ to be the image $\pi_2(\tcY(z))$, and 
let $\cY(z) \subseteq \P S_d $ denote its the Zariski closure. We may equally well 
denote these loci by $\cY(s)$ etc.~by referring to the sign sequence, and further shorten them 
to $\tcY, \cY$ etc.~if no confusion is likely. 

It is clear that if $d$ is odd, then $\cY(z) = \cY(-z)$. This may or may not hold for 
even $d$ (see \S\ref{d.even.ex} below). 

\subsection{The relation between $\cX$ and $\cY$} 
Assume that $A$ is a divisor on $C$ in involution with respect to 
$q = [Q] \in \P^2 \setminus C$. If $F \in S_d$ represents $A$, then 
$\sigma_{Q,\gamma}(F) = c \, F$ for some constant $c$. Formula (\ref{sigma2}) in 
\S\ref{section.defn.sigmaQF} implies that $c^2 = \Delta_Q^d$, i.e., $c = \pm \Delta_Q^{d/2}$. 
Hence, 
\[ \cX^\circ_d = \cY^\circ(\gamma) \cup \cY^\circ(-\gamma). \] 
For odd $d$, this reduces to $\cX^\circ_d = \cY^\circ(\gamma)$. 

\subsection{Canonical Forms} 
It turns out that an element in $\cY$ admits a canonical form analogous to 
Proposition~\ref{proposition.canonical1}. Let $s$ be a sign sequence for $d$, and 
$z =z(s)$. Assume $Q = x_1 \, x_2$ after a change of variables, and consider the condition 
\begin{equation} \langle x_1 \, x_2, F \rangle  \in \tcY(z). 
\label{condition.tcY} \end{equation} 
\begin{Theorem} \sl \label{can.form} 
\begin{enumerate} 
\item 
Assume $d$ to be even. Then (\ref{condition.tcY}) holds, if and only if $F$ is a linear combination of 
terms in the set 
\[ \{x_1^{d-i} \, x_2^i: s_i = 1\}. \] 
\item 
Assume $d$ to be odd. Then (\ref{condition.tcY}) holds, if and only if $F$ is a linear combination of 
terms in either one of the sets 
\[ \{x_1^{d-i} \, x_2^i: s_i = 1\} \quad \text{or} \quad 
\{x_1^{d-i} \, x_2^i: s_i = -1\}. \] 
\end{enumerate} \label{theorem.canonical2} \end{Theorem} 
For instance, if $s = (+,-,-,+,-,-,+)$, then such an $F$ is of the form 
$\Box \, x_1^6 + \Box \, x_1^3 \, x_2^3 + \Box \, x_2^6$. If $s = (+,-,-,+,+,-)$, then 
it is of the form 
\[ \Box \, x_1^5 + \Box \, x_1^2 \, x_2^3 + \Box \, x_1 \, x_2^4, \quad \text{or} \quad 
\Box \, x_2^5 + \Box \, x_2^2 \, x_1^3 + \Box \, x_2 \, x_1^4. \] 
In general, $\cY^\circ(z)$ is a union of $SL_2$-orbits of such forms. 

\smallskip 

\demo See~\S\ref{can.form.proof}. \qed 

\subsection{} \label{neg.z} 
Let $p(s)$ be the number of $+$ signs in $s$. If $d$ is odd, then $p(s) = \frac{1}{2}(d+1)$. 
By Theorem~\ref{theorem.canonical2}, each fibre of the projection 
$\tcY \lra \P S_2 \setminus C$ has dimension $p(s)-1$. Hence 
\begin{equation} \dim \cY \le \min \, \{ \, p(s)+1,d \, \}. 
\label{dim.cY} \end{equation} 
This inequality may be strict. For instance, let $d=3$, and  $s = (+,+,-,-)$. 
A typical element in $\cY$ can be written as $x_1^2 \, \ell$ up to a change of variables, 
hence $\cY \subseteq \P^3$ is the discriminant surface of degree $4$. 

In general, if $d$ is odd, and 
\[ s = (\underbrace{+,\dots,+}_{\text{$\frac{d+1}{2}$ times}},-,\dots,-), \] 
then the canonical form shows that a typical element in $\cY$ is a binary $d$-ic with a root of 
multiplicity $\frac{d+1}{2}$. Hence $\dim \cY = \frac{d+1}{2}$, and (\ref{dim.cY}) is strict. 

However, if $z = \gamma$, then (\ref{dim.cY}) is an equality for $d \ge 3$. Indeed, 
$p(\gamma) = \frac{d+2}{2}$ or $\frac{d+1}{2}$ according to whether $d$ is even or odd, 
and the right-hand side reduces to $n+2$. 

\subsection{} 
Let $d=4$. A binary quartic $F$ has covariants 
\[ A_F = (F,F)_4, \quad B_F = (F,(F,F)_2)_4. \] 
(See~\cite[\S 89]{GY}; however our notation differs from theirs.) 
Usually, $B_F$ is called the catalecticant. If $F$ has distinct roots, its $j$-invariant is 
defined to be (cf.~\cite[\S 171]{GY}) 
\[ j(F) = \frac{A_F^3}{A_F^3 - 6 \, B_F^2}. \] 
Recall that two binary quartics (with distinct roots) are in the same $SL_2$-orbit, exactly 
when they have the same $j$-invariant (see e.g.,~\cite[Example 10.12]{Harris}). 
Now let 
\[ s = (+,-,-,-,+), \quad t = (-,+,-,+,-), \] 
so the canonical forms are respectively 
\[ G_s = \Box \, x_1^4 + \Box \, x_2^4, \qquad 
G_t = x_1 \, x_2 \, (\Box \, x_1^2 + \Box \, x_2^2). \] 
A straightforward calculation shows that $j(G_s) = j(G_t) = 1$, hence 
$\cY(s) = \cY(t)$ is the cubic hypersurface defined by the equation $B_F=0$. 
(Classically, these were called the harmonic binary quartics.) 

In general, for even $d$ and $s = (+,-,\dots,-,+)$, the variety $\cY(s)$ is 
the chordal threefold (union of secant lines) of the rational normal $d$-ic curve. 

\subsection{} Let $d$ be even, and $s = (-,\dots,-,+,-,\dots,-)$. The canonical form is 
$(x_1 \, x_2)^{\frac{d}{2}}$, i.e., $\cY$ is the variety of $d$-ics which are expressible as 
powers of quadratic forms. It is shown in~\cite{Advances}, that its ideal is 
generated by (an explicitly given) list of cubic polynomials. 

\subsection{} \label{d.even.ex} 
Let $s = (+,-,-,-,+)$ and $s' = (-,+,+,+,-)$. A general quartic can be written as 
$x_1 \, x_2 \, (\Box \, x_1^2 + \Box \, x_1 \, x_2 + \Box \, x_2^2)$ up to a change of 
variables, hence $\cY(s') = \P^4$. Since $\cY(s)$ is a threefold, $\cY(s) \neq \cY(s')$. 

By contrast, let $d=2$, and $u = (+,-,+), u' = (-,+,-)$. 
Then $\cY(u,2) = \cY(u',2) = \P^2$. 

\subsection{} In general, one should like to know formulae for the dimension and degree of 
$\cY$ as a function of $s$; moreover, it would be of interest to be able to write down 
a set of $SL_2$-equivariant defining equations for $\cY$. Our next example shows 
(if nothing else), that such questions can be rather involved. 

Let $d = 6$, and $s = (-,-,+,-,+,-,-)$. Up to a change of variables, a form in $\cY^\circ(s)$ can 
be written as $G_s = x_1^2 \, x_2^2 \, (x_1^2 + x_2^2)$. This corresponds to the 
divisor $2 \, [\ell_1] + 2 \, [\ell_2] + [\ell_3] + [\ell_4]$, where 
\[ \ell_1 = x_1, \quad \ell_2 = x_2, \quad \ell_3 = x_1 + \sqrt{-1} \, x_2, \quad 
\ell_4 = x_1 - \sqrt{-1} \, x_2. \] 
Thus we can describe $\cY^\circ$ as the set of divisors 
$2 \, p_1 + 2 \, p_2 + p_3 + p_4$ on $C$, such that the intersection point $[Q]$ 
of tangents at $p_1,p_2$ falls on the line $\overline{p_3 \, p_4}$. 
The entire configuration is determined 
by $Q$ and an arbitrary line through it, hence $\dim \cY = 3$. Let 
\[ B = \Sym^\bullet \, (S_6) = \bigoplus\limits_{r \ge 0} \, \Sym^r \, (S_6) \] 
denote the c{\"o}ordinate ring of $\P S_6$. (Since each $S_m$ 
is canonically isomorphic to its dual, here it is unnecessary to distinguish between the two.) 
We calculated the defining ideal $J \subseteq B$ of $\cY \subseteq \P^6$ 
using straightforward elimination in Macaulay-2. 
It turns out that the degree of $\cY$ as a variety is $18$. Furthermore, 
$J$ is generated by one form in degree $4$, and $36$ forms in degree $5$. 
More precisely, its minimal resolution begins as 
\[ 0 \la B/J \la B \la B(-4) \otimes M_1 \oplus B(-5) \otimes M_{36} \la \dots \] 
where $M_i$ is an $SL_2$-representation of dimension $i$. 

Clearly $M_1 \simeq S_0$. Since the complete minimal system of binary sextics is known 
(see~\cite[\S 132-134]{GY}), it is a mechanical task to find the irreducible decomposition of 
$M_{36}$. Indeed, $M_{36} \subseteq \Sym^5 (S_6)$, and the latter can be decomposed into 
irreducibles by the Cayley-Sylvester formula (see~\cite[Corollary 4.2.8]{Sturmfels}). 
Thus we only need to identify those degree $5$ covariants of sextics which vanish 
when specialised to $G_s$. The decomposition turns out to be 
\[ M_{36} \simeq S_{14} \oplus S_{10} \oplus S_6 \oplus S_2. \]  
Let $\theta_{rn}$ stand for the covariant of degree-order $(r,n)$ as given in the table 
on~\cite[p.~156]{GY}; for instance, $\theta_{38} = (F,(F,F)_4)_1$. 
We will explicitly write down the covariants corresponding to these representations. They are, 
in degree $4$, 
\begin{equation} 
\text{order $0$} \leadsto \; 7 \, \theta_{20}^2 - 50 \, \theta_{40},  
\label{eq.Y.deg4} \end{equation} 
and in degree $5$, 
\begin{equation} 
\begin{aligned} 
\text{order $14$} & \leadsto \; 
20 \, \theta_{16} \, \theta_{24}^2 -21 \, \theta_{16} \, \theta_{20} \, \theta_{28} + 
10 \, \theta_{16}^2 \, \theta_{32} + 90 \, \theta_{28} \, \theta_{36},  \\ 
\text{order $10$} & \leadsto \; 
2 \, \theta_{16} \, \theta_{20} \, \theta_{24} 
-6 \, \theta_{16} \, \theta_{44} -27 \, \theta_{28} \, \theta_{32}, \\ 
\text{order $6$} & \leadsto \; 
50 \, \theta_{16} \, \theta_{40} -42 \, \theta_{20} \, \theta_{36} 
-105 \, \theta_{24} \, \theta_{32}, \\ 
\text{order $2$} & \leadsto \; \theta_{20} \theta_{32} - 10 \, \theta_{52}. 
\end{aligned} \label{eq.Y.deg5} \end{equation} 
In conclusion, a binary sextic $F$ belongs to $\cY(s)$, if and only if all the covariants in 
(\ref{eq.Y.deg4}) and (\ref{eq.Y.deg5}) vanish on $F$. It is not clear whether 
a simpler set of equations may be found.

\section{The locus of centres of involution} \label{section.vartheta.locus} 
Fix an involutor $z$, and let $F$ be a $d$-ic. The locus of centres of involution for $F$ 
often has interesting geometric structure, especially when $\cY^\circ$ is dense in $\P S_d$.  
We adduce a few such examples. 

Write $Q = (q_0, q_1, q_2 \cbrac x_1,x_2)^2$, and let 
$R = \complex[q_0,q_1,q_2]$ denote the c{\"o}ordinate ring of $\P S_2$. For $F \in S_d$, define 
\[ \vartheta(F)  = \{ [Q] \in \P^2: \Delta_Q \neq 0, \; 
\text{and (\ref{equation.xtilda}) holds} \}. \] 

Suppose that $\ell$ is a linear factor in $F$, and we let $Q = \ell^2$. Then 
it is easy to see that $(Q^d,F)_d$ vanishes. Indeed, by equivariance we may assume 
$\ell = x_1$, and then 
\[ (Q^d,F)_d = \text{constant} \, \times \, x_1^d \; \frac{\partial^d F}{\partial x_2^d} =0. \] 
Since $\sigma_{Q,z}(F)$ and $\Delta_Q$ are zero, (\ref{equation.xtilda}) is satisfied. 
It follows that when $\vartheta(F)$ is nonempty, its Zariski closure will always contain 
the points $[\ell^2] \in C$ corresponding to the linear factors of $F$. Recall that by 
our definition, such points do not count as centres of involution. 

\subsection{} \label{section.d4.conic} 
Let $d=4$, and $s = (+,+,-,+,+)$, then $z = (-12,24/7,3/5)$ by the formulae 
in \S\ref{signseq.involutors}. 
\begin{Proposition} \sl We have $\cY^\circ(z) = \P^4$. Moreover, for a general binary quartic 
$F$, the locus $\overline{\vartheta(F)}$ is a smooth conic in $\P^2$. 
\end{Proposition} 
\demo 
For arbitrary $Q$ and $F$, there is an identity (see \S\ref{identity.quartic} below) 
\begin{equation} \sigma_{Q,z}(F) - \Delta^2 \, F = -12 \, Q^2 \, (Q^2,F)_4. 
\label{id.QFd4} \end{equation} 
As a result, $[Q] \in \vartheta(F)$, if and only if $[Q] \notin C$ and $(Q^2,F)_4=0$. 
The latter is a quadratic equation in the $q_i$ whose coefficients depend upon $F$. 
Hence $\vartheta(F) \neq \emptyset$ 
for any $F$, and thus $\cY^\circ = \P^4$. The discriminant of this 
quadratic is $B_F = (F,(F,F)_2)_4$ (see \S\ref{claim.BF}). 
Hence $\overline{\vartheta(F)}$ is a smooth conic when $F$ is not harmonic. \qed 

\smallskip 

By what we have said above, this conic passes through the four points on $C$ 
corresponding to the linear factors of $F$. 

\subsection{} \label{identity.quartic} 
Identity~(\ref{id.QFd4}) is an easy consequence of Proposition~\ref{proposition.omega}
in \S\ref{section.omega}. Indeed, the left-hand side is 
\[ -12 \, (Q^2,F)_4 + \frac{24}{7} \, \Delta \, (Q^2,F)_2 - \frac{2}{5} \, \Delta^2 \, F, \] 
which is a rescaled expansion of $(Q^2,(Q^2,F)_4)_0$. 
Identities (\ref{id.QFd6})--(\ref{d6.id1}) can 
all be proved using a similar technique; this is left as an exercise for the reader. 

\subsection{} \label{claim.BF} 
The claim about $B_F$ can be easily established by writing down the determinantal 
expression for the discriminant of a ternary quadratic, followed by a straightforward expansion. 
A more elegant way to do this calculation is to use the macroscopic to microscopic 
rewriting as in~\cite[Eq. 9]{spinnet}.
The embedding $SL(V) \hookrightarrow SL(S_{n-1} \, V) \simeq SL_n$ allows us to reformulate 
the invariant theory of $SL_n$ entirely in terms of that of $SL_2$. 
In the symbolic formalism, that amounts to replacing $n$-ary brackets by homogenized binary
Vandermonde determinants, i.e., products of $\frac{n(n-1)}{2}$ binary brackets.
If we carry out the procedure on the example at hand, we get the well-known 
symbolic expression $B_F = (a \, b)^2 \, (a \, c)^2 \, (b \, c)^2$ for the catalecticant.  

\smallskip 

\subsection{} \label{section.d6.cubic} 
This example is similar to the previous one. Let $d=6$, 
and $s=(+,+,+,-,+,+,+)$, then $z = (40,-180/11,20/7,5/7)$. 
\begin{Proposition} \sl We have $\cY^\circ = \P^6$. Moreover, for a general binary sextic 
$F$, the locus $\overline{\vartheta(F)}$ is a smooth cubic in $\P^2$. 
\end{Proposition} 
\demo 
The argument is parallel to the previous proposition. We have an identity 
\begin{equation} \sigma_{Q,z}(F) - \Delta^3 \, F = 40 \, Q^3 \, (Q^3,F)_6.  
\label{id.QFd6} \end{equation} 
The equation $(Q^3,F)_6=0$ defines a planar cubic curve, whose discriminant 
is a degree $12$ invariant of $F$. 
It is not necessary to calculate it explicitly; for our purposes it would suffice to check that it is not 
identically zero. Specialise to $F = x_1^6 + x_2^6 + x_1^2 \, x_2^4$, when 
\[ (Q^3,F)_6 = q_0^3 + \frac{4}{5} \, q_0 \, q_1^2 + \frac{1}{5} \, q_0^2 \, q_2 + q_2^3. \] 
This is easily seen to be a nonsingular curve. \qed 

\smallskip 

\noindent The following lemma will be used in the next section. 
\begin{Lemma} \sl 
Assume $A$ and $B$ to be nonzero binary forms of the same order. If 
$(A,B)_1 =0$, then the forms are equal up to a multiplicative constant. 
\end{Lemma} 
\demo See \cite[Lemma 2.2]{Goldberg}. \qed 

\subsection{} \label{section.d4.3p} 
In the next two examples in this section, $\vartheta(F)$ is a finite set of 
points. The simplest such case is that of the geometric involution for $d=4$. Let 
$F \in S_4$ correspond to the divisor $A = p_1 + p_2 + p_3 + p_4$ 
consisting of four distinct points on $C$. Now 
$A$ has three centres of involution, namely the pairwise intersections of lines 
\[ \overline{p_1 p_2} \cap \overline{p_3 p_4}, \quad 
\overline{p_2 p_3} \cap \overline{p_1 p_4}, \quad 
\overline{p_1 p_3} \cap \overline{p_2 p_4} \, . \] 
One knows that the ideal of three non-collinear planar points is generated by a 
net of conics; we will see that this ideal can be written down in terms of $F$. (The reader may refer 
to \cite[Ch.~3]{Eisenbud} for generalities on ideals of finite point-sets in $\P^2$.) 
There is an identity 
\begin{equation} 
\sigma_{Q,\gamma}(F) - \Delta^2 \, F = 
Q \, [ \, \underbrace{16 \, (Q^3,F)_4 + \frac{24}{5} \, (Q,F)_2 \, \Delta}_{\alpha} \, ]. 
\label{id.d4.geom} \end{equation} 
Hence, $[Q] \in \vartheta(F)$ implies that $\alpha = 0$. The problem, as usual, 
is that $\alpha$  also vanishes 
if $Q = \ell^2$ with $\ell | F$. We need an expression which would force $\alpha$ to be zero, 
without itself vanishing on such points. The useful identity is 
\begin{equation}  \alpha = - 32 \, ( \, \underbrace{(Q^2,F)_3}_\beta,Q)_1. 
\end{equation} 
Once $F$ is fixed, $\beta$ is of the form $(\varphi_0,\varphi_1,\varphi_2 \cbrac x_1,x_2)^2$, where 
$\varphi_i(q_0,q_1,q_2)$ are homogeneous degree $2$ expressions in the $q_i$. 
Let $I = (\varphi_0,\varphi_1,\varphi_2) \subseteq R$ 
be the ideal generated by the coefficients of $\beta$. 
\begin{Proposition} \sl 
Assume that $F$ has no repeated linear factors. Then the ideal of the three 
points $\vartheta(F)$ is $I$. 
\end{Proposition} 
\demo Let $\Theta$ denote the ideal of $\vartheta(F)$. Assume that $[Q] \in \vartheta(F)$, then 
$(\beta,Q)_1$ vanishes at $[Q]$, and the lemma above implies that 
$\beta = c \, Q$ for some constant $c$. Since 
$[Q] \notin C$, we may assume $Q = x_1 \, x_2$ after a change of variables. If 
\begin{equation} F = (a_0,\dots,a_4 \cbrac x_1,x_2)^4, 
\label{F.gen4} \end{equation} 
then a direct calculation shows that 
$\beta = (Q^2,F)_3 = \frac{1}{2} \, (a_1 \, x_1^2 - a_3 \, x_2^2)$, which forces $c=0$. 
Hence $I \subseteq \Theta$. 

Alternately, assume $\beta = 0$ at $[Q]$. Then the left-hand side of~(\ref{id.d4.geom}) is 
zero. We claim that $[Q] \notin C$, and hence $[Q] \in \vartheta(F)$. Indeed, if 
$[Q] \in C$, then we may assume $Q = x_1^2$ and $F$ is as in (\ref{F.gen4}) above. 
Then $(Q^2,F)_3 = a_3 \, x_1^2 + a_4 \, x_1 \, x_2 =0$. Hence $a_3=a_4 =0$, which 
would force $F$ to have a repeated linear factor. 

Thus the zero locus of $I$ is $\vartheta(F)$. The forms $\{\varphi_0,\varphi_1,\varphi_2\}$ 
must be linearly independent (otherwise $I$ would either define a conic or a scheme 
of length $4$), which implies that $I = \Theta$. \qed 

\smallskip 

The first syzygies of $I$ can also be written down using transvectants. They correspond to 
the identities 
\begin{equation} (\beta,Q)_2 = (\beta, (Q,F)_2)_2 = 0. 
\label{syz.beta} \end{equation} 
Altogether, this accounts for the minimal free resolution 
\[ 0 \la R/I \la R \la R(-2)^3 \la R(-3)^2 \la 0. \] 
The correspondence between transvectant identities and syzygies is 
explained in \cite[\S 1.7,4.1]{JC_Hermite}. 

\subsection{} \label{section.d6.7p} 
This example is very similar to the previous one, hence we will keep the arguments brief. 
Let $d=6$, and $s = (+,+,-,+,-,+,+)$, then $z = (-60,-60/11, 30/7, 3/7)$. 
It turns out that the set $\vartheta(F)$ consists of $7$ points, whose ideal can 
be written down as follows. 

We have an identity 
\[ 
\sigma_{Q,z}(F) - \Delta^3 \, F = -60 \, Q^2 \, 
[ \, \underbrace{(Q^4,F)_6 + \frac{2}{7} \, \Delta \, (Q^2,F)_4}_\mu \, ], 
\] 
and furthermore 
\[ \mu = -2 \, (\lambda,Q)_1, \quad \text{where $\lambda = (Q^3,F)_5$}. \] 
Write $\lambda = (\psi_0,\psi_1,\psi_2 \cbrac x_1,x_2)^2$, and 
$I = (\psi_0,\psi_1,\psi_2) \subseteq R$. 

\begin{Proposition} \sl 
Let $F$ be a general binary sextic. Then the ideal of $\vartheta(F)$ is $I$. 
\end{Proposition} 
\demo The fact that the zero locus of $I$ is $\vartheta(F)$ is proved exactly as in the 
previous proposition. Now specialise $F$ to $x_1^6+x_2^6$, when 
\[ \lambda = (-q_1 \, q_2^2, \frac{q_0^3 - q_2^3}{2}, q_0^2 \, q_1 \cbrac x_1,x_2)^2. \] 
It is easy to see that the coefficients of $\lambda$ are linearly independent, and that 
$I$ is a saturated radical ideal. (We verified this in Macaulay-2.) Thus the assertions 
remain true for a general $F$. \qed 

\smallskip 

We have identities 
\begin{equation} 
(\lambda,Q)_2 = (\lambda,(Q^2,F)_4)_2 = 0, 
\label{d6.id1} \end{equation} 
and using the techniques of \cite[Ch.~VI]{GY} one verifies that each syzygy satisfied by 
$\lambda$ is a polynomial in these two. Hence $I$ has two first syzygies, and 
we have a minimal resolution 
\[ 0 \la R/I \la R \la R(-3)^3 \la R(-4) \oplus R(-5) \la 0. \] 
This allows us to calculate the Hilbert function of $R/I$, which shows that 
$\deg \vartheta(F)=7$. 

\smallskip 

The pairs of examples in \S\ref{section.d4.conic}, \ref{section.d6.cubic} 
and \S\ref{section.d4.3p}, \ref{section.d6.7p} suggest that there is an infinite 
class of results of each type. It may be noticed that all the considerations in 
\S\ref{section.varieties} make use of the sign sequence, whereas 
those in \S\ref{section.vartheta.locus} make use of the involutor. The connection between the 
two is mediated by a rather complicated formula, and each of them seems to contain 
a kind of algebraic information on its surface which is not easy to extract from the other. 

We are at a boundary which marks a sharp change in the texture of this paper. 
Henceforth all the calculations will make very substantial use of the graphical calculus. 
\section{Coefficients of the geometric involutor} 
\label{proof.geometric.formula} 
In this section we prove formula~(\ref{geom.involution}) from \S\ref{system.sigmad}, 
which gives a transvectant series expression for the geometric involution. 
The reader is referred to~\cite[\S2]{spinnet} for an explanation of the notation used.

\subsection{} Let $\mathcal{M}_d$ be the space of maps
\[
\begin{array}{rl}
\psi: & S_d\times S_2\lra S_d \\
 & (F,Q)\longmapsto \psi(F,Q)
\end{array}
\]
which are linear in $F$ and homogeneous of degree $d$ in $Q$.
Let $\mathcal{M}_d^{SL_2}$ be the subspace of equivariant maps.

\begin{Lemma}\label{first.dim.lemma}
${\rm dim}\ \mathcal{M}_d^{SL_2}=n+1$.
\end{Lemma}
\demo
Successively using self-duality, Hermite reciprocity and the Clebsch-Gordan
decomposition, we have 
\begin{align*}
\mathcal{M}_d & = {\rm Hom}\left(
S_d(S_2^\vee)\otimes S_d, S_d\right) \\
 & = S_d(S_2)\otimes S_d\otimes S_d \\
 & = S_2(S_d)\otimes S_d\otimes S_d \\
 & = \left(\bigoplus_{i=0}^n S_{2d-4i}\right)\otimes S_d\otimes S_d\ 
\end{align*}
as $SL_2$-representations.
Therefore,
\[
\mathcal{M}_d^{SL_2}= \bigoplus_{i=0}^n \left(S_{2d-4i}\otimes S_d\otimes S_d\right)^{SL_2}
\]
where each summand has dimension $1$.
\qed

\subsection{} 
Let $x = \{x_1,x_2\},y = \{y_1,y_2\}$ be two sets of point coordinates, and $Q$ a binary 
quadratic. We will use the customary notation of symbolic calculus, where 
\[ 
y \, \partial_x = y_1 \frac{\partial}{\partial x_1} + y_2 \frac{\partial}{\partial x_2}, 
\quad \text{and} \quad (x \, y) = x_1 y_2 - x_2 y_1. \] 
Let $\mathcal{C}$ denote the space of joint covariants $C(Q,x,y)$ 
which are homogeneous of degree $d$ in each of $Q,x,y$. 
\begin{Lemma} We have an isomorphism 
$\mathcal{M}_d^{SL_2}\simeq \mathcal{C}$.
\end{Lemma}
\demo
There is a linear map $\mathcal{L}:C\mapsto\psi$ given by
\[ \psi(F,Q)=(C,F(y))_d \]
where the transvection is with respect to the $y$-variables.
This map has an inverse given by the so-called evectant in the $y$-variables.
Namely, if $F$ is written as $(a_0,\ldots,a_d\cbrac x_1,x_2)^d$, then
one recovers $C$ from $\psi$ by letting 
\[
C(Q,x,y)=\mathcal{E}_{y\ra F}\ \psi(F,Q), 
\]
where
\[ \mathcal{E}_{y\ra F}=\sum_{i=0}^d \, \binom{d}{i} 
y_1^{i} (-y_2)^{d-i} \frac{\partial}{\partial a_i}. 
\]
\qed
\subsection{} 
We now define a particular collection of elements in $\mathcal{C}$. 
Let $Q_{xy}=\frac{1}{2}(y \, \partial_x) \, Q(x)$ denote the first polar in $y$. Define 
\[ 
\mathfrak{p}_i=\Delta_Q^i \, (x \, y)^{2i} \, Q_{xy}^{d-2i} \, \quad 
(0\le i\le n). \] 
\begin{Lemma}
The $\{ \mathfrak{p}_i: 0\le i\le n\}$  form a linear basis of $\mathcal{C}$.
\end{Lemma}
\demo
On account of the previous lemmas, it is enough to check linear independence.
Suppose one has a linear relation $\sum_{i=0}^{n} \, \lambda_i \, \mathfrak{p}_i=0$.
Then specialise to $Q(x)=x_1 \, x_2$, and also to $x_2=y_1=y_2=1$.
This gives $\Delta_Q=1$, $(x \, y)=x_1-1$ and
\[
Q_{xy}=\frac{1}{2}(x_1 \, y_2+x_2 \, y_1)=\frac{1}{2}(x_1+1)\ .
\] 
The linear relation therefore becomes
\[ \sum_{i=0}^n \, \lambda_i \, (x_1-1)^{2i} \, \left(\frac{x_1+1}{2}\right)^{d-2i}=0
\]
identically in $x_1$, which easily implies that all the $\lambda$'s vanish. \qed

\smallskip 

Note that $\mathfrak{p}_i=\Delta_Q^i\times \widehat{\mathfrak{p}}_i$
where, using the graphical notation of~\cite{spinnet},
\[
\widehat{\mathfrak{p}}_i=(x \, y)^{2i} Q_{xy}^{d-2i}=
\parbox{3cm}{\psfrag{m}{$\scriptstyle{d}$}
\psfrag{n}{$\scriptstyle{d-2i}$}\psfrag{k}{$\scriptscriptstyle{2i}$}
\psfrag{q}{$\scriptstyle{Q}$}\psfrag{x}{$\scriptstyle{x}$}
\psfrag{y}{$\scriptstyle{y}$}
\includegraphics[width=3cm]{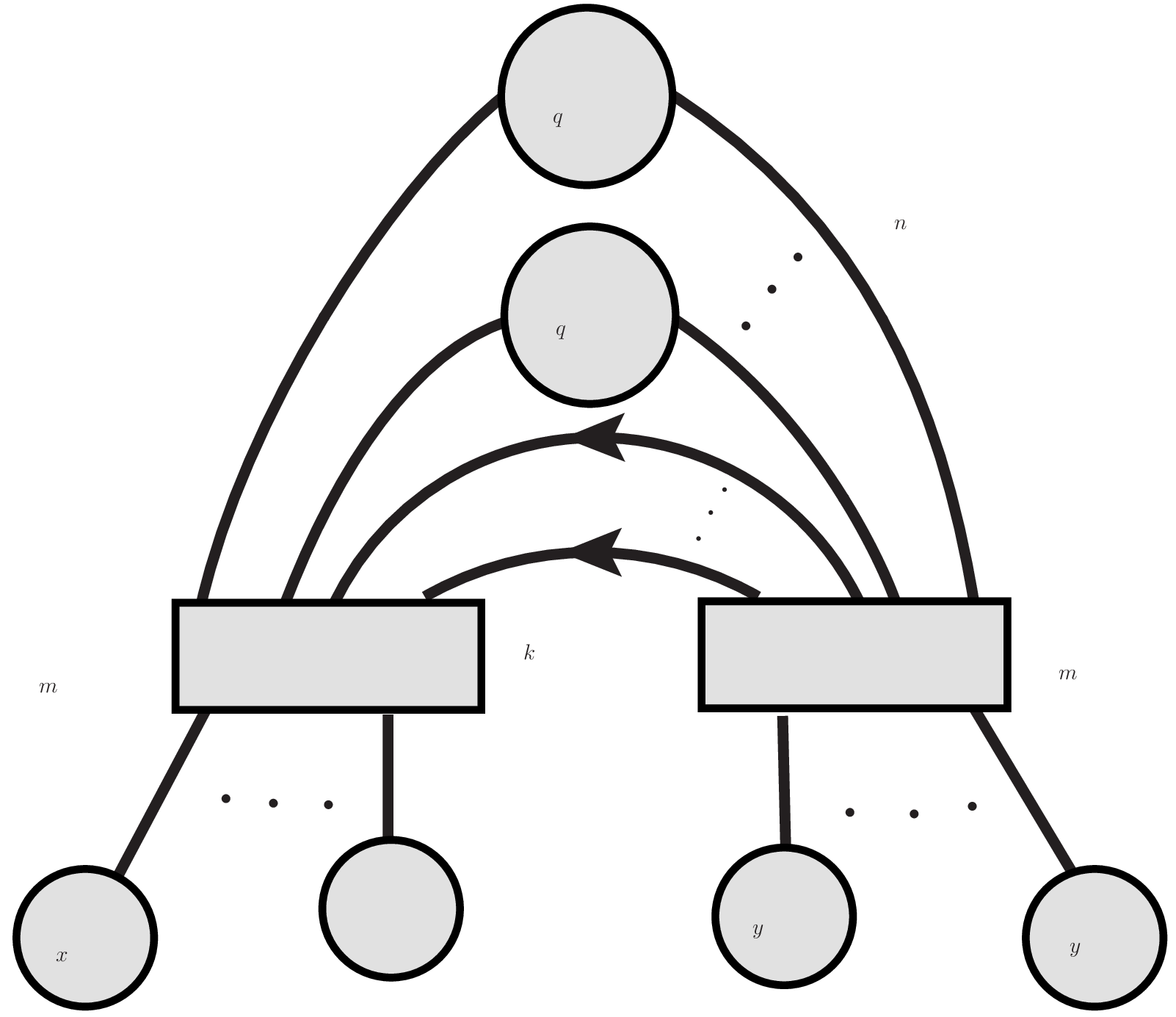}}\ \ .
\]
The graphical notation also allows for an easy visualization of $\mathcal{L}(\mathfrak{p}_i)$,
namely as the map
\begin{equation}
(F,Q)\longmapsto
\Delta_Q^i \times\ 
\parbox{3cm}{\psfrag{f}{$\scriptstyle{F}$}
\psfrag{m}{$\scriptstyle{d-2i}$}\psfrag{n}{$\scriptstyle{2i}$}
\psfrag{q}{$\scriptstyle{Q}$}\psfrag{x}{$\scriptstyle{x}$}
\includegraphics[width=3cm]{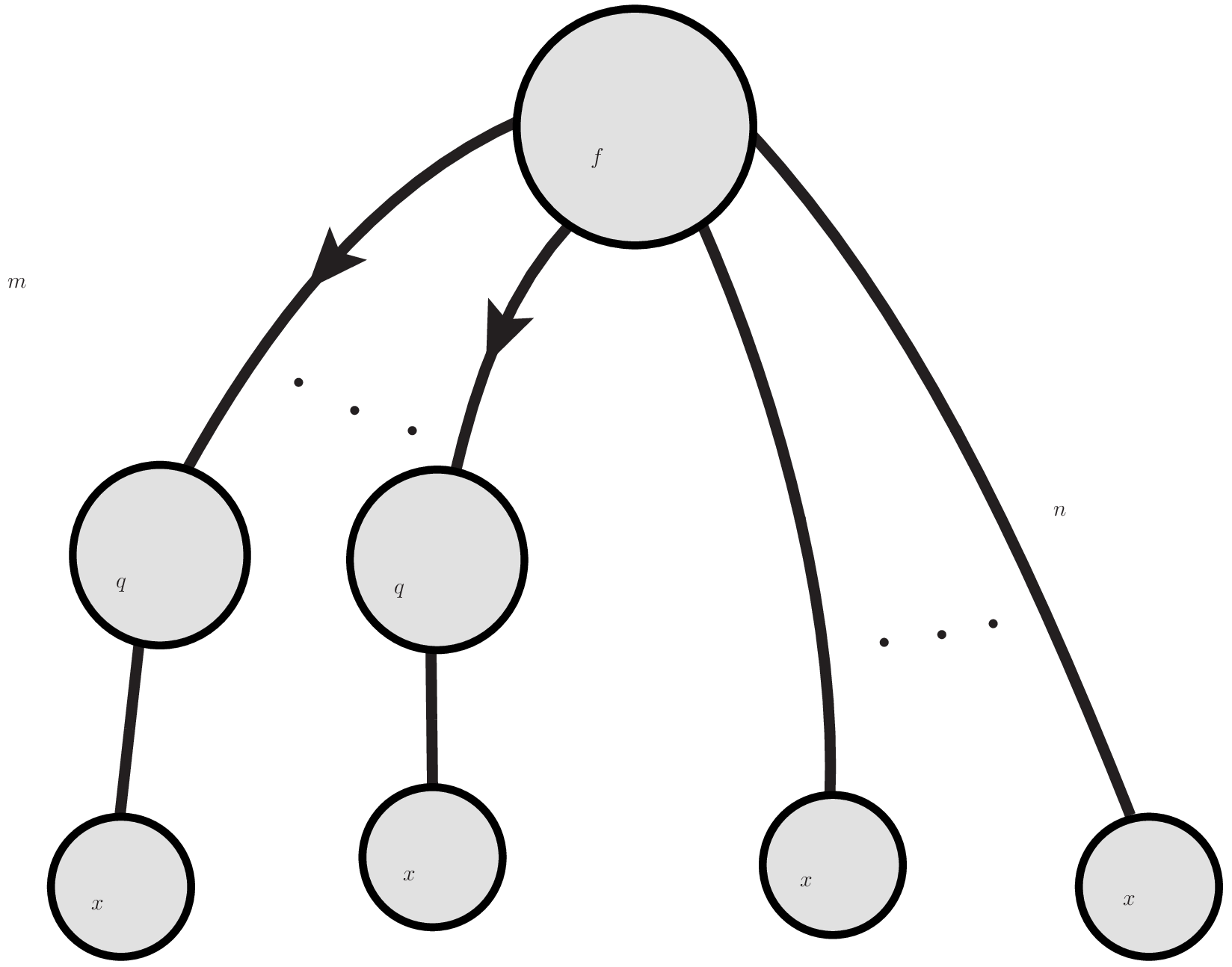}}\ \ .
\label{P.graphical.def}
\end{equation}
Note that the $Q$'s are mounted in parallel, hence our choice of letter
`$\mathfrak{p}$' for the notation.
Since
\[
\parbox{0.8cm}{\psfrag{f}{$\scriptstyle{F}$}
\includegraphics[width=0.8cm]{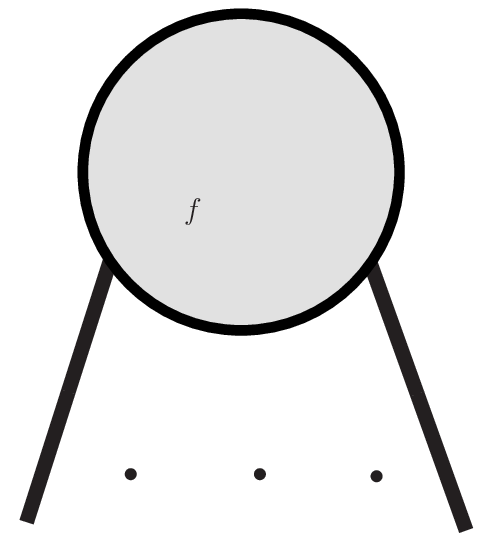}}
=\parbox{1.8cm}{\psfrag{1}{$\scriptscriptstyle{l_1}$}
\psfrag{2}{$\scriptscriptstyle{l_d}$}
\includegraphics[width=0.8cm]{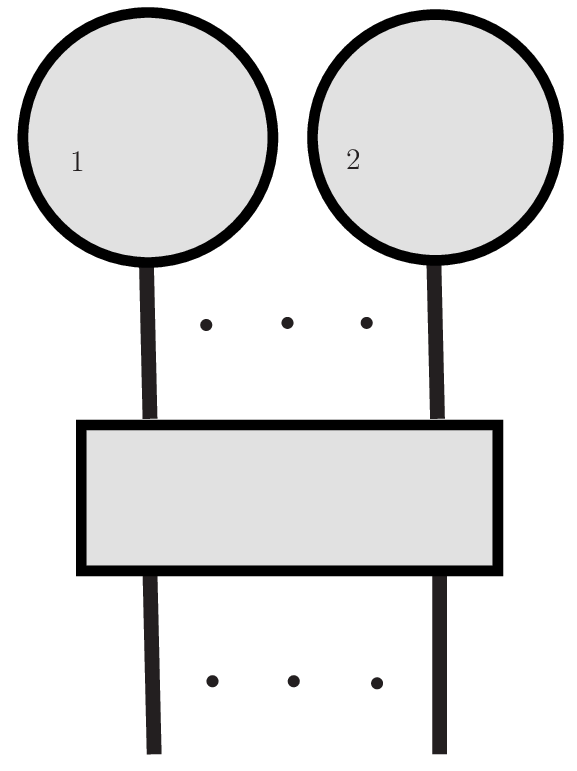}}, 
\]
for the particular case of the geometric involutor we have
\begin{equation}\label{sigmavsU}
\sigma_Q(F)=2^d \, \mathcal{L}(\mathfrak{p}_0)(F,Q). 
\end{equation}

\subsection{} 
Now define a new collection $\{ \mathfrak{t}_i: 0\le i\le n\}$ in $\mathcal{C}$
given by the $\mathcal{L}$ inverse images of the maps
\[
(F,Q)\longmapsto \Delta_Q^i (Q^{d-2i},F)_{d-2i}\ .
\]
The letter `$\mathfrak{t}$' stands for transvection. Since
\[
(Q^{d-2i},F)_{d-2i}= 
\parbox{3.8cm}{\psfrag{1}{$\scriptstyle{d-2i}$}
\psfrag{2}{$\scriptstyle{2i}$}
\psfrag{f}{$\scriptstyle{F}$}\psfrag{x}{$\scriptstyle{x}$}
\psfrag{q}{$\scriptstyle{Q}$}
\includegraphics[width=3.8cm]{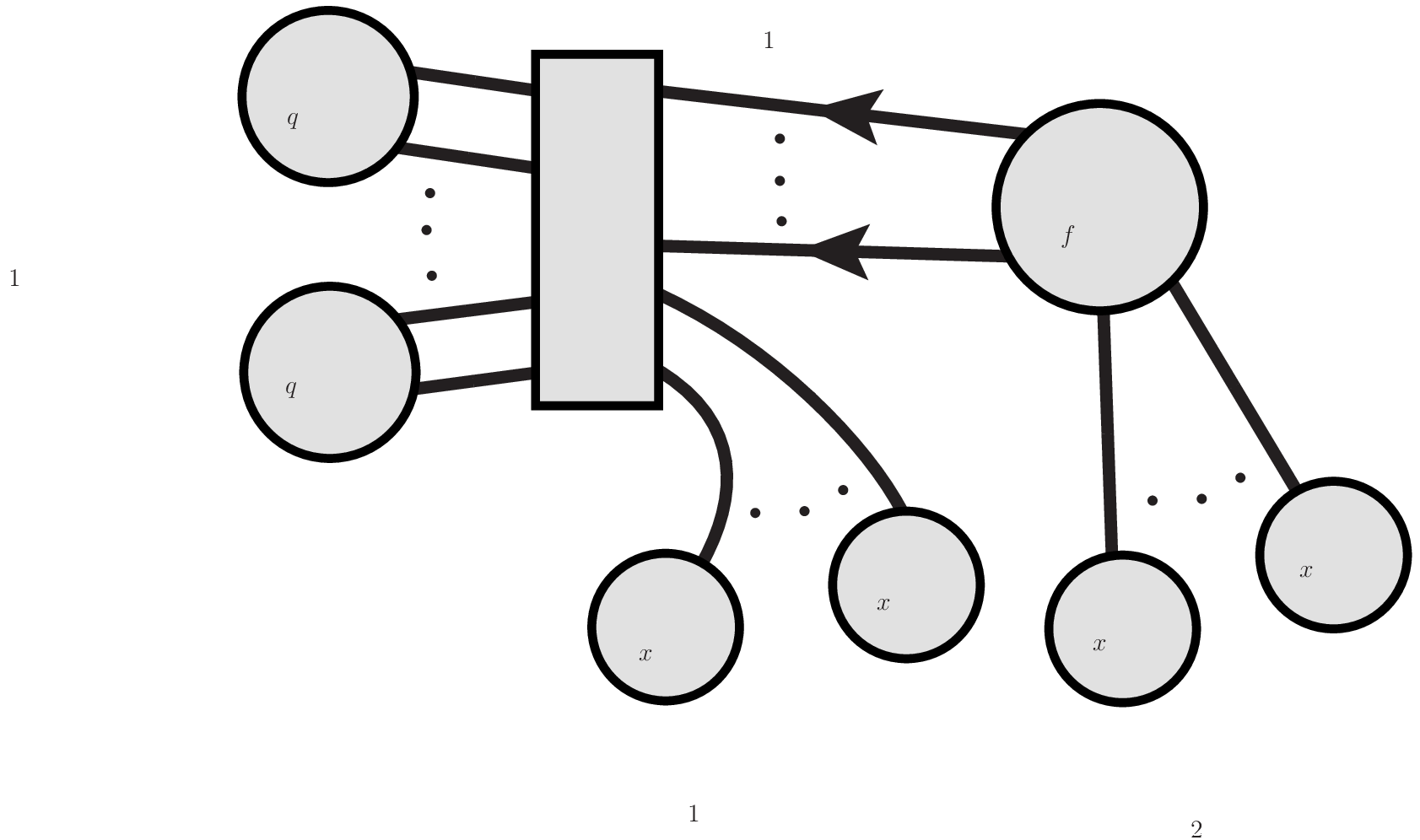}}
\]
and since the evectant amounts to the graphical substitution
\[
\parbox{0.8cm}{\psfrag{f}{$\scriptstyle{F}$}
\includegraphics[width=0.8cm]{Sec6fig3.eps}}
\lra
\parbox{1.8cm}{\psfrag{y}{$\scriptstyle{y}$}
\includegraphics[width=1cm]{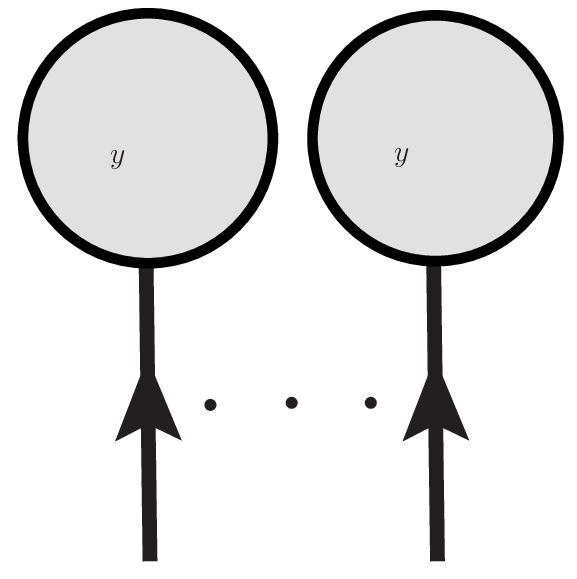}}\ \ ,
\]
we have
\[
\mathfrak{t}_i=\Delta_Q^i (x \, y)^{2i} \ 
\parbox{2.3cm}{\psfrag{q}{$\scriptstyle{Q}$}
\psfrag{x}{$\scriptstyle{x}$}\psfrag{y}{$\scriptstyle{y}$}
\psfrag{1}{$\scriptstyle{d-2i}$}
\includegraphics[width=2.3cm]{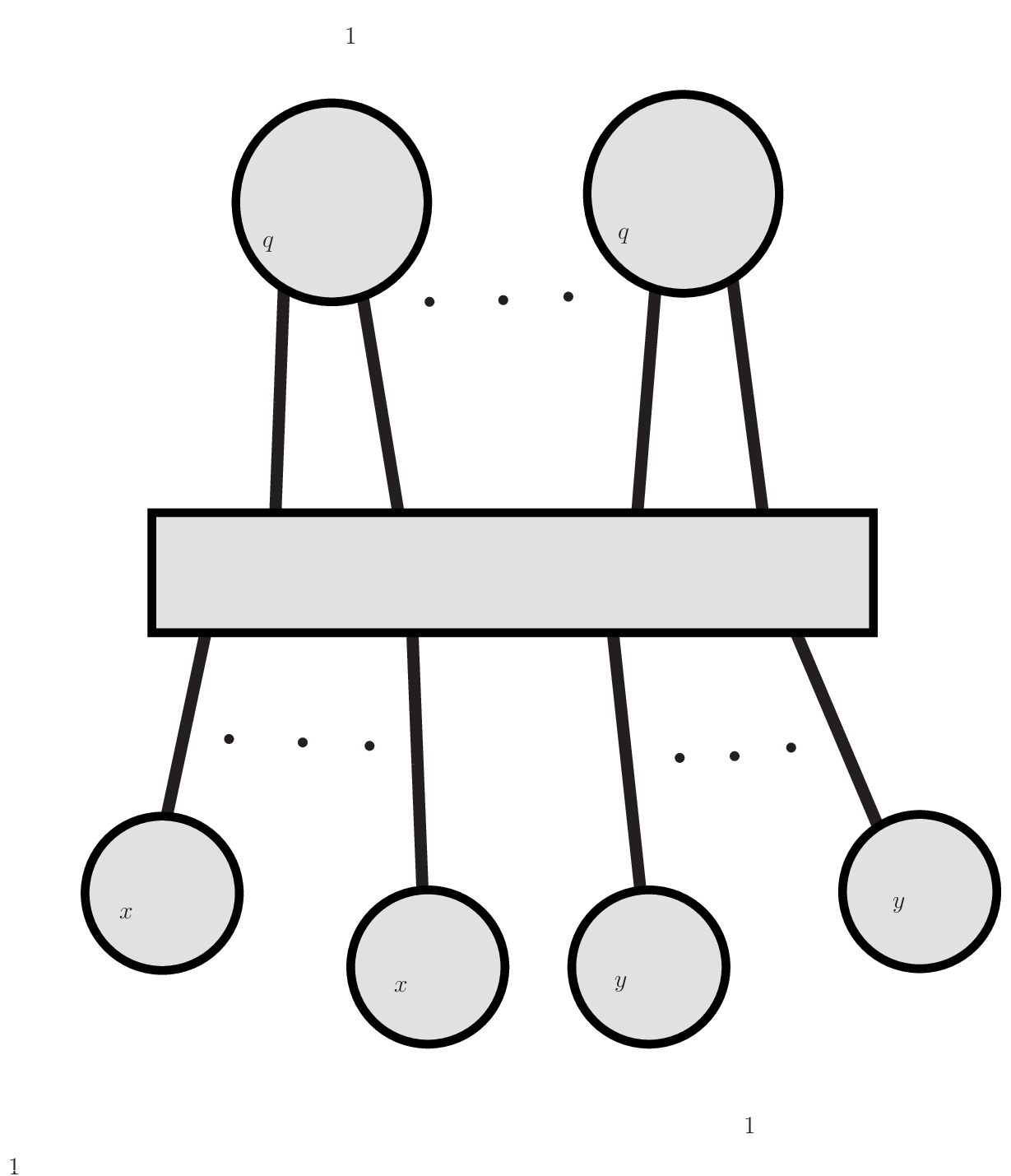}}\ \ .
\]
The following proposition gives the transition matrix for the two bases. 
\begin{Proposition}\label{triangprop}
For any $i$ in the range $0\le i\le n$, we have 
\[
\mathfrak{p}_i=\sum_{j=i}^{n} \, G_{i,j} \, \mathfrak{t}_j \, , 
\]
where
\[
G_{i,j}=\frac{(d-2i)!\ (2d-4j+1)!\ (d-i-j)!}{4^{j-i}\ (d-2j)!^2\ (2d-2i-2j+1)!\ (j-i)!}\ .
\]
Therefore, $\{\mathfrak{t}_i: 0\le i\le n\}$ is also a linear basis for $\mathcal{C}$.
\end{Proposition}
\demo
By the idempotence of symmetrizers, one can write
\[
\widehat{\mathfrak{p}}_i=
\parbox{2.6cm}{\psfrag{q}{$\scriptstyle{Q}$}
\psfrag{x}{$\scriptstyle{x}$}\psfrag{y}{$\scriptstyle{y}$}
\psfrag{n}{$\scriptstyle{d-2i}$}\psfrag{k}{$\scriptstyle{2i}$}
\psfrag{m}{$\scriptstyle{d}$}
\includegraphics[width=2.6cm]{Sec6fig8.eps}}
\]
\vskip 1cm
\noindent
where the dotted lines indicate where we next use the Clebsch-Gordan identity~\cite[Eq. 12]{spinnet}.
The latter gives
\[
\widehat{\mathfrak{p}}_i=\sum_{k=0}^d
\frac{\left(\begin{array}{c} d \\ k \end{array}\right)^2}{\left(
\begin{array}{c} 2d-k+1 \\ k \end{array}\right)}
\parbox{2.6cm}{\psfrag{q}{$\scriptstyle{Q}$}
\psfrag{x}{$\scriptstyle{x}$}\psfrag{y}{$\scriptstyle{y}$}
\psfrag{k}{$\scriptstyle{k}$}\psfrag{i}{$\scriptstyle{2i}$}
\psfrag{m}{$\scriptstyle{2d-2k}$}\psfrag{d}{$\scriptstyle{d}$}
\psfrag{n}{$\scriptstyle{d-2i}$}
\includegraphics[width=2.6cm]{Sec6fig9.eps}}
\]
\vskip 1cm
\[
=\sum_{k=0}^d
\frac{d!^2(2d-2k+1)!}{k! (d-k)!^2 (2d-k+1)!}
\ (x \, y)^k\times
\parbox{2.6cm}{\psfrag{q}{$\scriptstyle{Q}$}
\psfrag{x}{$\scriptstyle{x}$}\psfrag{y}{$\scriptstyle{y}$}
\psfrag{k}{$\scriptstyle{k}$}\psfrag{i}{$\scriptstyle{2i}$}
\psfrag{m}{$\scriptstyle{2d-2k}$}\psfrag{d}{$\scriptstyle{d}$}
\psfrag{n}{$\scriptstyle{d-2i}$}\psfrag{p}{$\scriptstyle{d-k}$}
\includegraphics[width=2.6cm]{Sec6fig10.eps}}
\]
\vskip 1cm
The last graphical expression is the $(d-k)$-th polar in $y$ of
\[
R_{i,k}(x)=
\parbox{2.6cm}{\psfrag{q}{$\scriptstyle{Q}$}
\psfrag{x}{$\scriptstyle{x}$}
\psfrag{k}{$\scriptstyle{k}$}\psfrag{i}{$\scriptstyle{2i}$}
\psfrag{d}{$\scriptstyle{d}$}
\psfrag{n}{$\scriptstyle{d-2i}$}\psfrag{p}{$\scriptstyle{d-k}$}
\includegraphics[width=2.6cm]{Sec6fig11.eps}}\ \ ,
\]
namely, the result of applying $\frac{(d-k)!}{(2d-2k)!}(y\partial_x)^{d-k}$ to $R_{i,k}(x)$.
Now,
\[
R_{i,k}(x)=
\parbox{3.3cm}{\raisebox{7ex}{\psfrag{q}{$\scriptstyle{Q}$}
\psfrag{x}{$\scriptstyle{x}$}
\psfrag{k}{$\scriptstyle{k}$}\psfrag{i}{$\scriptstyle{2i}$}
\psfrag{d}{$\scriptstyle{d}$}
\psfrag{n}{$\scriptstyle{d-2i}$}\psfrag{p}{$\scriptstyle{d-k}$}
\includegraphics[width=3.3cm]{Sec6fig12.eps}}}
\]
where the dotted double arrow again indicates where~\cite[Eq. 12]{spinnet} will be used next.
Indeed,
\[
R_{i,k}(x)=\sum_{p=0}^{d-2i}
\frac{\left(\begin{array}{c} d-2i \\ p \end{array}\right)^2}{\left(
\begin{array}{c} 2d-4i-p+1 \\ p \end{array}\right)}
\parbox{4.2cm}{\raisebox{12ex}{\psfrag{q}{$\scriptstyle{Q}$}
\psfrag{x}{$\scriptstyle{x}$}
\psfrag{k}{$\scriptstyle{k}$}\psfrag{i}{$\scriptstyle{2i}$}
\psfrag{d}{$\scriptstyle{d}$}
\psfrag{n}{$\scriptstyle{d-2i}$}\psfrag{p}{$\scriptstyle{d-k}$}
\psfrag{1}{$\scriptstyle{p}$}\psfrag{2}{$\scriptstyle{2d-4i-2p}$}
\includegraphics[width=4.2cm]{Sec6fig13.eps}}}
\]
\noindent
by the Clebsch-Gordan identity.
Letting the $d$-symmetrizers absorb the bottom
$(d-2i)$-symmetrizers as shown by the dotted arrows, we obtain
\[
R_{i,k}(x)=
\sum_{p=0}^{d-2i}
\frac{(d-2i)!^2(2d-4i-2p+1)!}{p!(d-2i-p)!^2(2d-4i-p+1)!}
\parbox{4.2cm}{\raisebox{14ex}{\psfrag{q}{$\scriptstyle{Q}$}
\psfrag{x}{$\scriptstyle{x}$}
\psfrag{k}{$\scriptstyle{k}$}\psfrag{1}{$\scriptstyle{p}$}
\psfrag{d}{$\scriptstyle{d}$}
\psfrag{n}{$\scriptstyle{d-2i}$}\psfrag{3}{$\scriptstyle{2i+p}$}
\psfrag{2}{$\scriptstyle{2d-4i-2p}$}\psfrag{m}{$\scriptstyle{2d-2k}$}
\includegraphics[width=4.2cm]{Sec6fig14.eps}}}
\]
\noindent
where we put an extra symmetrizer at the bottom, which is allowed since the $x$ `blobs' are identical.
Now replace the portion indicated by the dotted double arrow, using~\cite[Eq. 13]{spinnet}.
The notation $\bbone\{C\}$ will be used for the characteristic function; i.e., given a logical 
condition $C$ on a parameter, $\bbone\{C\}$ is $1$ if $C$ is satisfied and $0$ otherwise. 
The result is 
\[
R_{i,k}(x)=
\sum_{p=0}^{d-2i}
\frac{(d-2i)!^2(2d-4i-2p+1)!}{p!(d-2i-p)!^2(2d-4i-p+1)!}
\times\bbone\{k=2i+p\}
\]
\[
\times 
\frac{(2i+p)!(2d-2i-p+1)!(d-2i-p)!^2}{d!^2\ (2d-4i-2p+1)!}
\times
\parbox{3.4cm}{\psfrag{q}{$\scriptstyle{Q}$}
\psfrag{x}{$\scriptstyle{x}$}
\psfrag{2}{$\scriptstyle{p}$}
\psfrag{1}{$\scriptstyle{d-2i-p}$}\psfrag{n}{$\scriptstyle{d-2i}$}
\includegraphics[width=3.4cm]{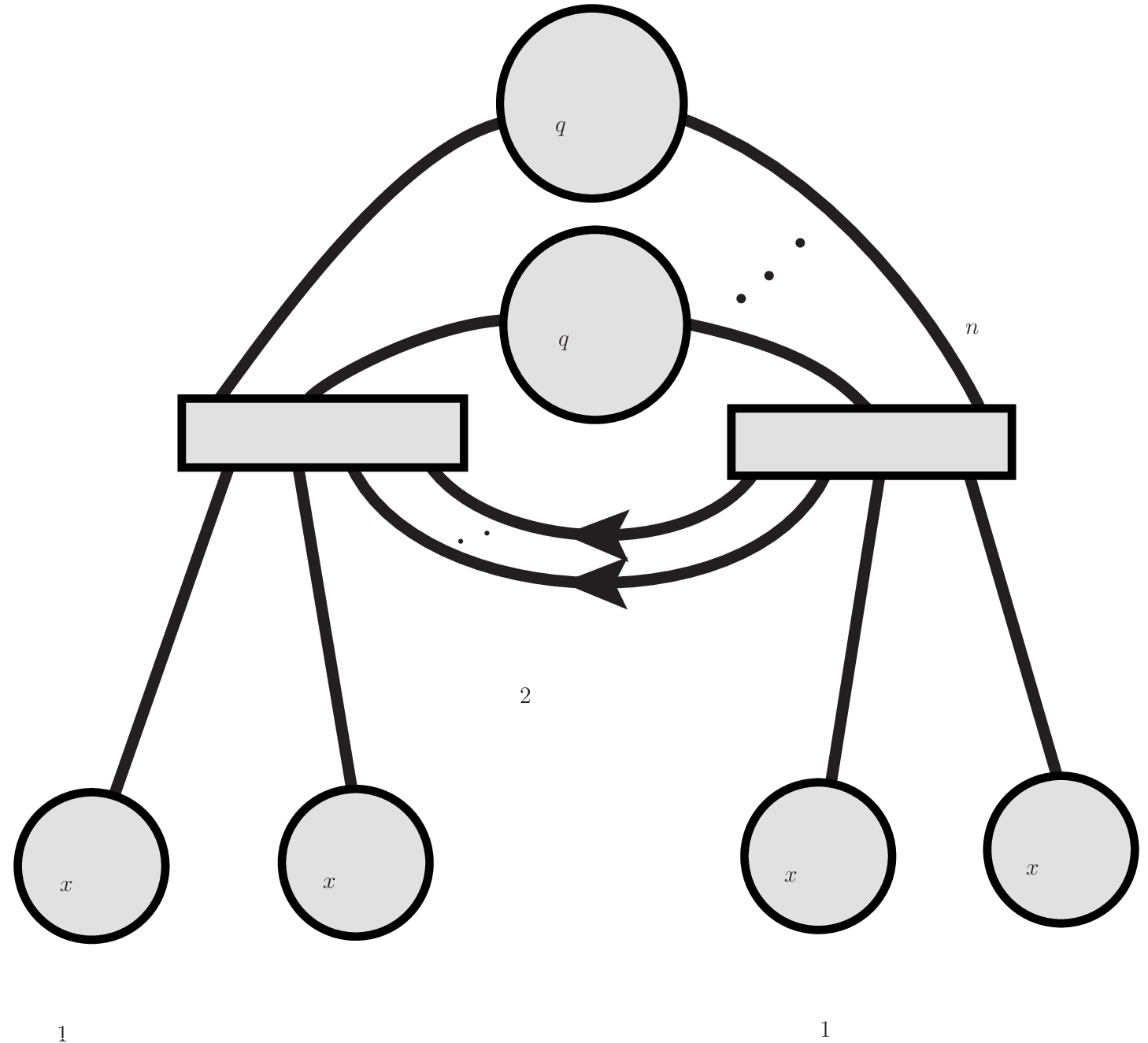}}\ \ .
\]
After simplifying, 
\[
R_{i,k}(x)=\bbone\left\{\begin{array}{c} k-2i\le d-2i \\ k-2i\ge 0 \end{array}\right\}
\frac{(d-2i)!^2 k!(2d-k+1)!}{(k-2i)!(2d-2i-k+1)! d!^2}\times
S_{i,k}
\]
where
\[
S_{i,k}=
\parbox{3.4cm}{\psfrag{q}{$\scriptstyle{Q}$}
\psfrag{x}{$\scriptstyle{x}$}
\psfrag{2}{$\scriptstyle{k-2i}$}
\psfrag{1}{$\scriptstyle{d-k}$}\psfrag{n}{$\scriptstyle{d-2i}$}
\includegraphics[width=3.4cm]{Sec6fig15.eps}}\ \ .
\]
The latter is a covariant of the quadratic $Q$ which is of degree $d-2i$ and order $2d-2k$.
From the known structure of the ring of covariants for a quadratic binary 
form (see~\cite[\S85]{GY}), we see that $S_{i,k}$ must be of the form
\[ S_{i,k}=\bbone\{k\ {\rm even}\} \, \alpha_{i,k}\ \Delta_Q^{\frac{k-2i}{2}}\ Q^{d-k}, 
\]
for some scalar $\alpha_{i,k}$. Let $k=2j$ for some new index $j$.
In order to determine $\alpha_{i,k}$, specialise to $Q=x_1^2+x_2^2$. 
The `blob' of $Q$ becomes
\[
\parbox{0.9cm}{\psfrag{1}{$\scriptstyle{1}$}
\includegraphics[width=0.9cm]{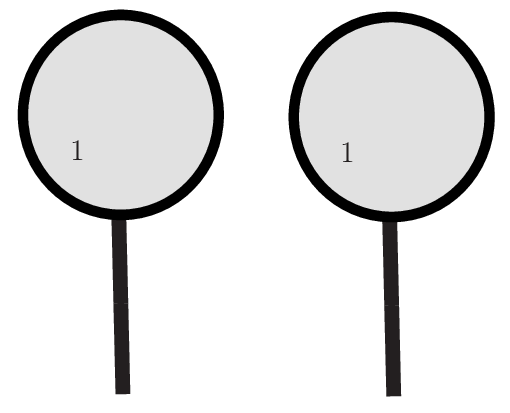}}
+
\parbox{0.9cm}{\psfrag{1}{$\scriptstyle{2}$}
\includegraphics[width=0.9cm]{Sec6fig16.eps}}
\]
where $\parbox{0.7cm}{\psfrag{1}{$\scriptstyle{1}$}
\includegraphics[width=0.7cm]{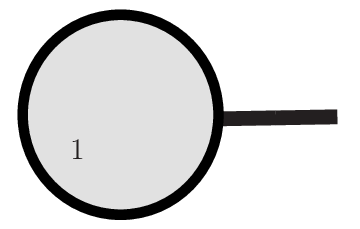}}$
and $\parbox{0.7cm}{\psfrag{1}{$\scriptstyle{2}$}
\includegraphics[width=0.7cm]{Sec6fig17.eps}}$
are the graphical representations of the canonical basis
vectors. By expanding the sums which give each of the $d-2i$ copies of $Q$, we get 
\[
S_{i,k}=\sum_{r=0}^{d-2i}
\left(\begin{array}{c} d-2i \\ r \end{array}\right)
\parbox{2.9cm}{\psfrag{1}{$\scriptstyle{d-2j}$}
\psfrag{2}{$\scriptstyle{2j-2i}$}\psfrag{a}{$\scriptstyle{1}$}
\psfrag{b}{$\scriptstyle{2}$}\psfrag{x}{$\scriptstyle{x}$}
\psfrag{r}{$\scriptstyle{r}$}\psfrag{d}{$\scriptstyle{d-2i-r}$}
\includegraphics[width=2.9cm]{Sec6fig18.eps}}\qquad\qquad .
\]

\bigskip \bigskip 

We then further specialise to $x_1=1$ and $x_2=0$, which gives
\[
S_{i,k}=\sum_{r=0}^{d-2i}
\left(\begin{array}{c} d-2i \\ r \end{array}\right)
\parbox{3.2cm}{\psfrag{1}{$\scriptstyle{d-2j}$}
\psfrag{2}{$\scriptstyle{2j-2i}$}\psfrag{a}{$\scriptstyle{1}$}
\psfrag{b}{$\scriptstyle{2}$}
\psfrag{r}{$\scriptstyle{r}$}\psfrag{d}{$\scriptstyle{d-2i-r}$}
\includegraphics[width=3.2cm]{Sec6fig19.eps}}\qquad\qquad\qquad .
\]
In order to compute the last graphical expression, note the following identity
\[
\parbox{2.5cm}{\psfrag{1}{$\scriptstyle{d-2j}$}
\psfrag{2}{$\scriptstyle{2j-2i}$}\psfrag{a}{$\scriptstyle{1}$}
\psfrag{b}{$\scriptstyle{2}$}
\psfrag{r}{$\scriptstyle{r}$}\psfrag{d}{$\scriptstyle{d-2i-r}$}
\includegraphics[width=2.5cm]{Sec6fig20.eps}}
=
\parbox{2.5cm}{\psfrag{1}{$\scriptstyle{d-2j}$}
\psfrag{2}{$\scriptstyle{2j-2i}$}\psfrag{a}{$\scriptstyle{1}$}
\psfrag{b}{$\scriptstyle{2}$}
\psfrag{r}{$\scriptstyle{r}$}\psfrag{d}{$\scriptstyle{d-2i-r}$}
\includegraphics[width=2.5cm]{Sec6fig21.eps}}
\]
\vskip 0.6cm
\[
=\bbone\{r\ge d-2j\}\ \frac{r!\ (2j-2i)!}{(d-2i)!(r-d+2j)!}
\ \times\ 
\parbox{3cm}{\psfrag{1}{$\scriptstyle{d-2j}$}
\psfrag{2}{$\scriptstyle{2j-2i}$}\psfrag{a}{$\scriptstyle{1}$}
\psfrag{b}{$\scriptstyle{2}$}
\psfrag{r}{$\scriptstyle{r-d+2j}$}\psfrag{d}{$\scriptstyle{d-2i-r}$}
\includegraphics[width=3cm]{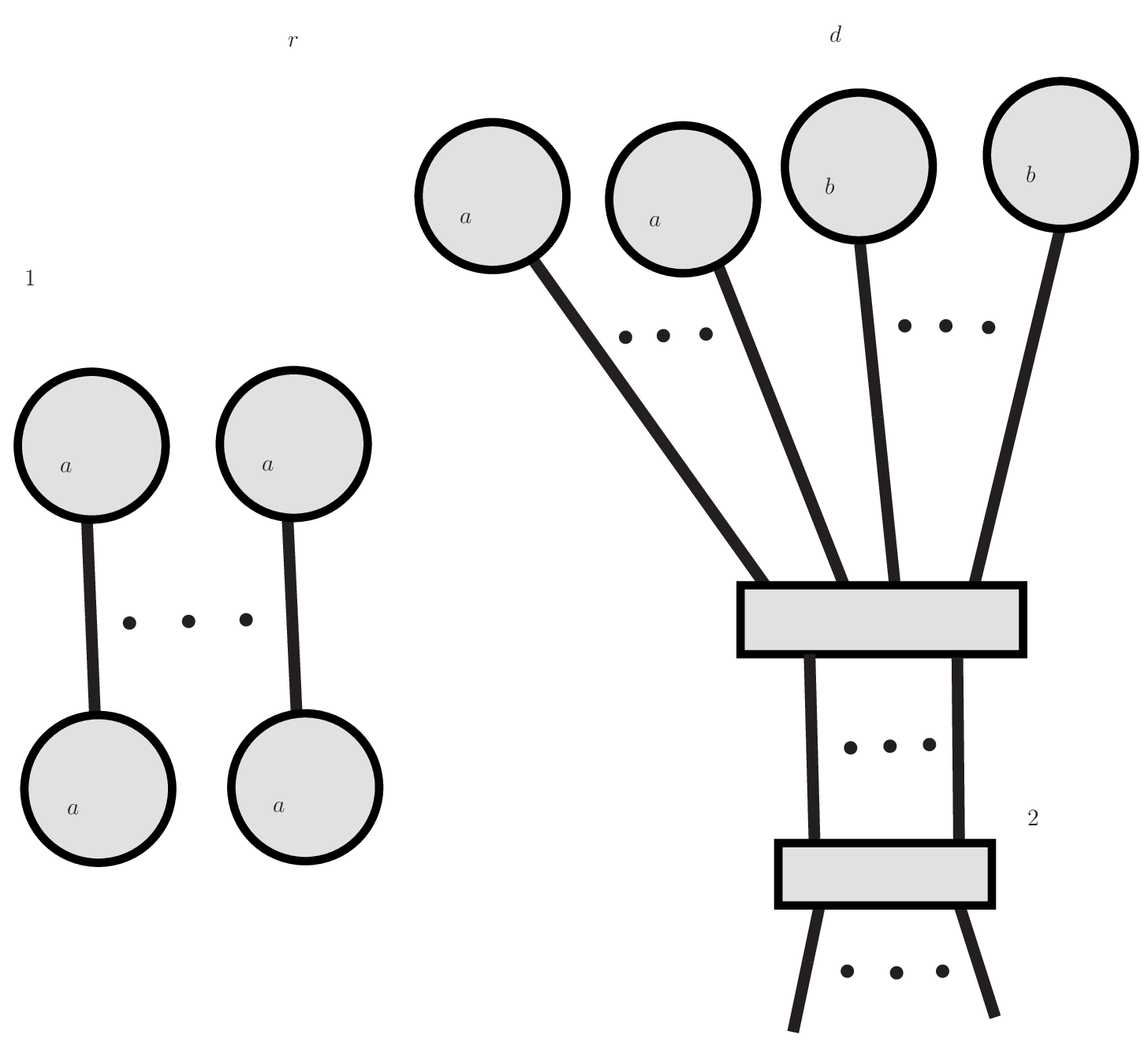}}\qquad\ .
\]
Indeed, one can expand the top symmetrizer as a sum $\frac{1}{(d-2i)!}\sum_{\tau}$
over permutations $\tau$ of $d-2i$ elements. Only those permutations which connect the bottom
$\parbox{0.7cm}{\psfrag{1}{$\scriptstyle{1}$}
\includegraphics[width=0.7cm]{Sec6fig17.eps}}$
to the top ones survive.  There are $\frac{r!}{(r-d+2j)!}\times (2j-2i)!$ such permutations.  
Because of the newly introduced bottom symmetriser, all of them give the same contribution. 
As a result, we have
\[ \begin{aligned} 
S_{i,k}=\sum_{r=0}^{d-2i}
\left(\begin{array}{c} d-2i \\ r \end{array}\right)
\times & \bbone\{r\ge d-2j\}\left(\frac{r!\ (2j-2i)!}{(d-2i)!(r-d+2j)!}\right)^2 \\ 
\times & 
\parbox{5cm}{\psfrag{1}{$\scriptstyle{d-2j}$}
\psfrag{2}{$\scriptstyle{2j-2i}$}\psfrag{a}{$\scriptstyle{1}$}
\psfrag{b}{$\scriptstyle{2}$}
\psfrag{r}{$\scriptstyle{r-d+2j}$}\psfrag{d}{$\scriptstyle{d-2i-r}$}
\includegraphics[width=5cm]{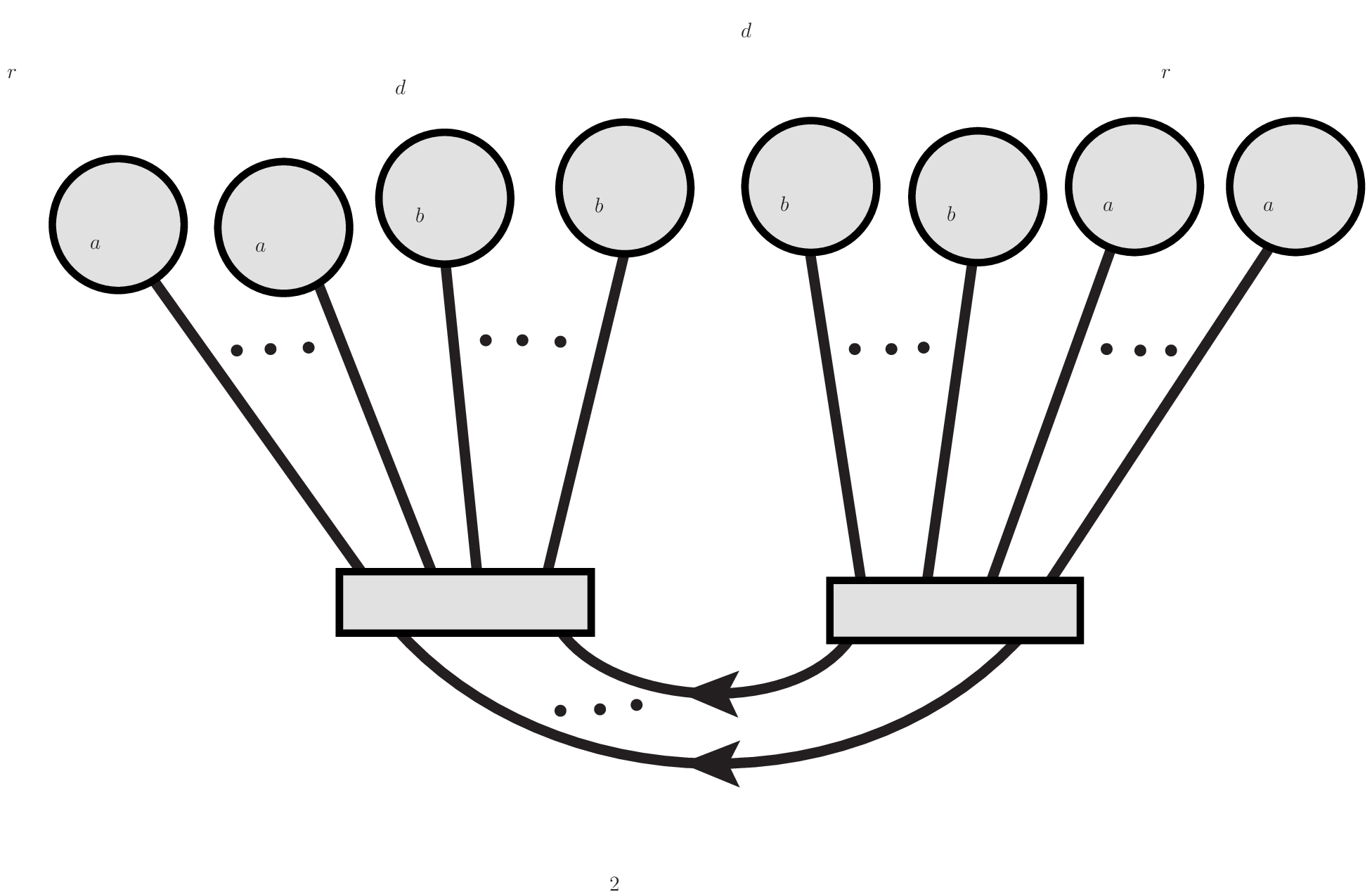}}
\end{aligned} \] 
which is equal to 
\[ \begin{aligned} 
\sum_{r} \, & \bbone\{d-2j\le r\le d-2i\}
\ \frac{r!\ (2j-2i)!^2}{(d-2i)!(d-2i-r)!(r-d+2j)!^2} \\ 
& \times (-1)^{r-d+2j}\times 
\parbox{3cm}{\raisebox{6ex}{\psfrag{a}{$\scriptstyle{1}$}
\psfrag{b}{$\scriptstyle{2}$}
\psfrag{r}{$\scriptstyle{r-d+2j}$}\psfrag{d}{$\scriptstyle{d-2i-r}$}
\includegraphics[width=3cm]{Sec6fig24.eps}}}. 
\end{aligned} \] 
After evaluating the last graphical expression, this is seen to be equal to  
\[ \begin{aligned} 
=\sum_{r} \; & \bbone\{d-2j\le r\le d-2i\}
\ \frac{r!\ (2j-2i)!^2}{(d-2i)!(d-2i-r)!(r-d+2j)!^2} \\ 
\times & \bbone\{d-2i-r=r-d+2j\}\ (-1)^{r-d}\ \frac{(d-2i-r)!^2}{(2j-2i)!}. 
\end{aligned} \]
Solving for $r=d-i-j$ and tidying up the final result, we get
\[
S_{i,k}=(-1)^{i+j}
\frac{(d-i-j)!\ (2j-2i)!}{(d-2i)!\ (j-i)!}\ .
\]
On the other hand, $\Delta_Q=-4$ and $Q(x)=1$ for the assumed specialisation, 
so that $S_{i,k}=(-4)^{j-i} \, \alpha_{i,k}$.
Thus, we have identified the proportionality constant:
\[
\alpha_{i,k}=\frac{(d-i-j)!\ (2j-2i)!}{4^{j-i}\ (d-2i)!\ (j-i)!}\ .
\]
Putting everything together, we obtain
\[ \begin{aligned} 
\widehat{\mathfrak{p}}_i= & 
\sum_{j=0}^{n}\ \frac{d!^2\ (2d-4j+1)!}{(2j)!(d-2j)!^2(2d-2j+1)!}\,  (x \, y)^{2j} \\ 
& \times \bbone\{j\ge i\}\times
\frac{(d-2i)!^2 (2j)!(2d-2j+1)!}{(2j-2i)!(2d-2i-2j+1)! d!^2} \\ 
& \times
\frac{(d-i-j)!\ (2j-2i)!}{4^{j-i}\ (d-2i)!\ (j-i)!}
\ \times\ \Delta_Q^{j-i}\ \times\ 
\parbox{2.3cm}{\psfrag{q}{$\scriptstyle{Q}$}
\psfrag{x}{$\scriptstyle{x}$}\psfrag{y}{$\scriptstyle{y}$}
\psfrag{1}{$\scriptstyle{d-2j}$}
\includegraphics[width=2.3cm]{Sec6fig7.eps}}\ \ .
\end{aligned} \] 
We finally get the required formula after multiplying by $\Delta_Q^i$, 
and simplifying the result. 
Now (\ref{geom.involution}) is an immediate consequence of (\ref{sigmavsU})
and Proposition~\ref{triangprop}. Indeed, $g_i=2^d \, G_{0,i}$. \qed

\section{A recoupling formula for transvectants}

This section gives formulae for the coefficients $\omega$ which are needed to define the 
system $\SYS(d)$ in \S\ref{section.Sd}. However, the results proved here are substantially 
more general. First we recall some preliminaries on 6-j symbols, and then 
calculate the renormalisation coefficient for the tetrahedron graph 
(see~Proposition~\ref{tetraCG}). This calculation is then used to prove a recoupling 
formula for transvectants in Theorem~\ref{Trecoupling}. 
Then we specialise this formula to the case where $A$ and $B$ are powers of a quadratic form 
$Q$. The resulting coefficients are precisely the ones needed to build $\SYS(d)$. 
 
\subsection{Preliminaries on 6-j symbols}
In what follows we assume familiarity with the graphical formalism of~\cite{spinnet},
as well as the definitions and results about Wigner symbols collected in~\cite[\S7]{AIF}.
Let $j_1, j_2, j_3, j_{12}, j_{23}, J$ be elements in $\frac{1}{2}\NN$
where $(j_1,j_2,j_{12})$, $(j_2,j_3,j_{23})$, $(j_{12},j_3,J)$ and $(j_1,j_{23},J)$ are triads in 
the sense of~\cite[\S7.6]{AIF}. Using the notations of~\cite[\S3]{spinnet}, we have a map 
$S_{2J} \lra S_{2J}$ given by the diagram 
\[
\left\langle
\parbox{2.7cm}{\ \ \ 
\psfrag{a}{$\scriptstyle{2j_1}$}
\psfrag{b}{$\scriptstyle{2j_2}$}
\psfrag{c}{$\scriptstyle{2j_{12}}$}
\psfrag{d}{$\scriptstyle{2j_3}$}
\psfrag{e}{$\scriptstyle{2J}$}
\psfrag{f}{$\scriptstyle{2j_{23}}$}
\includegraphics[width=2cm]{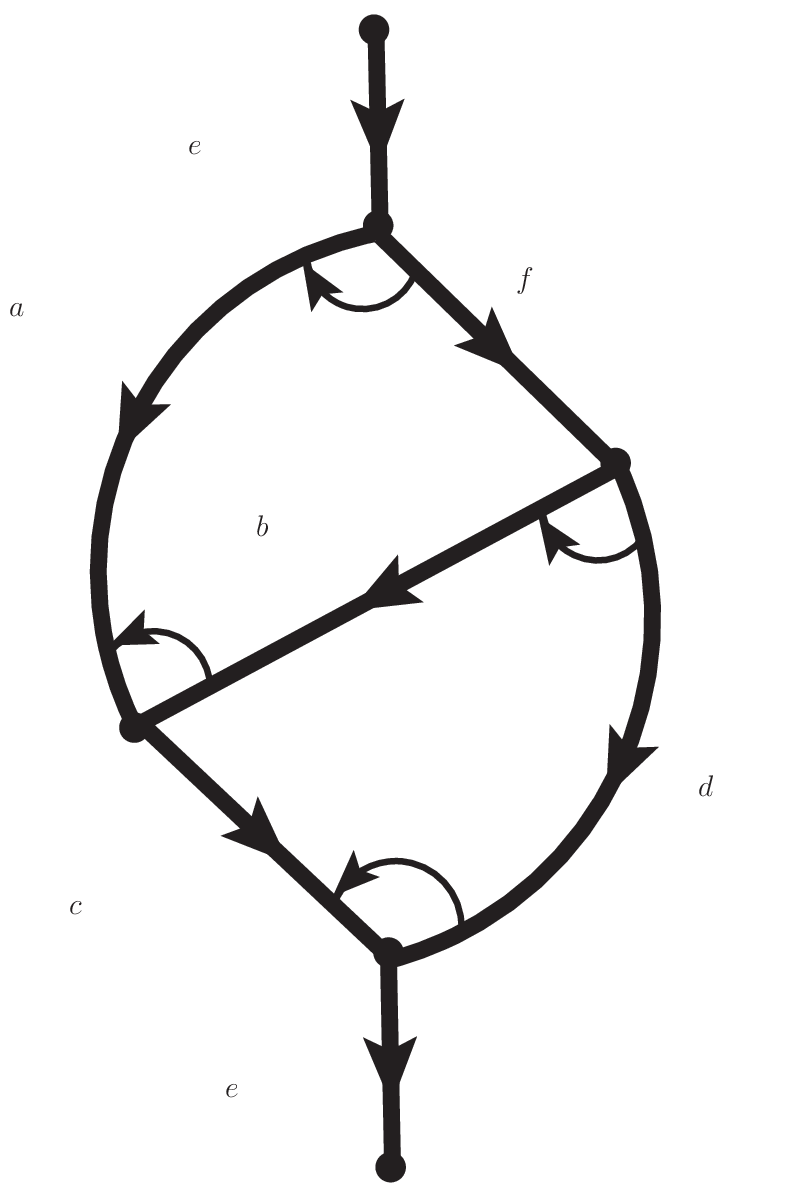}}
\right\rangle^{CG}
\] 
This map is equal to the one defined in~\cite[\S7.7]{AIF}, 
times the combinatorial normalisation factor
\[
\frac{(j_1+j_{12}-j_2)!(j_2+j_{12}-j_1)!(j_{12}+J-j_3)!(j_3+J-j_{12})!}
{(2j_1)!(2j_2)!(2j_3)!(2j_{12})!(2j_{23})!(2J)!}\ .
\]
Besides, the map considered in~\cite[\S7.7]{AIF}
is a multiple $\widetilde{\alpha}$ of the identity, where
\[
\widetilde{\alpha}=
(-1)^{j_1+j_2+j_3+J} \ \times \frac{1}{2J+1} \times 
\sqrt{\frac{P_2\, P_3}{P_1}} \times
\left\{ \begin{array}{ccc} j_1 & j_2 & j_{12} \\
j_3 & J & j_{23} \end{array} \right\}
\]
with
\[ 
\begin{aligned} 
P_1 = & \, (j_1+j_{12}-j_2)! \, (j_2+j_{12}-j_1)! \, 
(j_{12}+J-j_3)! \, (j_3+J-j_{12})!, \\
P_2 = & \, (j_1+j_{23}-J)! \, (j_1+J-j_{23})! \, 
(j_{23}+J-j_1)! \, (j_2+j_3-j_{23})! \, \times \\ 
& \, (j_2+j_{23}-j_3)! \, (j_3+j_{23}-j_2)! \, (j_1+j_2-j_{12})! \, 
(j_{12}+j_3-J)!, \\ 
P_3 = & \, (j_1+j_2+j_{12}+1)! \, (j_2+j_3+j_{23}+1)! \, \times\\ 
& \, (j_1+j_{23}+J+1)! \, (j_{12}+j_3+J+1)!. 
\end{aligned} \] 
Here $\left\{ \begin{array}{ccc} j_1 & j_2 & j_{12} \\
j_3 & J & j_{23} \end{array} \right\}$
is the standard Wigner 6-j symbol from the quantum theory of angular momentum.
It is given by Racah's single sum formula (see~\cite[Appendix B]{Racah}). 
If one writes
\[
\Delta(a,b,c)=\sqrt{\frac{(a+b-c)!(a+c-b)!(b+c-a)!}{(a+b+c+1)!}}\, ,
\]
for any triad $(a,b,c)$ of half integers, then Racah's formula is
\[
\left\{ \begin{array}{ccc} j_1 & j_2 & j_{12} \\
j_3 & J & j_{23} \end{array} \right\}=
\Delta(j_1,j_2,j_{12})
\Delta(j_2,j_3,j_{23})
\Delta(j_1,j_{23},J)
\Delta(j_{12},j_3,J)
\]
\[
\times
\sum_n (-1)^n (n+1)! \times
\left[(n-j_1-j_2-j_{12})!)(n-j_2-j_3-j_{23})!
\right. 
\]
\[
\times (n-j_1-j_{23}-J)!(n-j_{12}-j_3-J)!
(j_1+j_2+j_3+J-n)!
\]
\[
\left.
(j_2+j_{12}+j_{23}+J-n)!
(j_1+j_3+j_{12}+j_{23}-n)!
\right]^{-1}\ .
\]
The range of summation is over all integers $n$ for which
the arguments of the seven factorials in the denominator are nonnegative.

The following result is an immediate corollary of~\cite[\S3,Proposition 2]{spinnet}. 

\begin{Proposition}\label{tetraCG}
The Clebsch-Gordan normalisation of the tetrahedron graph is given by
the formula
\[
\left\langle
\parbox{3.5cm}{\ \ 
\psfrag{a}{$\scriptstyle{2j_1}$}
\psfrag{b}{$\scriptstyle{2j_2}$}
\psfrag{c}{$\scriptstyle{2j_{12}}$}
\psfrag{d}{$\scriptstyle{2j_3}$}
\psfrag{e}{$\scriptstyle{2J}$}
\psfrag{f}{$\scriptstyle{2j_{23}}$}
\includegraphics[width=2.7cm]{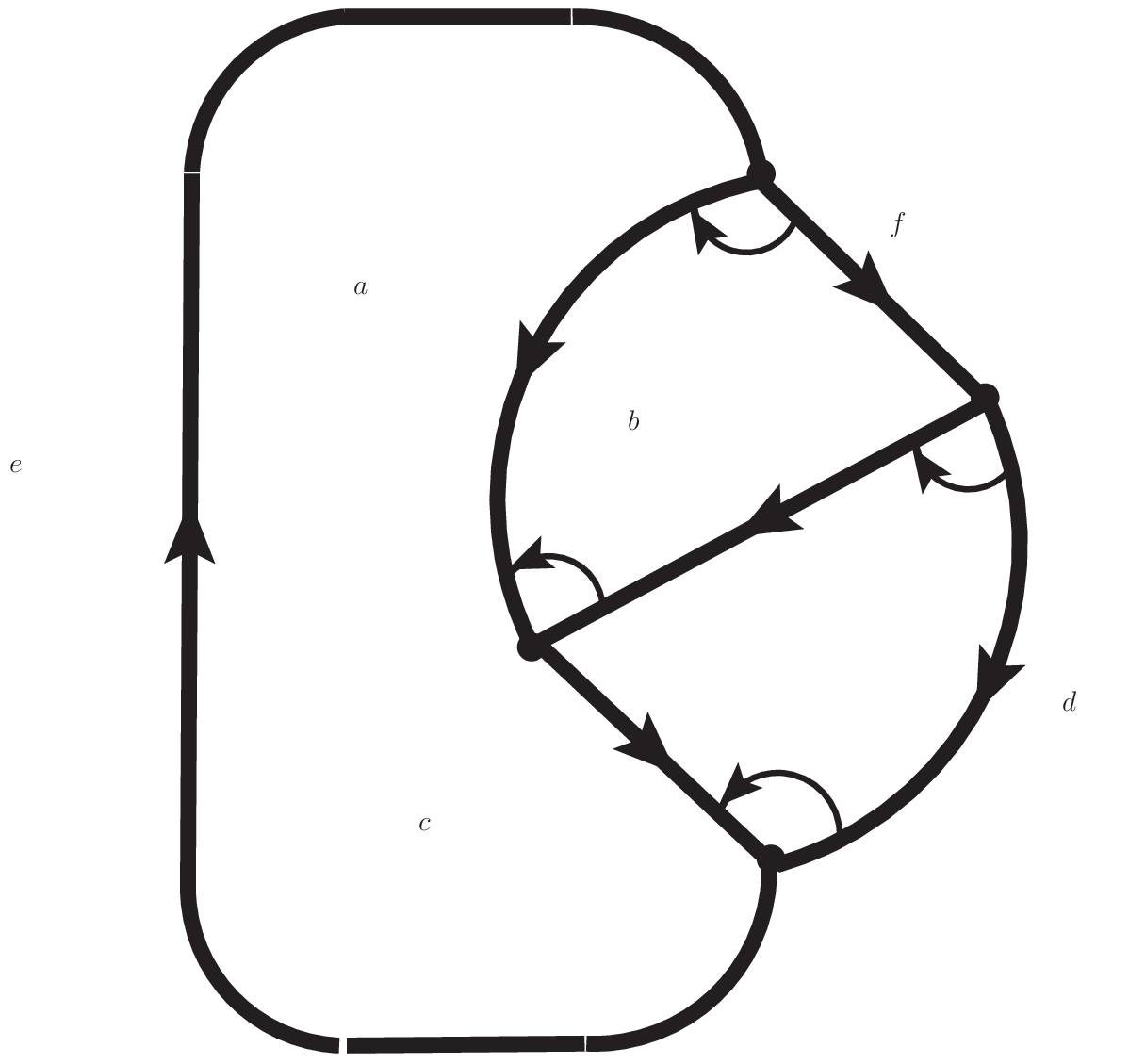}}
\right\rangle^{CG}
=(-1)^{j_1+j_2+j_3+J}
\times\frac{\mathcal{V}}{\mathcal{E}}
\]
\[
\times\sum\limits_{n=\max(T_1,T_2,T_3,T_4)}^{\min(S_1,S_2,S_3)}
\frac{(-1)^n (n+1)!}{(n-T_1)!(n-T_2)!(n-T_3)!(n-T_4)!(S_1-n)!(S_2-n)!(S_3-n)!}. 
\]
Here 
$\mathcal{E}$ is the product of the factorials of the edge labels, namely,
\[
\mathcal{E}=(2j_1)! \, (2j_2)! \, (2j_3)! \, (2j_{12})! \, (2j_{23})! \, (2J)! \, , 
\]
while $\mathcal{V}$ is the product over vertices with edge labels $2a,2b,2c$
of the factorial combination $(a+b-c)!(a+c-b)!(b+c-a)!$, i.e.,
\[ \begin{aligned}
\mathcal{V}= & (j_1+j_2-j_{12})!(j_1+j_{12}-j_2)!(j_2+j_{12}-j_1)!\\
{ } & \times (j_2+j_3-j_{23})!(j_2+j_{23}-j_3)!(j_3+j_{23}-j_3)!\\
{ } & \times (j_1+j_{23}-J)!(j_1+J-j_{23})!(j_{23}+J-j_1)!\\
{ } & \times (j_{12}+j_3-J)!(j_{12}+J-j_3)!(j_3+J-j_{12})!\ ,
\end{aligned}
\]
and 
\[
\begin{array}{ll} 
S_1= j_1+j_2+j_3+J\ , & T_1= j_1+j_2+j_{12}\ , \\ 
S_2= j_2+j_{12}+j_{23}+J\ , & T_2= j_2+j_3+j_{23}\ ,\\
S_3= j_1+j_3+j_{12}+j_{23}\ , & T_3= j_1+j_{23}+J\ ,\\ 
{} & T_4= j_{12}+j_3+J. \end{array} \]

\end{Proposition}

\subsection{Sketch of an alternate proof} 
We sketch an alternate derivation of Proposition \ref{tetraCG}, 
which uses the formula for the tetrahedral spin network in~\cite[Lemma 6.3]{GvdV}. 
The quickest derivation
of the latter is via the so-called chromatic method of Penrose and Moussouris 
(see e.g.~\cite{KL,Mous,West}).

Proposition \ref{tetraCG} immediately follows from~\cite[Lemma 6.3]{GvdV} and
the negative dimensionality theorem~\cite[Theorem 4]{spinnet}, except for the determination 
of the sign $\mu$. 
From the proof of~\cite[Theorem 4]{spinnet}, we know that the sign 
is given by $\mu=\tilde{\mu}\times (-1)^{\frac{k}{2}}$
where $k$ is total number of epsilon arrows in the microscopic
Clebsch-Gordan diagram, namely,
\[
\begin{aligned}
k= & (j_1+j_2-j_{12})+(j_2+j_3-j_{23})+(j_1+j_{23}-J)+(j_{12}+j_3-J)\\
= & 2(j_1+j_2+j_3-J), 
\end{aligned}
\]
and $\tilde{\mu}=(-1)^{C(\vec{\sigma})+N(\vec{\sigma})+B_+(\vec{\sigma})}$
is an invariant of the choice $\vec{\sigma}$ of branching permutations of the strands
at each edge of the tetrahedron graph.
Let us choose the identity for all these permutations so that there is no crossing, 
hence $C(\vec{\sigma})=0$. If we label the finite faces of the graph (drawn on the plane) as in:
\[
\parbox{2.7cm}{
\psfrag{a}{$\scriptstyle{\alpha}$}
\psfrag{b}{$\scriptstyle{\beta}$}
\psfrag{c}{$\scriptstyle{\gamma}$}
\psfrag{e}{$\scriptstyle{\ }$}
\includegraphics[width=2.7cm]{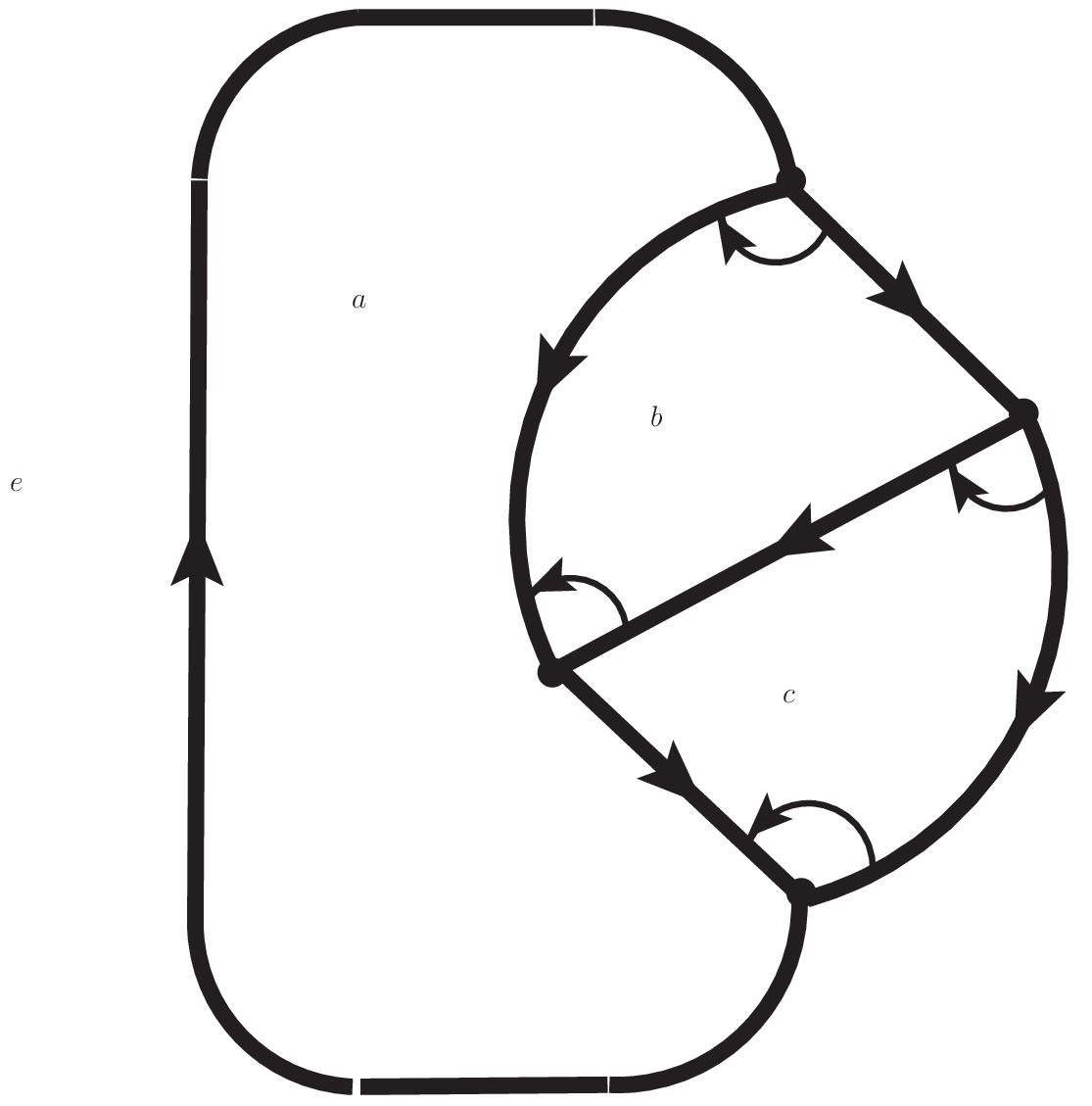}}
\]
then it is easy to see that there can only be seven types of closed curves present,
denoted by $\alpha$, $\beta$, $\gamma$, $\alpha\beta$, $\alpha\gamma$, $\beta\gamma$ and
$\alpha\beta\gamma$ according to which union of faces has the given curve as a boundary.
Let us denote by $n_\alpha,\ldots,n_{\alpha\beta\gamma}$
the number of curves present for each type. Then 
\[
N(\vec{\sigma})=n_\alpha+n_\beta+\cdots+n_{\alpha\beta\gamma}, 
\]
and
\[
B_+(\vec{\sigma})=n_\alpha B_+(\alpha)+n_\beta B_+(\beta)+\cdots+n_{\alpha\beta\gamma}
B_+(\alpha\beta\gamma), 
\]
where $B_+(\cdot)$ is number of good gate crossings for the given type of curve.
From the diagram, we can read off the following values: 
\begin{equation} \begin{array}{llll} 
B_+(\alpha)=0, & B_+(\beta)=1, & B_+(\gamma)=1, & B_+(\alpha\beta)=0, \\ 
B_+(\alpha\gamma)=0, & B_+(\beta\gamma)=1, & B_+(\alpha\beta\gamma)=0 
\end{array} \end{equation} 
for clockwise direction of travel along the curves.
Therefore, the exponent of $-1$ in $\tilde{\mu}$ is congruent modulo $2$ to the number
of curves which contain $\alpha$:
\[
n_\alpha+n_{\alpha\beta}+n_{\alpha\gamma}+n_{\alpha\beta\gamma}=2J\ ,
\]
i.e., the number of strands through the edge with label $2J$.
Thus, the sign which was left undetermined in~\cite[Theorem 4]{spinnet}
is here given by $\mu=(-1)^{j_1+j_2+j_3+J}$.

\subsection{} 
We now prove a recoupling formula for transvectants (cf.~\cite[Fig.~4]{GvdV}). 

Let $A$, $B$, $C$ be three binary forms of respective orders $a$, $b$ and $c$.
Ler $r,s$ be two integers such that $0\le r\le \min(b,c)$ 
and $0\le s\le \min(a, b+c-2r)$. Let 
\[ T_1 = \frac{r!(b-r)!(c-r)!s!(a-s)!(b+c-2r-s)!}{(b+c-2r)!}. \] 
For $k \in \ZZ$, define 
\[ \begin{aligned} 
T_2(k) = \sum_{z\in\ZZ} & \bbone\left\{
\begin{array}{c}
{\rm arguments\ of}\\
{\rm factorials\ in}\\
{\rm denominator\ are\ }\ge 0 \end{array} \right\} (-1)^z (z+1)! \times  \\ 
 [ \, & (z-a-b+k)! (z-b-c+r)!(z-a-b-c+2r+s)! \\ 
& (z-a-b-c+r+s+k)!(a+b+c-r-s-z)! \\ 
& (a+b+c-r-k-z)!(a+2b+c-2r-s-k-z)! \, ]^{-1}, 
\end{aligned} \] 
and finally 
\[ \begin{aligned} 
\theta_k = & (-1)^{a+b+c+r+s}  \, T_1 \times 
\ \bbone\left\{
\begin{array}{ll}
k\ge 0, & k\ge r+s-c\\
k\le a,  & k\le a+b-r-s\\
k\le b,  & k\le r+s
\end{array} \right\}  \\ 
& \times\frac{(a+b-2k+1)! \; T_2(k)}
{(a+b-k+1)!(a+b+c-r-s-k+1)!}. 
\end{aligned} \] 
\begin{Theorem} 
With notation as above, we have an expansion 
\[ (A,(B,C)_r)_s= \sum\limits_{k \in \ZZ} \; \theta_k \, ((A,B)_k,C)_{r+s-k}. \] 
\label{Trecoupling} \end{Theorem} 
\demo 
First recall from~\cite[Eq. 16]{spinnet} the macroscopic form of the Clebsch-Gordan series
\begin{equation}
\left\langle
\parbox{1.3cm}{\psfrag{m}{$\scriptstyle{m}$}
\psfrag{n}{$\scriptstyle{n}$}
\includegraphics[width=1.3cm]{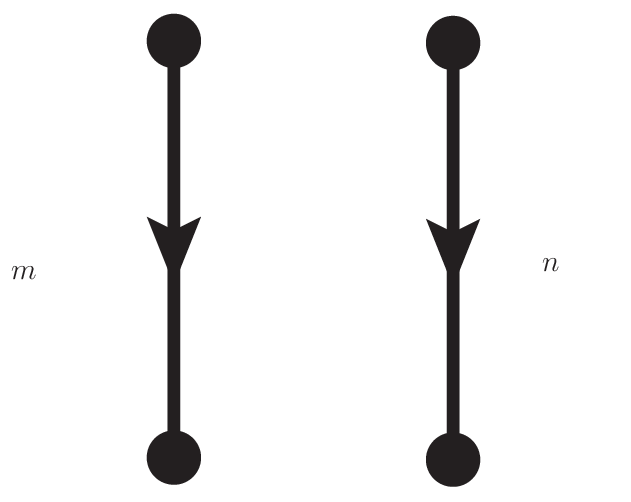}}
\right\rangle^{CG} =
\sum\limits_{k=0}^{\min(m,n)}
\frac{ \left( \begin{array}{c} m\\ k \end{array} \right) 
\left( \begin{array}{c} n\\ k \end{array} \right) }
{ \left( \begin{array}{c} m+n-k+1\\ k \end{array} \right) }
\left\langle
\ \ 
\parbox{2cm}{\psfrag{m}{$\scriptstyle{m}$}
\psfrag{n}{$\scriptstyle{n}$}
\psfrag{j}{$\scriptstyle{m+n-2k}$}
\includegraphics[width=1.4cm]{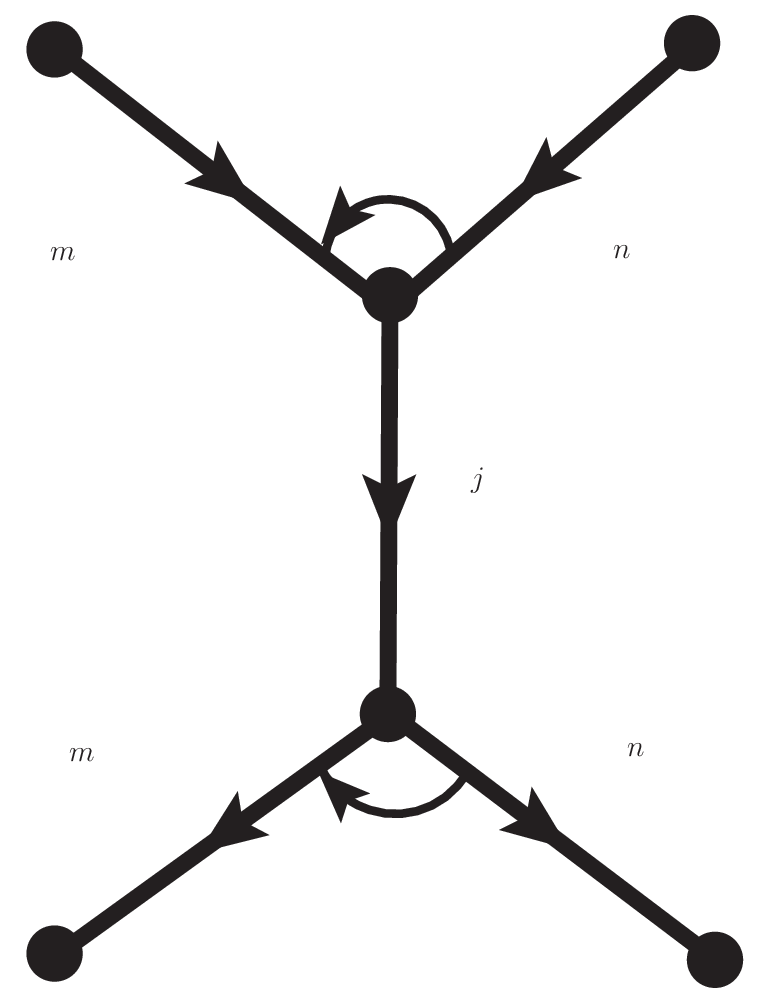}}
\right\rangle^{CG}\ .
\label{CGdirect}
\end{equation}
Using the same trick as in~\cite[p. 51]{spinnet}, we have the following variant 
form of this identity: 
\begin{equation}
\left\langle
\parbox{1.3cm}{\psfrag{m}{$\scriptstyle{m}$}
\psfrag{n}{$\scriptstyle{n}$}
\includegraphics[width=1.3cm]{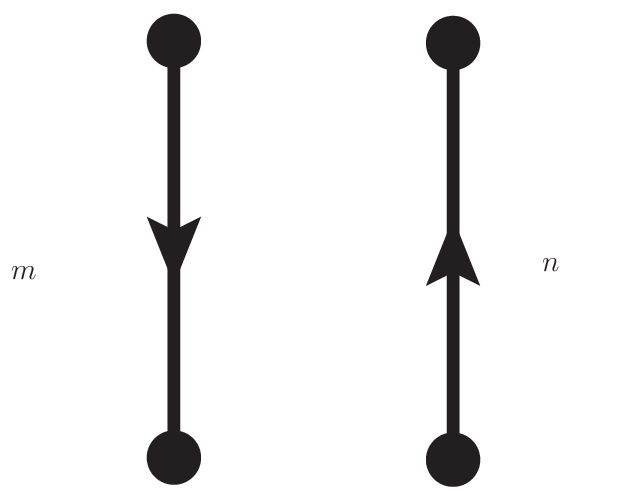}}
\right\rangle^{CG} =
\sum\limits_{k=0}^{\min(m,n)}
\frac{
\left(
\begin{array}{c}
m\\
k
\end{array}
\right)
\left(
\begin{array}{c}
n\\
k
\end{array}
\right)
}{
\left(
\begin{array}{c}
m+n-k+1\\
k
\end{array}
\right)
}
\left\langle
\ \ 
\parbox{2cm}{\psfrag{m}{$\scriptstyle{m}$}
\psfrag{n}{$\scriptstyle{n}$}
\psfrag{j}{$\scriptstyle{m+n-2k}$}
\includegraphics[width=1.4cm]{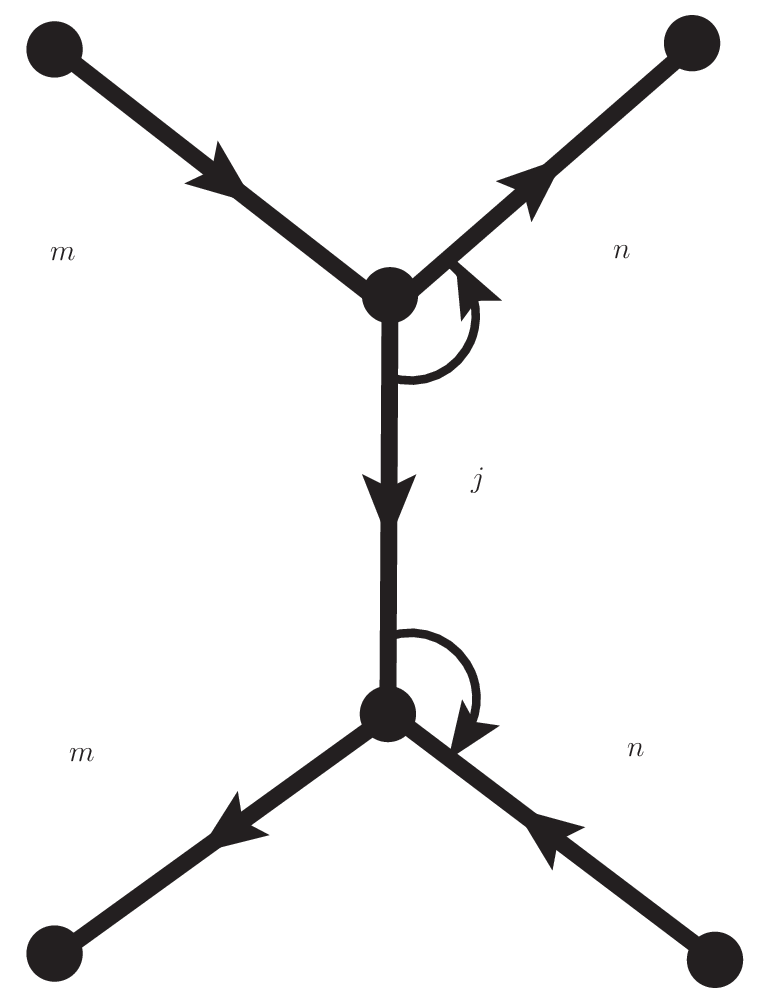}}
\right\rangle^{CG}
\label{CGvariant}
\end{equation}
Now, as explained in~\cite[p. 22--23]{spinnet},
one can write
\[
(A,(B,C)_r)_s=
\parbox{3.5cm}{\ \ 
\psfrag{a}{$\scriptstyle{a}$}
\psfrag{b}{$\scriptstyle{b}$}
\psfrag{c}{$\scriptstyle{c}$}
\psfrag{1}{$\scriptstyle{b+c-2r}$}
\psfrag{2}{$\scriptstyle{a+b+c-2r-2s}$}
\psfrag{x}{$\scriptstyle{x}$}
\psfrag{A}{$\scriptstyle{A}$}
\psfrag{B}{$\scriptstyle{B}$}
\psfrag{C}{$\scriptstyle{C}$}
\includegraphics[width=3.5cm]{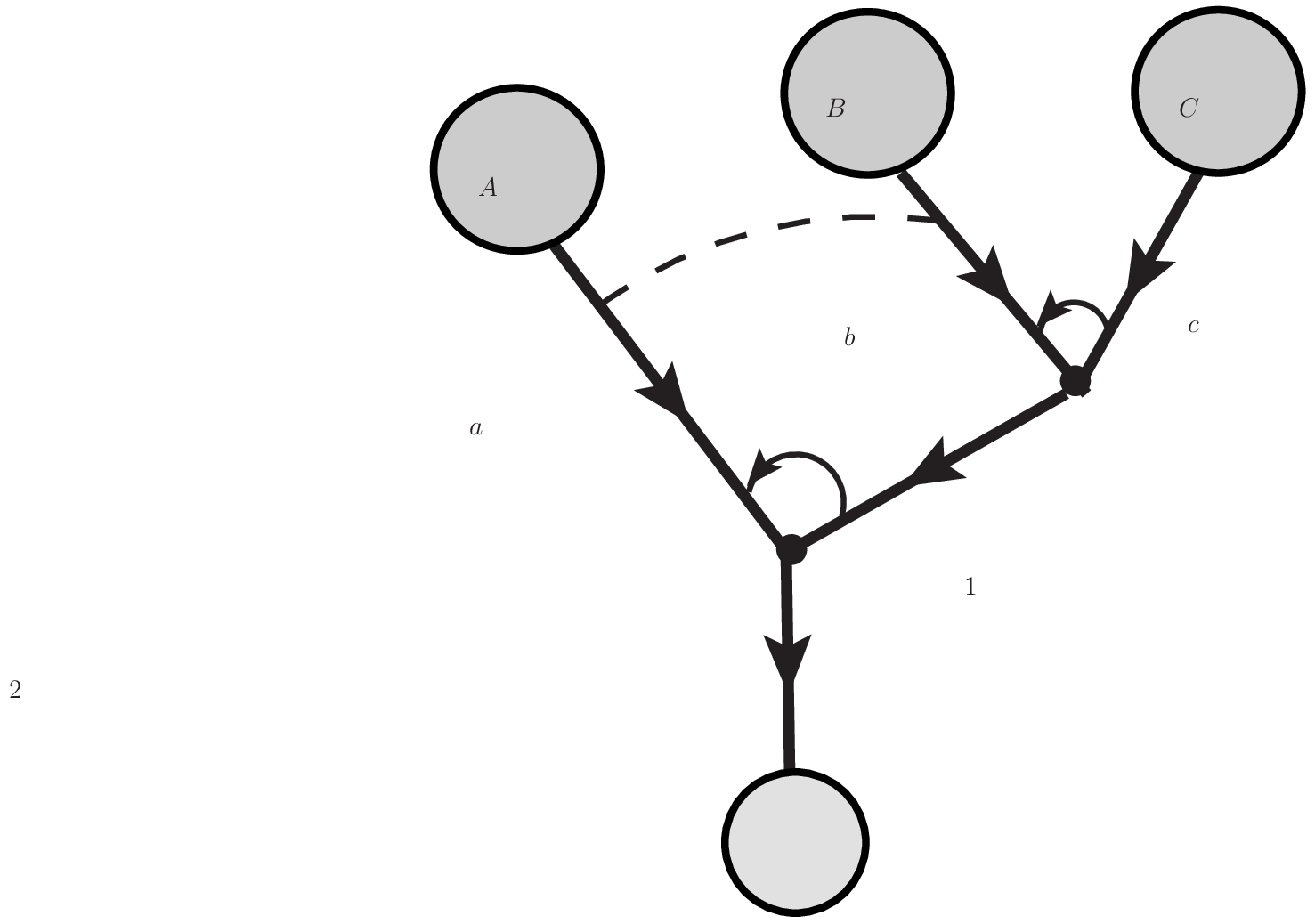}}
\]
where the dotted lines indicate where (\ref{CGdirect}) will next be used. The outcome is
\[
(A,(B,C)_r)_s=
\sum\limits_{k_1=0}^{\min(a,b)}
\frac{
\left(
\begin{array}{c}
a\\
k_1
\end{array}
\right)
\left(
\begin{array}{c}
b\\
k_1
\end{array}
\right)
}{
\left(
\begin{array}{c}
a+b-k_1+1\\
k_1
\end{array}
\right)
}
\times
\parbox{4.5cm}{\ \ 
\psfrag{a}{$\scriptstyle{a}$}
\psfrag{b}{$\scriptstyle{b}$}
\psfrag{c}{$\scriptstyle{c}$}
\psfrag{0}{$\scriptstyle{a+b-2k_1}$}
\psfrag{1}{$\scriptstyle{b+c-2r}$}
\psfrag{2}{$\scriptstyle{a+b+c-2r-2s}$}
\psfrag{x}{$\scriptstyle{x}$}
\psfrag{A}{$\scriptstyle{A}$}
\psfrag{B}{$\scriptstyle{B}$}
\psfrag{C}{$\scriptstyle{C}$}
\includegraphics[width=3.8cm]{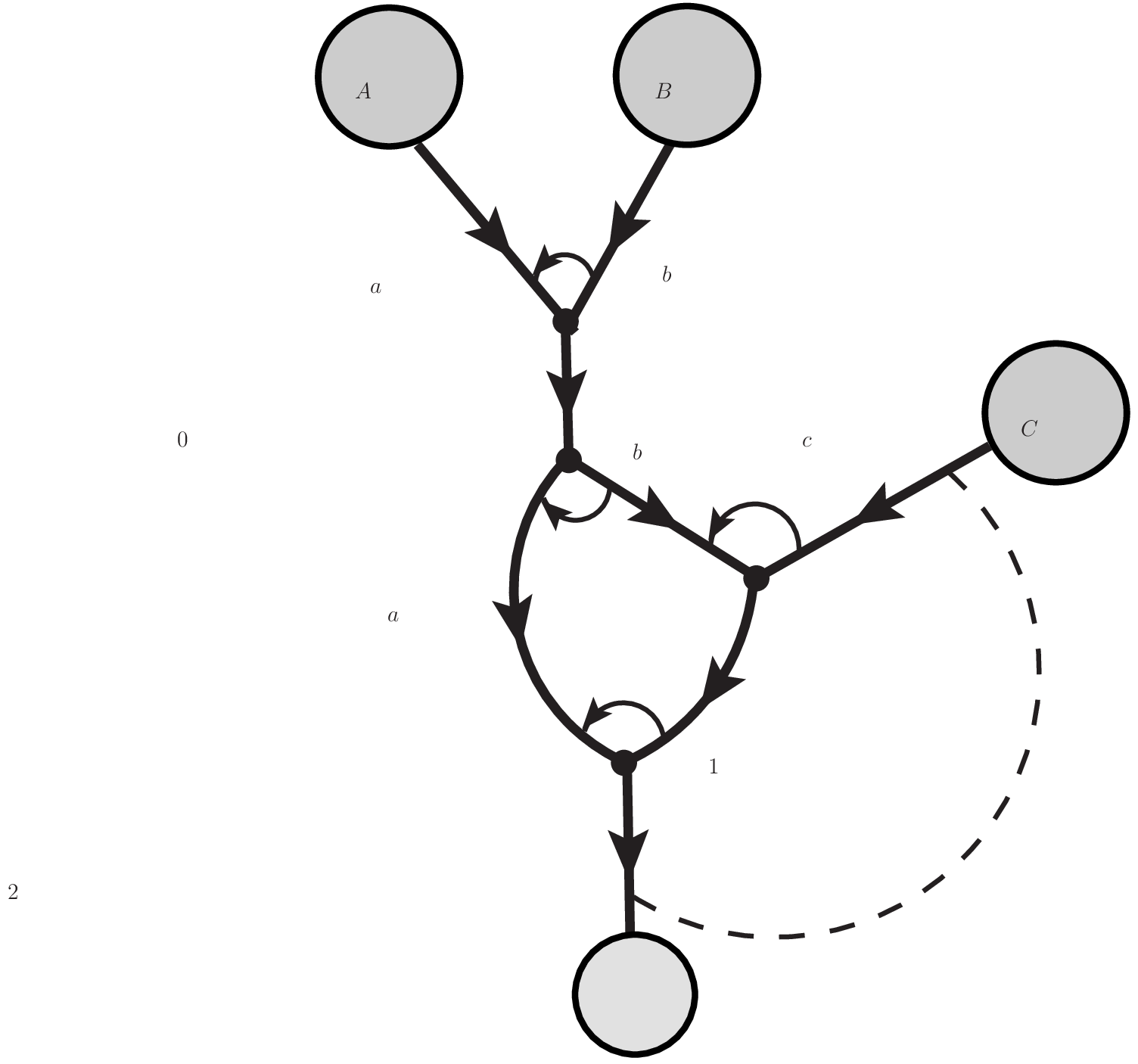}}
\]
\[
=\sum\limits_{k_1=0}^{\min(a,b)}
\frac{
\left(
\begin{array}{c}
a\\
k_1
\end{array}
\right)
\left(
\begin{array}{c}
b\\
k_1
\end{array}
\right)
}{
\left(
\begin{array}{c}
a+b-k_1+1\\
k_1
\end{array}
\right)
}
\times
\sum\limits_{k_2=0}^{\min(c,a+b+c-2r-2s)}
\frac{
\left(
\begin{array}{c}
c\\
k_2
\end{array}
\right)
\left(
\begin{array}{c}
a+b+c-2r-2s\\
k_2
\end{array}
\right)
}{
\left(
\begin{array}{c}
a+b+2c-2r-2s-k_2+1\\
k_2
\end{array}
\right)
}
\]
\[
\times
\parbox{5cm}{
\psfrag{a}{$\scriptstyle{a}$}
\psfrag{b}{$\scriptstyle{b}$}
\psfrag{c}{$\scriptstyle{c}$}
\psfrag{0}{$\scriptstyle{a+b-2k_1}$}
\psfrag{1}{$\scriptstyle{b+c-2r}$}
\psfrag{2}{$\scriptstyle{a+b+c-2r-2s}$}
\psfrag{3}{$\scriptstyle{a+b+2c-2r-2s-2k_2}$}
\psfrag{x}{$\scriptstyle{x}$}
\psfrag{A}{$\scriptstyle{A}$}
\psfrag{B}{$\scriptstyle{B}$}
\psfrag{C}{$\scriptstyle{C}$}
\includegraphics[width=3.8cm]{Fig6.eps}}
\]
\vskip 2cm
\noindent
by (\ref{CGvariant}). Now by the graphical Schur lemma (see~\cite[Proposition 2]{spinnet}), 
this becomes 
\[ (A,(B,C)_r)_s
=\sum\limits_{k_1=0}^{\min(a,b)}
\frac{ \left( \begin{array}{c}
a\\ k_1 \end{array} \right) \left( \begin{array}{c}
b\\ k_1 \end{array} \right) }{ \left( \begin{array}{c}
a+b-k_1+1\\ k_1 \end{array} \right) } \]
\[
\times
\sum\limits_{k_2=0}^{\min(c,a+b+c-2r-2s)}
\frac{
\left(
\begin{array}{c}
c\\
k_2
\end{array}
\right)
\left(
\begin{array}{c}
a+b+c-2r-2s\\
k_2
\end{array}
\right)
}{
\left(
\begin{array}{c}
a+b+2c-2r-2s-k_2+1\\
k_2
\end{array}
\right)
}
\]
\[
\times \frac{1}{a+b-2k_1+1}
\times \bbone
\{a+b-2k_1=a+b+2c-2r-2s-2k_2\}
\]
\[
\times\left\langle
\parbox{5.7cm}{\ \ \ 
\psfrag{a}{$\scriptstyle{a}$}
\psfrag{b}{$\scriptstyle{b}$}
\psfrag{c}{$\scriptstyle{c}$}
\psfrag{0}{$\scriptstyle{a+b-2k_1}$}
\psfrag{1}{$\scriptstyle{b+c-2r}$}
\psfrag{2}{$\scriptstyle{a+b+c-2r-2s}$}
\includegraphics[width=5cm]{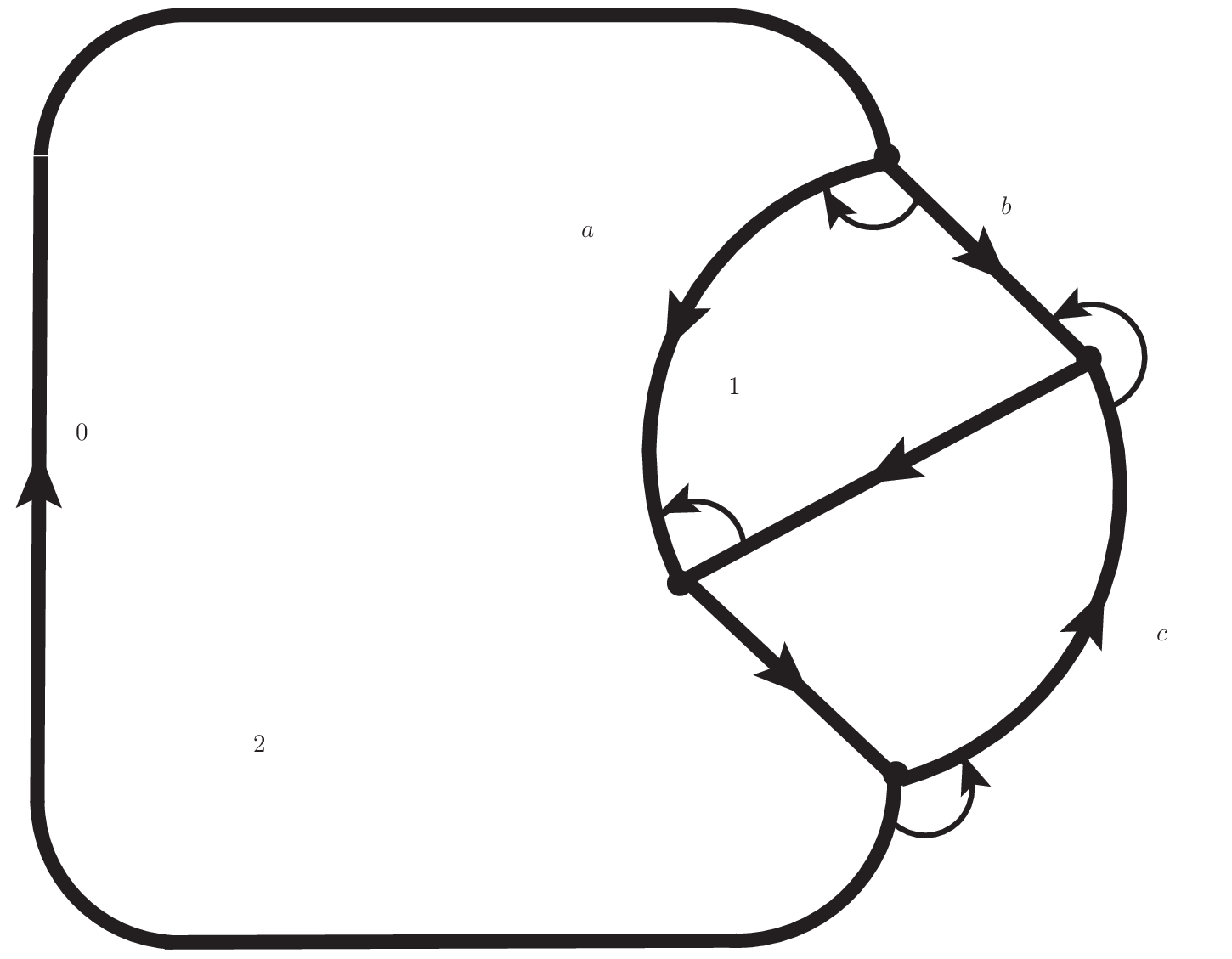}}
\right\rangle^{CG}
\ \ \times
\ \parbox{4cm}{\ \ \ 
\psfrag{a}{$\scriptstyle{a}$}
\psfrag{b}{$\scriptstyle{b}$}
\psfrag{c}{$\scriptstyle{c}$}
\psfrag{x}{$\scriptstyle{x}$}
\psfrag{A}{$\scriptstyle{A}$}
\psfrag{B}{$\scriptstyle{B}$}
\psfrag{C}{$\scriptstyle{C}$}
\psfrag{0}{$\scriptstyle{a+b-2k_1}$}
\psfrag{2}{$\scriptstyle{a+b+c-2r-2s}$}
\includegraphics[width=4cm]{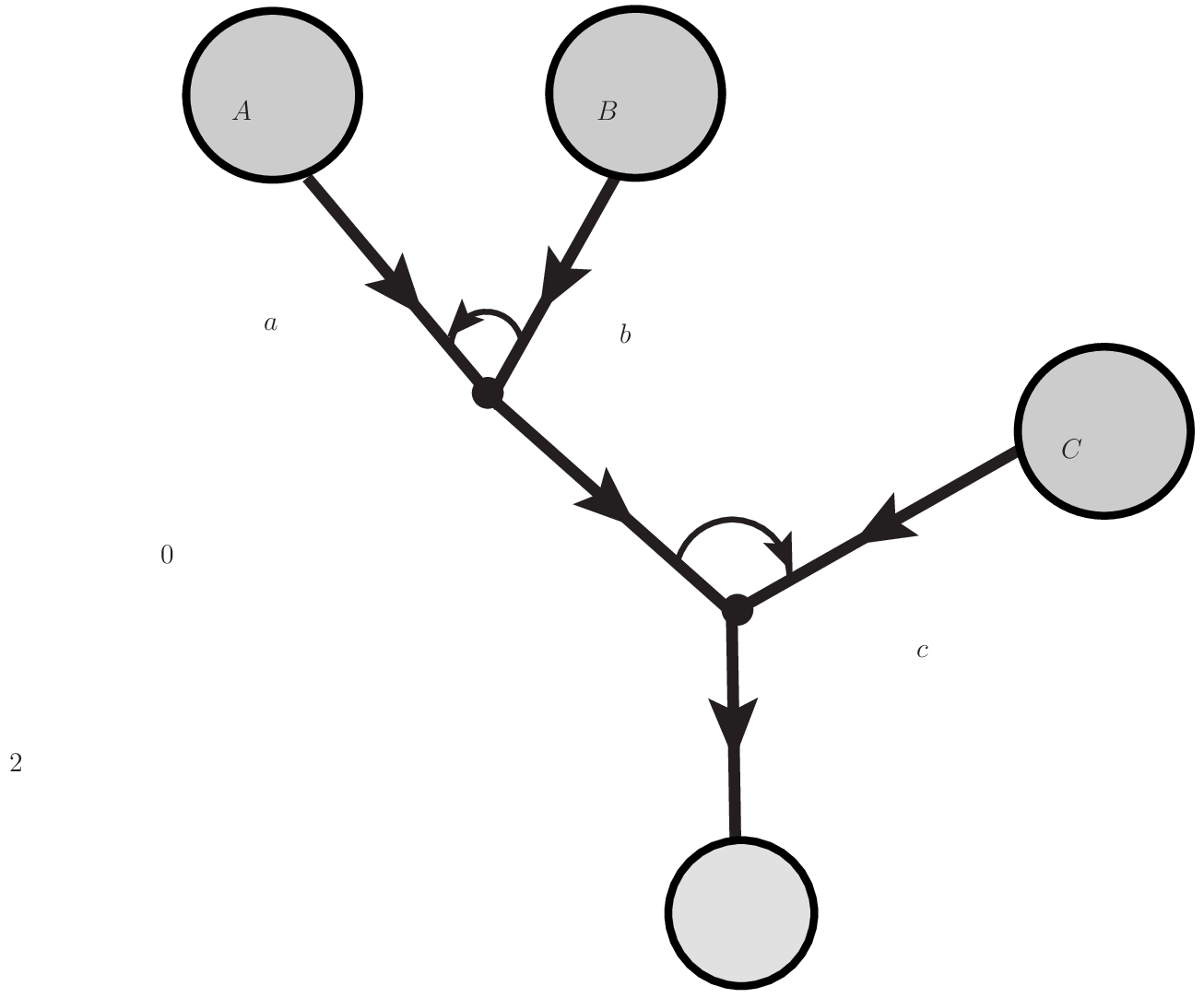}}\ .
\]
By reversing the epsilon arrows in the lower transvection vertex, we have 
\[
\parbox{4cm}{\ \ \ 
\psfrag{a}{$\scriptstyle{a}$}
\psfrag{b}{$\scriptstyle{b}$}
\psfrag{c}{$\scriptstyle{c}$}
\psfrag{x}{$\scriptstyle{x}$}
\psfrag{A}{$\scriptstyle{A}$}
\psfrag{B}{$\scriptstyle{B}$}
\psfrag{C}{$\scriptstyle{C}$}
\psfrag{0}{$\scriptstyle{a+b-2k_1}$}
\psfrag{2}{$\scriptstyle{a+b+c-2r-2s}$}
\includegraphics[width=4cm]{Fig8.eps}}
\ \ = (-1)^{r+s-k_1} \times ((A,B)_{k_1},C)_{r+s-k_1}\ .
\]
Furthermore, by inserting the matrix identity $\epsilon \epsilon^{\rm T}={\rm id}$ in the strands
going through the edge with label $c$, we get 
\[
\left\langle
\parbox{4.7cm}{\ \ \ 
\psfrag{a}{$\scriptstyle{a}$}
\psfrag{b}{$\scriptstyle{b}$}
\psfrag{c}{$\scriptstyle{c}$}
\psfrag{0}{$\scriptstyle{a+b-2k_1}$}
\psfrag{1}{$\scriptstyle{b+c-2r}$}
\psfrag{2}{$\scriptstyle{a+b+c-2r-2s}$}
\includegraphics[width=4.3cm]{Fig7.eps}}
\right\rangle^{CG}\ 
=(-1)^r
\left\langle
\parbox{4.7cm}{\ \ \ 
\psfrag{a}{$\scriptstyle{a}$}
\psfrag{b}{$\scriptstyle{b}$}
\psfrag{c}{$\scriptstyle{c}$}
\psfrag{0}{$\scriptstyle{a+b-2k_1}$}
\psfrag{1}{$\scriptstyle{b+c-2r}$}
\psfrag{2}{$\scriptstyle{a+b+c-2r-2s}$}
\includegraphics[width=4.3cm]{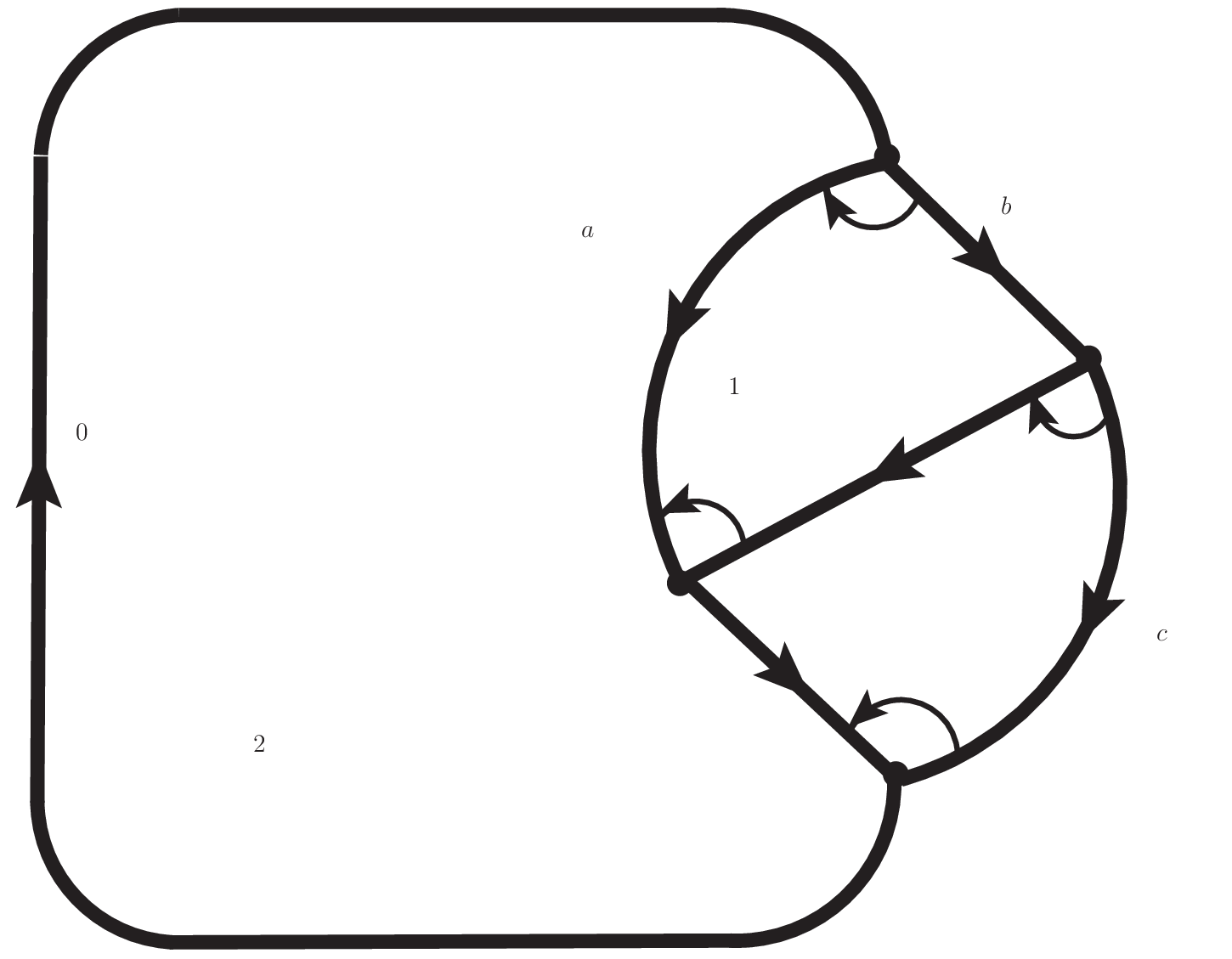}}
\right\rangle^{CG}
\]
\vskip 0.8cm
\noindent
Now evaluate the latter tetrahedral diagram using Proposition~\ref{tetraCG} with the 
values 
\[
\begin{array}{lll} 
2j_1= a, & 2j_2= b+c-2r, & 2j_3=  c,\\ 
2j_{12}= a+b+c-2r-2s, & 2j_{23}= b, & 2J= a+b-2k_1. 
\end{array} \]
Once the resulting expression is simplified, Theorem~\ref{Trecoupling} follows immediately. \qed

\subsection{The $\omega$ coefficients} \label{section.omega} 

Let $a,b,r,s$ denote nonnegative integers such that 
$r \le \min \, \{d,2b\}$, and $s \le \min \, \{2a,2b+d-2r\}$. Let $\Delta = - 2 \, (Q,Q)_2$. 
Given an arbitrary integer $t$ in the range 
\begin{equation}  \begin{aligned} 
{} &\max \, \{|a+b-r-s|,|a-b|\} \le t \le \min \, \{a+b+d-r-s,a+b\}, \\ 
& \; t \equiv a+b \; (\text{mod} \, 2); \end{aligned} 
\label{t.range} \end{equation} 
define rational numbers $P_1,P_2,P_3$ as follows. 
\[ \begin{aligned} 
P_1 = \, & \{a! \, b! \, r! \, s! \, (d-r)! \, (2b-r)! \, (2a-s)! \, (2b+d-2r-s)! \} \\ 
& \{ (2a)! \, (2b)! \, (2b+d-2r)! \}^{-1}, \\ 
P_2 = \, & \{ (2t+1)! \, (\frac{a+b+t}{2})! \, (a+b-t)! \, (a-b+t)! \, (b-a+t)! \} \, \times \\ 
& \{t! \, (a+b+t+1)! \, (a+b+d-r-s+t+1)! \, (\frac{a+b-t}{2})! \\ 
& (\frac{a-b+t}{2})! \, (\frac{b-a+t}{2})! \, \}^{-1}, 
\end{aligned} \] 
and 
\[ 
\begin{aligned} 
P_3 = \sum\limits_z \, & (-1)^z \, (z+1)! \, \times \\ 
& \{(z-a-b-t)! \, (z-2b-d+r)! \, (z-2a-2b-d+2r+s)! \, \\ 
& (z-a-b-d+r+s-t)! \, (a+b+d-r+t-z)! \, \\ 
& (2a+2b+d-r-s-z)! \, (a+3b+d-2r-s+t-z)!\}^{-1}, 
\end{aligned} \] 
where the sum is quantified over all integers $z$ such that the arguments of all 
factorials in the denominator are nonnegative. Now, let 
\begin{equation} \omega(a,b;r,s;t) = (-1)^{d+r+s+\frac{1}{2}(a+b-t)} \, P_1 \, P_2 \, P_3. 
\label{formula.omega} \end{equation} 
\begin{Proposition} \sl \label{proposition.omega} 
With notation as above, we have a formula 
\[ (Q^a,(Q^b,F)_r)_s = \sum\limits_t 
\omega(a,b;r,s;t) \, \Delta^{\frac{a+b-t}{2}} \, (Q^t,F)_{r+s-a-b+t}, \] 
where the sum is quantified over (\ref{t.range}). \end{Proposition} 
\demo 
Specialize~Theorem \ref{Trecoupling} to the case $A=Q^a$, $B=Q^b$, and $C=F$. The 
terms of the form $(Q^a,Q^b)$ can be simplified using~\cite[Proposition 6.1]{Advances}, 
and the result follows. \qed 

\smallskip 

For instance, if $d=5$, then $(Q^5,(Q^6,F)_2)_4$ can be expanded into 
\[ \frac{95}{286286} \, \Delta^3 \, (Q^5,F)_0 + 
\frac{575}{1123122} \, \Delta^2 \, (Q^7,F)_2 + 
\frac{95}{9438} \, \Delta \, (Q^9,F)_4. \] 

\bigskip 

\begin{Lemma} \sl 
In the notation of \S\ref{section.Sd}, the coefficient 
\[ \alpha_{i,i}^{(0)} = \omega(d-2i,d-2i,d-2i,d-2i,0) \] 
is always positive. 
\end{Lemma} 
It follows that the expression $\sum\limits_i \, \alpha_{i,i}^{(0)} \, z_i^2$ from 
(\ref{condition.inv.one}) is positive definite. 

\demo 
Apply the formula for $P_3$ and notice that the sum over $z$ reduces to the single term
$z=2d-2i$. This is forced by the nonnegativity of the arguments of the second and fifth
denominator factorials. One easily checks that the sum is not vacuous, i.e., the other factorials
have nonnegative arguments. For the relevant $z$, $(-1)^z$ times the
sign in (\ref{formula.omega}) gives $+1$.
\qed

\section{Proof of the main theorem} 
\label{main.proof.sec}

In this section we will prove Theorem~\ref{MainTheorem}. 

\subsection{ }
Picking up the thread from \S\ref{proof.geometric.formula},
let us use the notation $\mathcal{P}_i=\mathcal{L}(\mathfrak{p}_i)$
and $\mathcal{T}_i=\mathcal{L}(\mathfrak{t}_i)$, for $0\le i\le n$.
We know that 
\[ \{ \mathcal{P}_i: 0 \le i \le n \}, \quad \{ \mathcal{T}_i: 0 \le i \le n \} \] 
are two linear bases for the space $\mathcal{M}_d^{SL_2}$, with transition matrix
given by Proposition~\ref{triangprop}.
We now introduce a new space $\mathcal{M}_{2d}^{SL_2}$ where compositions of maps
in $\mathcal{M}_d^{SL_2}$ will reside.

Let $\mathcal{M}_{2d}$
 be the space of maps
\[
\begin{array}{rl}
\psi: & S_d\times S_2\lra S_d \\
 & (F,Q)\longmapsto \psi(F,Q)
\end{array}
\]
which are linear in $F$ and homogeneous of degree $2d$ in $Q$.
Let $\mathcal{M}_d^{SL_2}$ be the subspace of equivariant maps.
Composition (with respect to $F$) defines a bilinear map:
\[
\begin{array}{rl}
\circ: & \mathcal{M}_d\times \mathcal{M}_d\lra \mathcal{M}_{2d} \\
 & (\phi,\psi)\longmapsto \phi\circ\psi
\end{array}
\]
where
\[
(\phi\circ\psi)\ (F,Q)=\phi(\psi(F,Q),Q)\ .
\]
Clearly, this restricts to a bilinear map $\mathcal{M}_d^{SL_2}\times \mathcal{M}_d^{SL_2}
\lra \mathcal{M}_{2d}^{SL_2}$.
Note that there is a distinguished element $\ID$ in $\mathcal{M}_{2d}^{SL_2}$
given by $\ID(F,Q)=\Delta_Q^d \, F$.
We will show that the system $\SYS(d)$ is equivalent to solving the equation
\[
\psi\circ\psi=\ID
\]
in $\mathcal{M}_{d}^{SL_2}$.
This will involve the interplay of several bases for the two spaces $\mathcal{M}_{d}^{SL_2}$
and $\mathcal{M}_{2d}^{SL_2}$.

\begin{Lemma}
${\rm dim}\ \mathcal{M}_{2d}^{SL_2}=n+1$.
\end{Lemma}
\demo
Proceeding as in Lemma \ref{first.dim.lemma},
one can write:
\begin{align*}
\mathcal{M}_{2d} & = {\rm Hom}\left(
S_{2d}(S_2^\vee)\otimes S_d, S_d\right) \\
 & = S_{2d}(S_2)\otimes S_d\otimes S_d \\
 & = S_2(S_{2d})\otimes S_d\otimes S_d \\
 & = \left(\bigoplus_{i=0}^d S_{4d-4i}\right)\otimes S_d\otimes S_d\ 
\end{align*}
as $SL_2$-representations.
Therefore,
\[
\mathcal{M}_{2d}^{SL_2}= \bigoplus_{i=0}^d 
\left(S_{4d-4i}\otimes S_d\otimes S_d\right)^{SL_2}, 
\]
where each summand either vanishes or has dimension $1$.
The nonvanishing condition is the triangle inequality $4d-4i\le 2d$. The latter is equivalent to
$d-n\le i\le d$, which holds for exactly $n+1$ values of $i$. \qed

\subsection{ }
We now introduce a first `parallel' basis $\{ \widetilde{\mathcal{P}}_i: 0 \le i \le n\}$ 
for $\mathcal{M}_{2d}^{SL_2}$. By definition, $\widetilde{\mathcal{P}}_i$ is the map
\[
(F,Q)\longmapsto
\Delta_Q^{d-i} \times\ 
\parbox{3cm}{\psfrag{f}{$\scriptstyle{F}$}
\psfrag{m}{$\scriptstyle{2i}$}\psfrag{n}{$\scriptstyle{d-2i}$}
\psfrag{q}{$\scriptstyle{Q}$}\psfrag{x}{$\scriptstyle{x}$}
\includegraphics[width=3cm]{Sec6fig2.eps}}\ \ \ \ .
\]

\begin{Lemma} The collection $\{ \widetilde{\mathcal{P}}_i: 0 \le i \le n\}$ 
is indeed a linear basis for $\mathcal{M}_{2d}^{SL_2}$.
\end{Lemma}
\demo
It is enough to show linear independence.
Suppose there is a linear relation $\sum_{i=0}^n \, \mu_i \, \widetilde{\mathcal{P}}_i=0$.
Specialize to $Q=x_1^2+x_2^2$, and then set $x_1=1$, $x_2=0$. We get
\[
\sum_{i=0}^n \, \mu_i \, \widetilde{\mathcal{P}}_i(F,Q)=
\sum_{i=0}^n \, \mu_i \, (-4)^{d-i} \, a_{2i}=0, 
\]
where the $a$'s are the coefficients of $F$ in Cayley's notation.
For generic $F$ these coefficients are independent, and therefore the $\mu$'s must vanish.
\qed

\subsection{ }
We now introduce a new collection $\{ \widetilde{\mathcal{T}_i}: 0 \le i \le n\}$ 
of elements in $\mathcal{M}_{2d}^{SL_2}$ defined using transvection.
For any $i, F$ and $Q$, let 
\[
\widetilde{\mathcal{T}}_i (F,Q)=\Delta_Q^{d-i} \, (Q^{2i},F)_{2i}\ .
\]

\begin{Lemma} The collection 
$\{ \widetilde{\mathcal{T}_i}: 0 \le i \le n\}$ is also a linear basis for
$\mathcal{M}_{2d}^{SL_2}$.
\end{Lemma}
\demo
Assume that there is a linear relation $\sum_{i=0}^{n}\mu_i\widetilde{\mathcal{T}}_i=0$.
Then we have identically in $F$, $Q$ and $x$:
\[
\sum_{i=0}^{n} \, \mu_i \, 
\Delta_Q^{d-i}\times\qquad
\parbox{3.8cm}{\psfrag{1}{$\scriptstyle{2i}$}
\psfrag{f}{$\scriptstyle{F}$}\psfrag{x}{$\scriptstyle{x}$}
\psfrag{q}{$\scriptstyle{Q}$}
\includegraphics[width=3.8cm]{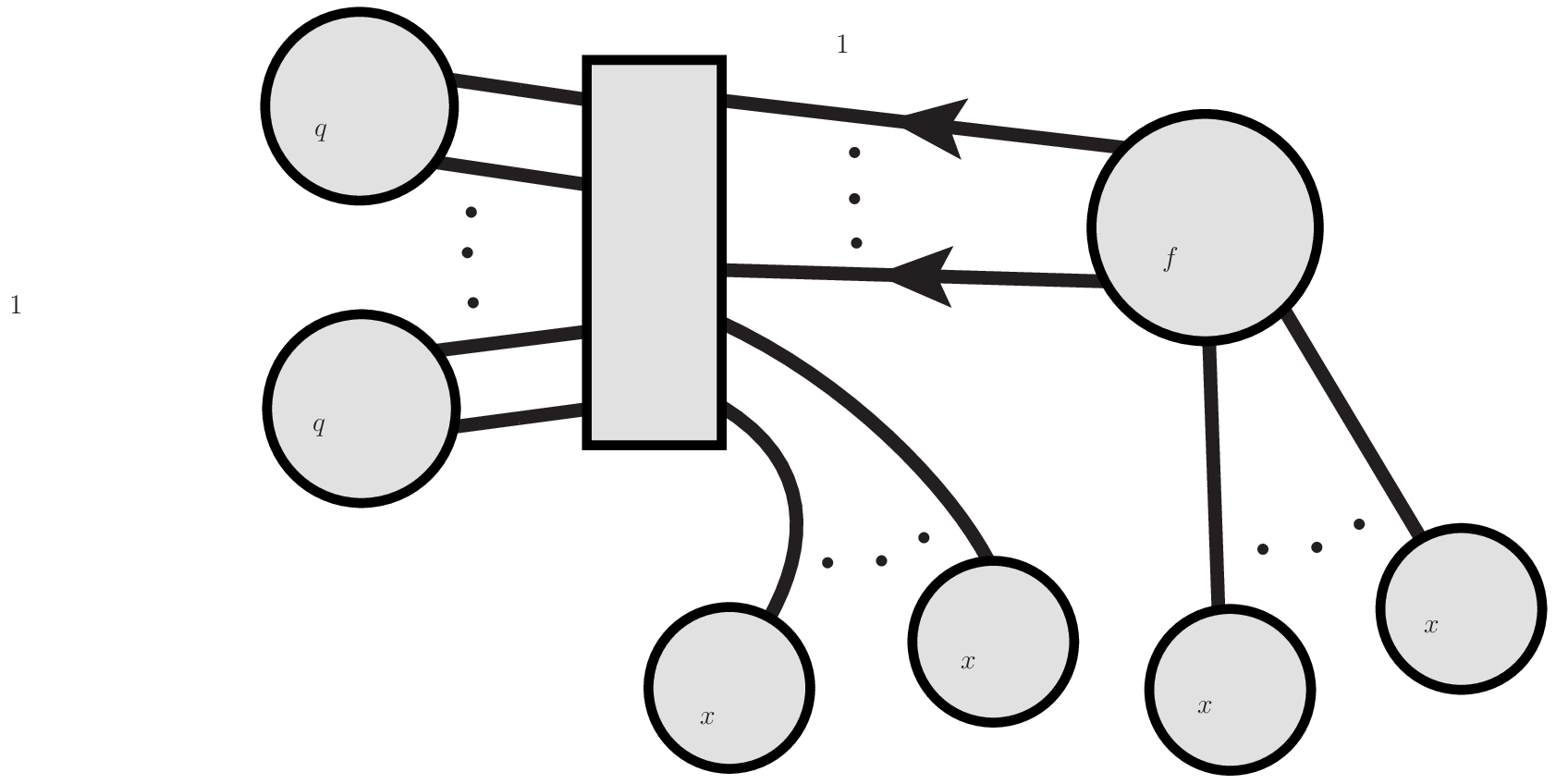}}
=0. 
\]
Now specialize $F$ to $F(x)=(x \, y)^d$, and 
write $Q$ as a product of linear forms $Q(x)=a_x \, b_x$.
Then $\Delta_Q=(a \, b)^2$, and hence 
\[
\sum_{i=0}^{n} \, \mu_i \, \Delta_Q^{d-i}
\times\qquad
\parbox{3.8cm}{\psfrag{1}{$\scriptstyle{2i}$}
\psfrag{f}{$\scriptstyle{F}$}\psfrag{x}{$\scriptstyle{x}$}
\psfrag{q}{$\scriptstyle{Q}$}
\includegraphics[width=3.8cm]{Efig1.eps}}
\] 
is equal to 
\[
(-1)^d \, \sum_{i=0}^{n} \, \mu_i
\times\qquad
\parbox{3cm}{\psfrag{x}{$\scriptstyle{x}$}\psfrag{y}{$\scriptstyle{y}$}
\psfrag{a}{$\scriptstyle{a}$}\psfrag{b}{$\scriptstyle{b}$}
\psfrag{d}{$\scriptstyle{2d}$}
\psfrag{m}{$\scriptstyle{4i}$}\psfrag{e}{$\scriptstyle{d}$}
\psfrag{f}{$\scriptstyle{d-2i}$}\psfrag{k}{$\scriptscriptstyle{2d-2i}$}
\includegraphics[width=3cm]{Efig2.eps}}\qquad\qquad .
\]

\bigskip \bigskip 

Differentiating out the variables $a$, $b$, $x$ and $y$, we get the
diagrammatic identity
\[
\sum_{i=0}^{n}\mu_i \quad 
\times\qquad\parbox{3.6cm}{\psfrag{d}{$\scriptstyle{2d}$}
\psfrag{m}{$\scriptstyle{4i}$}\psfrag{e}{$\scriptstyle{d}$}
\psfrag{f}{$\scriptstyle{d-2i}$}\psfrag{k}{$\scriptscriptstyle{2d-2i}$}
\includegraphics[width=3.6cm]{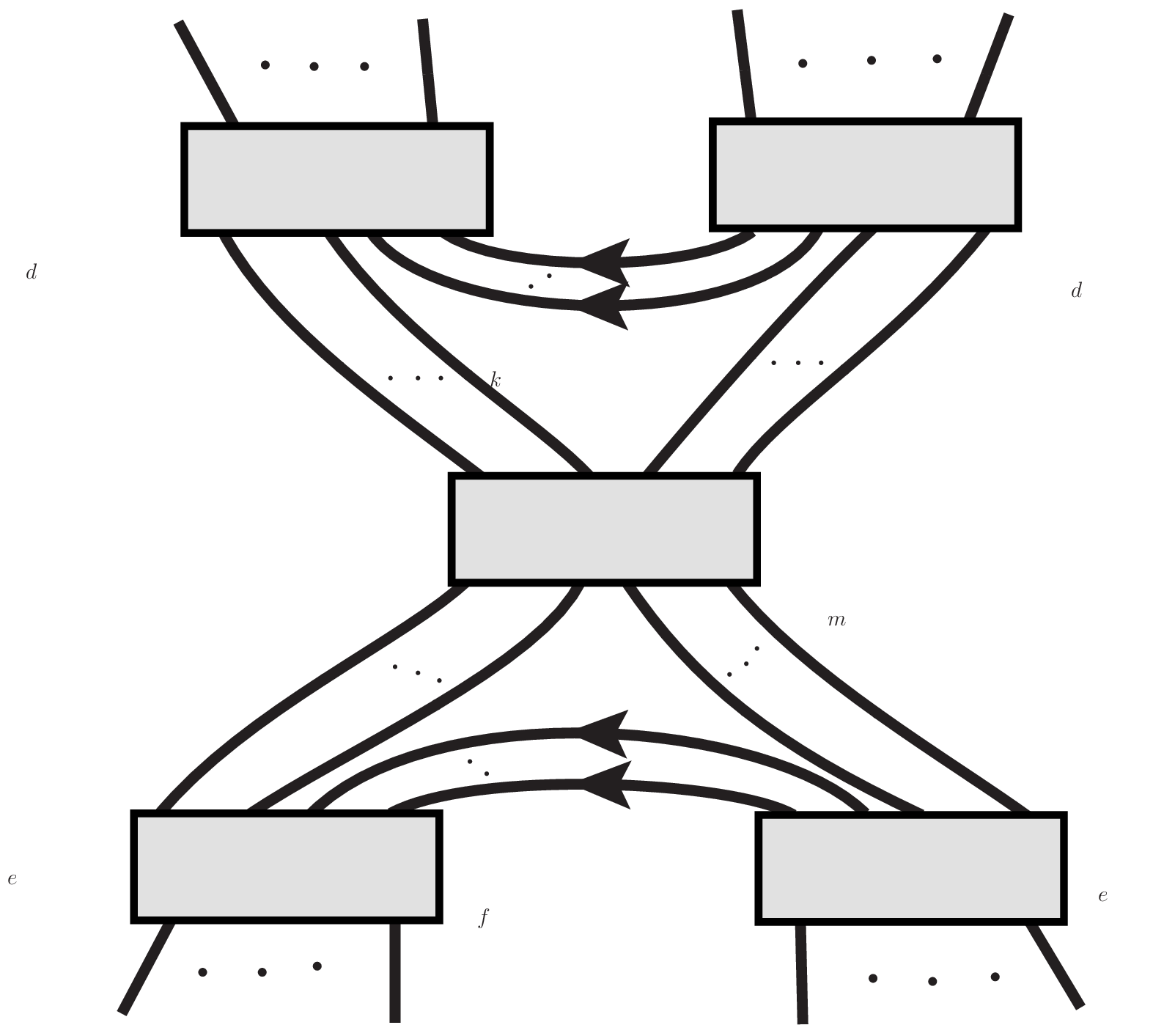}}
\qquad\ =0
\]
\vskip 0.8cm
\noindent
which holds for any choice of the $6d$ external leg indices in $\{1,2\}$.
Choose any $j$ in the range $0\le j\le n$. Contract on top with
\[
\parbox{3.6cm}{\psfrag{a}{$\scriptstyle{1}$}
\psfrag{m}{$\scriptstyle{4j}$}\psfrag{e}{$\scriptstyle{2d}$}
\psfrag{f}{$\scriptstyle{2d-2j}$}
\includegraphics[width=3.6cm]{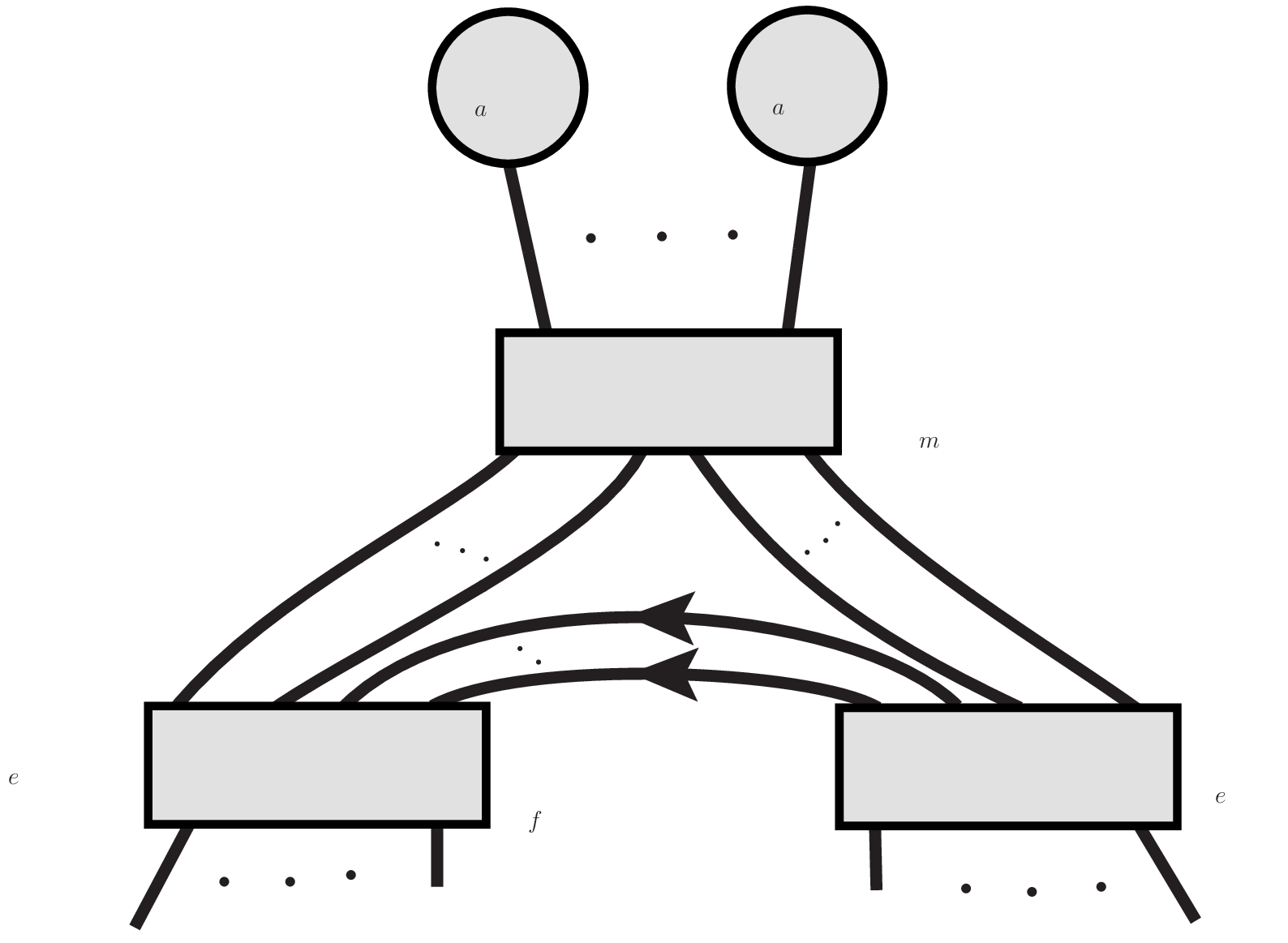}}
\]
and on the bottom with the same diagram upside down and $2d$ replaced by $d$.
Using the orthogonality relation~\cite[Eq. 13]{spinnet}
twice, we see that the sum reduces to the single term with $i=j$, which
gives $\mu_j$ times a nonzero coefficient. 
This extraction procedure shows that all $\mu$'s must vanish. \qed

There are several other ways to prove the last lemma. For instance, another possibility is to invoke 
Gordan's result to the effect that a joint covariant in $F$ and $Q$ (such 
as $\widetilde{\mathcal{T}}_i$) is a linear combination of iterated transvectants 
of the form 
\[
(Q,\cdots(Q,(Q,F)_{k_1})_{k_2}\cdots)_{k_d}\ .
\]
If we use Proposition \ref{proposition.omega} repeatedly 
while keeping track of the degrees, then it follows that the $\widetilde{\mathcal{T}}_i$ span
$\mathcal{M}_{2d}^{SL_2}$. 

\begin{Proposition}
The system $\SYS(d)$ is equivalent to saying that
the element $\psi=\sum_{i=0}^{n} \, z_i \, \mathcal{T}_i$ of
$\mathcal{M}_{d}^{SL_2}$ satisfies $\psi\circ\psi=\ID$.
\end{Proposition}
\demo
By bilinearity, 
\[
\psi\circ\psi=\sum\limits_{0 \le i, j \le n} \, z_i \, z_j\, 
\mathcal{T}_i\circ\mathcal{T}_j\ ,
\]
and from the calculation in \S\ref{section.Sd}, 
\[
\mathcal{T}_i\circ\mathcal{T}_j=
\sum_{k=0}^n \, \alpha_{i,j}^{(2k)} \; \widetilde{\mathcal{T}}_k\ .
\]
Therefore, $\SYS(d)$ is equivalent to saying that the coordinates 
of $\psi\circ\psi$ with respect to the basis
$\{\widetilde{\mathcal{T}}_k: 0\le k\le n\}$ are $1,0,\ldots,0$.
Since $\ID=\widetilde{\mathcal{T}}_0$, the claim follows. \qed

\subsection{ }
We will now solve the system via the introduction of more convenient bases $\mathcal{O}$
and $\widetilde{\mathcal{O}}$ for both spaces involved.
The basic identity we will need is
\begin{equation}
(Q \, \epsilon)^2=\frac{\Delta_Q}{4}Id
\label{mastereq}
\end{equation}
for the $2\times 2$ matrices
\[
Q=\left(\begin{array}{cc}
q_0 & q_1 \\
q_1 & q_2
\end{array} \right) \quad \text{and} \quad 
\epsilon=
\left(\begin{array}{cc}
0 & 1 \\
-1 & 0
\end{array}
\right)\ .
\]
Identity (\ref{mastereq}) is simply Cramer's rule for the matrix $Q$.
For any given quadratic $Q$ with complex coefficients let us make a choice of square root
$\sqrt{\Delta_Q}$. (How this choice varies with $Q$ will be irrelevant to the 
following discussion.) This allows us to factor the matrix equation (\ref{mastereq}) as
\begin{equation}
M_{+} M_{-}=M_{-} M_{+}=0, 
\label{MMortho} \end{equation}
where
\[
M_{\pm}=\frac{\sqrt{\Delta_Q}}{2}\pm Q \, \epsilon\ .
\]
Note that we have equations 
\[
M_{+}^2=\sqrt{\Delta_Q} M_{+}\qquad\text{and}\qquad
M_{-}^2=\sqrt{\Delta_Q} M_{-}\ .
\]
We will use the graphical notation
\[
\parbox{1.2cm}{\psfrag{i}{$\scriptstyle{i}$}
\psfrag{j}{$\scriptstyle{j}$}\psfrag{p}{$\scriptstyle{+}$}
\includegraphics[width=1.2cm]{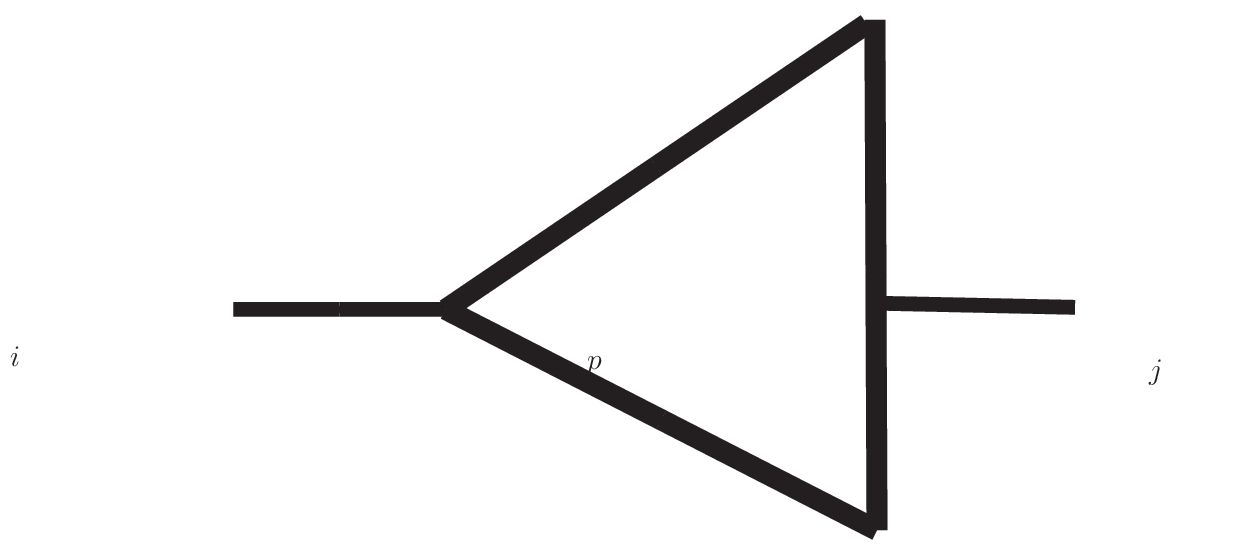}}
\ \ =(M_{+})_{ij}
\]
and likewise for $M_{-}$. For any $i$ in the range $0\le i\le d$, define the binary form
\begin{equation}
\mathcal{N}_i(F,Q)=\parbox{3cm}{\psfrag{f}{$\scriptstyle{F}$}
\psfrag{m}{$\scriptscriptstyle{-}$}\psfrag{p}{$\scriptscriptstyle{+}$}
\psfrag{i}{$\scriptstyle{i}$}\psfrag{x}{$\scriptstyle{x}$}
\psfrag{d}{$\scriptstyle{d-i}$}
\includegraphics[width=3cm]{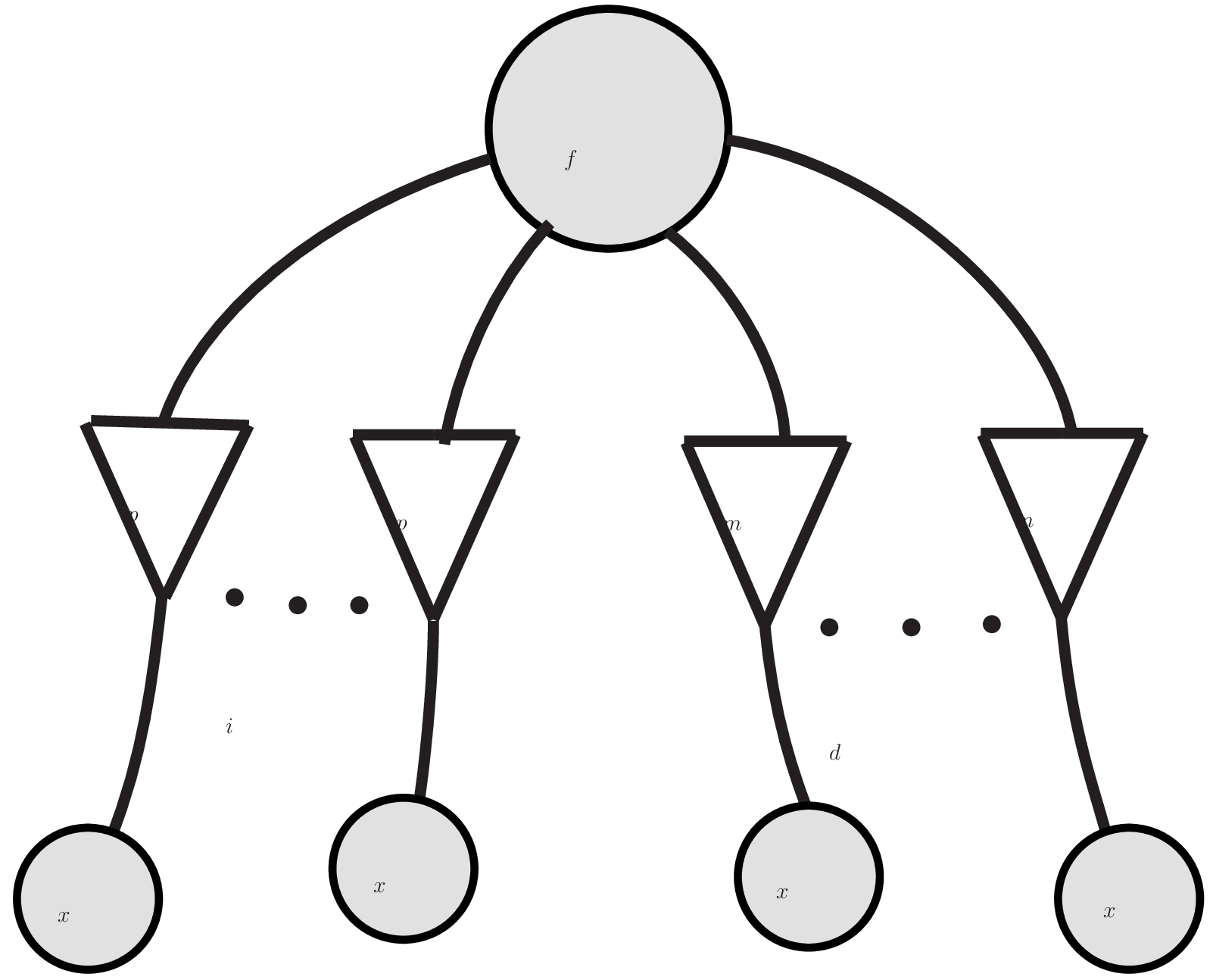}}\qquad\ .
\label{N.equation}
\end{equation}
The key observation is the following.
\begin{Lemma}\label{masterlemma}
For a fixed $Q$, the $\mathcal{N}$ can be seen as linear operators on binary forms $F$
of order $d$, which satisfy the relation 
\[
\mathcal{N}_i \, \mathcal{N}_j=\frac{\delta_{ij}}{\binom{d}{i}} \, 
\Delta_{Q}^{\frac{d}{2}} \, \mathcal{N}_i\ . \]
\end{Lemma}
\demo
Indeed,
\[
\mathcal{N}_i(\mathcal{N}_j(F,Q),Q)=
\parbox{2.8cm}{\psfrag{f}{$\scriptstyle{F}$}
\psfrag{m}{$\scriptscriptstyle{-}$}\psfrag{p}{$\scriptscriptstyle{+}$}
\psfrag{i}{$\scriptstyle{i}$}\psfrag{x}{$\scriptstyle{x}$}
\psfrag{d}{$\scriptstyle{d-i}$}
\psfrag{j}{$\scriptstyle{j}$}\psfrag{e}{$\scriptstyle{d-j}$}
\includegraphics[width=2.8cm]{Efig7.eps}}\qquad\qquad .
\]

\bigskip \bigskip \bigskip 

\noindent
When expanding the symmetriser, we see that because of (\ref{MMortho}), an $M_{+}$ can 
only contract to an $M_{+}$, and likewise for $M_{-}$. This forces $i=j$, and the binomial
accounts for the probability of having the branching permutation connect things properly,
i.e., the $M_{+}$ with the $M_{+}$ and the $M_{-}$ with the $M_{-}$. \qed

\subsection{ }
Now let, for $0\le i\le n$,
\begin{equation}
\mathcal{O}_i(F,Q)=\mathcal{N}_i(F,Q)+(-1)^{d} \, \mathcal{N}_{d-i}(F,Q)\ .
\label{O.convenient} \end{equation}
We will show that the latter expression only features even powers of $\sqrt{\Delta_Q}$, 
and therefore gives well defined elements in $\mathcal{M}_{d}^{SL_2}$.

If one expands the sums giving the $M_{\pm}$ matrices within the graphical formula
(\ref{N.equation}), then one obtains
\[
\mathcal{N}_i(F,Q)=
\sum_{p=0}^{i}\sum_{q=0}^{d-i}
\left(
\begin{array}{c} i \\ p 
\end{array}
\right)
\left(
\begin{array}{c} d-i \\ q 
\end{array}
\right)
\frac{(-1)^q\ \Delta_Q^{\frac{d-p-q}{2}}}{2^{d-p-q}}
\]
\begin{equation}
\times\qquad
\parbox{3cm}{\psfrag{f}{$\scriptstyle{F}$}
\psfrag{m}{$\scriptstyle{p+q}$}\psfrag{n}{$\scriptstyle{d-p-q}$}
\psfrag{q}{$\scriptstyle{Q}$}\psfrag{x}{$\scriptstyle{x}$}
\includegraphics[width=3cm]{Sec6fig2.eps}}\qquad\ .
\label{N.expansion}
\end{equation}
\vskip 0.6cm
Writing the same expansion for $\mathcal{N}_{d-i}(F,Q)$
with $p$ and $q$ exchanged, and adding the two contributions, we find that 
\[
\mathcal{O}_i(F,Q)=
\sum_{p=0}^{i}\sum_{q=0}^{d-i}
\left(
\begin{array}{c} i \\ p 
\end{array}
\right)
\left(
\begin{array}{c} d-i \\ q 
\end{array}
\right)
\frac{((-1)^q+(-1)^{p+d})\ \Delta_Q^{\frac{d-p-q}{2}}}{2^{d-p-q}}
\]
\[
\times\qquad
\parbox{3cm}{\psfrag{f}{$\scriptstyle{F}$}
\psfrag{m}{$\scriptstyle{p+q}$}\psfrag{n}{$\scriptstyle{d-p-q}$}
\psfrag{q}{$\scriptstyle{Q}$}\psfrag{x}{$\scriptstyle{x}$}
\includegraphics[width=3cm]{Sec6fig2.eps}}\qquad\ .
\]
\vskip 0.6cm
However, letting $s=p+q$, we see that $(-1)^q+(-1)^{p+d}$ vanishes unless $d-s$ is even.
As a result, introducing a new index $j$ for $\frac{d-s}{2}$, after simplification we get
\begin{equation}
\mathcal{O}_i(F,Q)=
\sum_{j=0}^n \, o_{i,j} \, \mathcal{P}_j(F,Q), \label{O.def} \end{equation}
with
\begin{equation}
o_{i,j}=
\frac{1}{2^{2j-1}} \, \sum_{p=0}^{i} \, \sum_{q=0}^{d-i}
\binom{i}{p} \, \binom{d-i}{q} \, \bbone\{p+q=d-2j\}\ (-1)^{q}\ .
\label{O.coeffs} \end{equation}
We can now use (\ref{O.def}) and (\ref{O.coeffs})
as the definition of the collection $\{\mathcal{O}_i : 0\le i\le n\}$ in
$\mathcal{M}_{d}^{SL_2}$, and regard the earlier (\ref{O.convenient}) simply
as a convenient representation which involves a choice of $\sqrt{\Delta_Q}$, yet
yields an outcome which is independent of that choice.

\begin{Lemma}\label{Oisabasis}
The collection $\{ \mathcal{O}_i: 0\le i\le n\}$ is a linear 
basis of $\mathcal{M}_{d}^{SL_2}$.
\end{Lemma}
\demo
It is enough to show that it generates the $\mathcal{P}_i$. We will also
write the explicit transition matrix.
For nonsingular $Q$, one can write the $2\times 2$ matrix identities
\[
Id=\frac{M_{+}+M_{-}}{\sqrt{\Delta_Q}}, \quad 
Q\epsilon=\frac{M_{+}-M_{-}}{2}, 
\]
and then expand the corresponding sums in the graphical
representation (\ref{P.graphical.def}) of $\mathcal{P}_i$ as
we did earlier for $\mathcal{N}$. We get
\[
\mathcal{P}_i(F,Q)=\frac{1}{2^{d-2i}} \, 
\sum_{p=0}^{d-2i} \, \sum_{q=0}^{2i} \, 
\binom{d-2i}{p} \, \binom{2i}{q} \, (-1)^p \,  \mathcal{N}_{d-p-q}(F,Q). 
\]
Rewrite this as 
\begin{equation}
\mathcal{P}_i(F,Q)=\sum_{l=0}^d \, \Gamma_{i,l} \,  \mathcal{N}_{l}(F,Q), 
\label{PintermsofN} \end{equation}
where
\begin{equation}
\Gamma_{i,l}=\frac{1}{2^{d-2i}} \, \sum_{p=0}^{d-2i} \, \sum_{q=0}^{2i} \, 
\binom{d-2i}{p}  \, \binom{2i}{q} \, (-1)^p \, \bbone\{d-p-q=l\}.
\label{Gamma.formula} \end{equation}
Now write the same formula for $\Gamma_{i,d-l}$, and make a change 
of indices $p\ra d-2i-p$ and $q\ra 2i-q$. This gives the relation
\[
\Gamma_{i,d-l}=(-1)^d \, \Gamma_{i,l}\ .
\]
Hence, one can fold the long sum (\ref{PintermsofN}) into the shorter one
\[
\mathcal{P}_i(F,Q)=\sum_{l=0}^n m_l\ \Gamma_{i,l}\ \mathcal{O}_{l}(F,Q), 
\]
where $m_l$ has been defined as in~\S\ref{z.formula.sec}. \qed

\subsection{ }
We now define a new collection of elements 
$\{ \widetilde{\mathcal{O}}_i : 0\le i\le n \}$ of
elements in the target space $\mathcal{M}_{2d}^{SL_2}$.
As before,  let 
\begin{equation}
\widetilde{\mathcal{O}}_i(F,Q)=
\Delta_Q^{\frac{d}{2}} \, \mathcal{N}_i(F,Q)
+\Delta_Q^{\frac{d}{2}} \, \mathcal{N}_{d-i}(F,Q). 
\label{Otilde.def} \end{equation}
Now use the expansion (\ref{N.expansion}) on both terms, 
exchanging the role of $p$ and $q$ in the second term. This produces the factor
\[
(-1)^q+(-1)^p=2(-1)^q \, \bbone\{p+q\ {\rm even}\}, 
\]
which forces the featured powers of $\sqrt{\Delta_Q}$ to be even. Therefore,
\[
\widetilde{\mathcal{O}}_i(F,Q)= 
\sum_{j=0}^n \, \Theta_{i,j} \, \widetilde{\mathcal{P}}_j(F,Q)
\]
where
\[ \Theta_{i,j}= 2 \, \times
\sum_{p=0}^{i} \, \sum_{q=0}^{d-i} \, \binom{i}{p} \, 
\binom{d-i}{q} (-1)^q \, \bbone\{p+q=d-2j\}\ .
\]
This shows that the collection $\widetilde{\mathcal{O}}$ is well defined in
$\mathcal{M}_{2d}^{SL_2}$.

\begin{Lemma} 
The collection $\{\widetilde{\mathcal{O}}_i : 0 \le i \le n \}$ is a linear basis of
$\mathcal{M}_{2d}^{SL_2}$.
\end{Lemma}
\demo
We proceed as in the proof of Lemma \ref{Oisabasis}. Using the expansion of the sums
for the matrices $Id$ and $Q \, \epsilon$ in terms of $M_{+}$ and $M_{-}$,
we have 
\[
\widetilde{\mathcal{P}}_i(F,Q)=\sum_{l=0}^d \, \Upsilon_{i,l} \, 
\Delta_Q^{\frac{d}{2}} \,  \mathcal{N}_{l}(F,Q), \qquad 
(0 \le i \le n), \]
where
\[
\Upsilon_{i,l}= \frac{1}{2^{2i}} \, \times
\sum_{p=0}^{2i} \, \sum_{q=0}^{d-2i} \, \binom{2i}{p} \, 
\binom{d-2i}{q} \, (-1)^p \, \bbone\{d-p-q=l\}.
\]
Again, the change of indices $p\ra 2i-p$ and $q\ra d-2i-q$
shows the relation $\Upsilon_{i,l}=\Upsilon_{i,d-l}$.
As a result, one has the `folded' sum representation
\begin{equation}
\widetilde{\mathcal{P}}_i=\sum_{l=0}^n \, m_l \, \Upsilon_{i,l} \, 
 \widetilde{\mathcal{O}}_l, 
\label{Upsilon.transition} \end{equation}
and the $\widetilde{\mathcal{O}}$ span $\mathcal{M}_{2d}^{SL_2}$. \qed

\begin{Proposition}\label{O.mult}
For $0\le i,j\le n$, we have
\[
\mathcal{O}_i\circ\mathcal{O}_j=
\frac{\delta_{ij}}{m_i \binom{d}{i}} \; \widetilde{\mathcal{O}}_i\ .
\]
\end{Proposition}
\demo
With obvious notations, 
\[
\mathcal{O}_i\circ\mathcal{O}_j=
\mathcal{N}_i \, \mathcal{N}_j
+(-1)^d \, \mathcal{N}_i \, \mathcal{N}_{d-j}
+(-1)^d \, \mathcal{N}_{d-i}\, \mathcal{N}_j
+\mathcal{N}_{d-i} \, \mathcal{N}_{d-j}\ .
\]
Now apply Lemma \ref{masterlemma} and use (\ref{Otilde.def}),
bearing in mind that $i,j$ range from $0$ to $n=\lfloor \frac{d}{2}\rfloor$.
Indeed, $i$ cannot equal $d-j$ and vice versa, except when $i=j=n=\frac{d}{2}$, which 
requires $d$ to be even. This accounts for the discrepancy factor $m_i$. \qed

\subsection{ }
We now solve the system in the $\mathcal{O}, \widetilde{\mathcal{O}}$ bases.
An element $\psi=\sum_{l=0}^n \, \rho_l \, \mathcal{O}_l$ of $\mathcal{M}_{d}^{SL_2}$
satisfies $\psi\circ\psi=\ID$, if and only if
there exists an initial segment $(s_0,\ldots,s_n)$ of a sign sequence $s$, 
such that 
\[ \rho_l=s_l \, m_l \, \binom{d}{l} \]
for all $l$ in the range $0\le l \le n$.
Indeed, on the one hand, by Proposition~\ref{O.mult}, 
\[
\psi\circ\psi=\sum_{l=0}^n \; \frac{\rho_l^2}{m_l \, \binom{d}{l}} \, 
\widetilde{\mathcal{O}}_l\ .
\]
On the other hand,
\[
\ID=\widetilde{\mathcal{P}}_0=\sum_{l=0}^n \, m_l \,  \Upsilon_{0,l} \, 
\widetilde{\mathcal{O}}_l = 
\sum_{l=0}^n \, m_l \,  \binom{d}{l} \, \widetilde{\mathcal{O}}_l \, , 
\]
by (\ref{Upsilon.transition}) and the formula for $\Upsilon$.
Note that for given $F$ and $Q$ one can `unfold' the expression of such a solution
$\psi$ as
\begin{equation}
\psi(F,Q)=\sum_{l=0}^d \, s_l \, \binom{d}{l} \, \mathcal{N}_l(F,Q), 
\label{unfold.eq} \end{equation}
where we have used the full sign sequence $s$, as in Definition~\ref{defofs}.

\subsection{ }
We now derive the formula for the involutors $z$. This is a simple change of basis calculation.
For a sign sequence $s$, the corresponding solution $\psi$ can be written
\[
\psi=\sum_{l=0}^n \, s_l \, m_l \, \binom{d}{l} \, \mathcal{O}_l
=\sum_{l=0}^n \, \sum_{e=0}^n \, o_{l,e} \, \mathcal{P}_e, 
\]
by (\ref{O.def}). Thus,
\[
\psi=\sum_{l=0}^n \, \sum_{e=0}^n \, \sum_{i=e}^{n} \, o_{l,e} \, G_{e,i} \, \mathcal{T}_i
\]
by Proposition~\ref{triangprop}.
Finally, if we extract the coefficient $z_i$ of $\mathcal{T}_i$, use our formulae
for the $o$ and $G$ transition matrices and simplify, then we get the required 
the expression for $z(s)$ in \S\ref{z.formula.sec}.
This completes the proof of Theorem \ref{MainTheorem}. \qed 

\section{Remaining computations} 
It only remains to prove Proposition~\ref{special.prop} and Theorem~\ref{can.form}. 
\subsection{Special sign sequences}\label{special.sec}
The geometric involutor $\sigma_Q$ corresponds to the element 
$2^d \, \mathcal{P}_0$ in $\mathcal{M}_{d}^{SL_2}$.
Therefore, by the proof of Lemma~\ref{Oisabasis}, 
\[
\sigma_Q(F)=2^d \, \sum_{l=0}^n \, m_l \, \Gamma_{0,l} \, \mathcal{O}_l(F,Q), 
\]
which reduces to 
\[
\sigma_Q(F)=\sum_{l=0}^n \, m_l  \, (-1)^{d-l} \binom{d}{l} \, \mathcal{O}_l(F,Q), 
\]
by (\ref{Gamma.formula}). 
Hence the corresponding sign sequence is $\gamma$, as defined in \S\ref{signseq.involutors}.

Now assume $d=2n$ is even. Then the improper involutor $(0,\ldots,0,1)$
corresponds to the element $\mathcal{T}_n=\mathcal{P}_n$ in 
$\mathcal{M}_{d}^{SL_2}$.
Using the explicit change of basis formulae in the proof of Lemma \ref{Oisabasis}, 
we find that 
\[
\mathcal{P}_n=\sum_{l=0}^n \, m_l \, \Gamma_{n,l} \, \mathcal{O}_l 
= \sum_{l=0}^n \, m_l  \binom{d}{l} \, \mathcal{O}_l\ .
\]
Hence, the corresponding sign sequence is $(+,\ldots,+)$.
This concludes the proof of Proposition~\ref{special.prop}. \qed 

\subsection{Canonical forms}\label{can.form.proof}
We now prove Theorem~\ref{can.form}. Let $z$ denote the involutor 
corresponding to a sign sequence $s$. Using (\ref{unfold.eq}), we have
\[
\sigma_{Q,z}(F)=\sum_{i=0}^d \, s_i \, \binom{d}{i} \, \mathcal{N}_i(F,Q) . 
\]
Let us specialize to $Q=x_1\, x_2$, and choose $\sqrt{\Delta_Q}=1$.
Then
\[ M_{+}=\left(\begin{array}{cc} 0 & 0 \\ 0 & 1 
\end{array}\right) \qquad \text{and} \qquad
M_{-}=\left(\begin{array}{cc} 1 & 0 \\ 0 & 0 
\end{array}\right)\ . \]
Therefore,
\[
\mathcal{N}_i(F,Q)=\ \ 
\parbox{3cm}{\psfrag{f}{$\scriptstyle{F}$}
\psfrag{i}{$\scriptstyle{i}$}\psfrag{d}{$\scriptstyle{d-i}$}
\psfrag{1}{$\scriptstyle{1}$}\psfrag{2}{$\scriptstyle{2}$}
\psfrag{x}{$\scriptstyle{x}$}
\includegraphics[width=3cm]{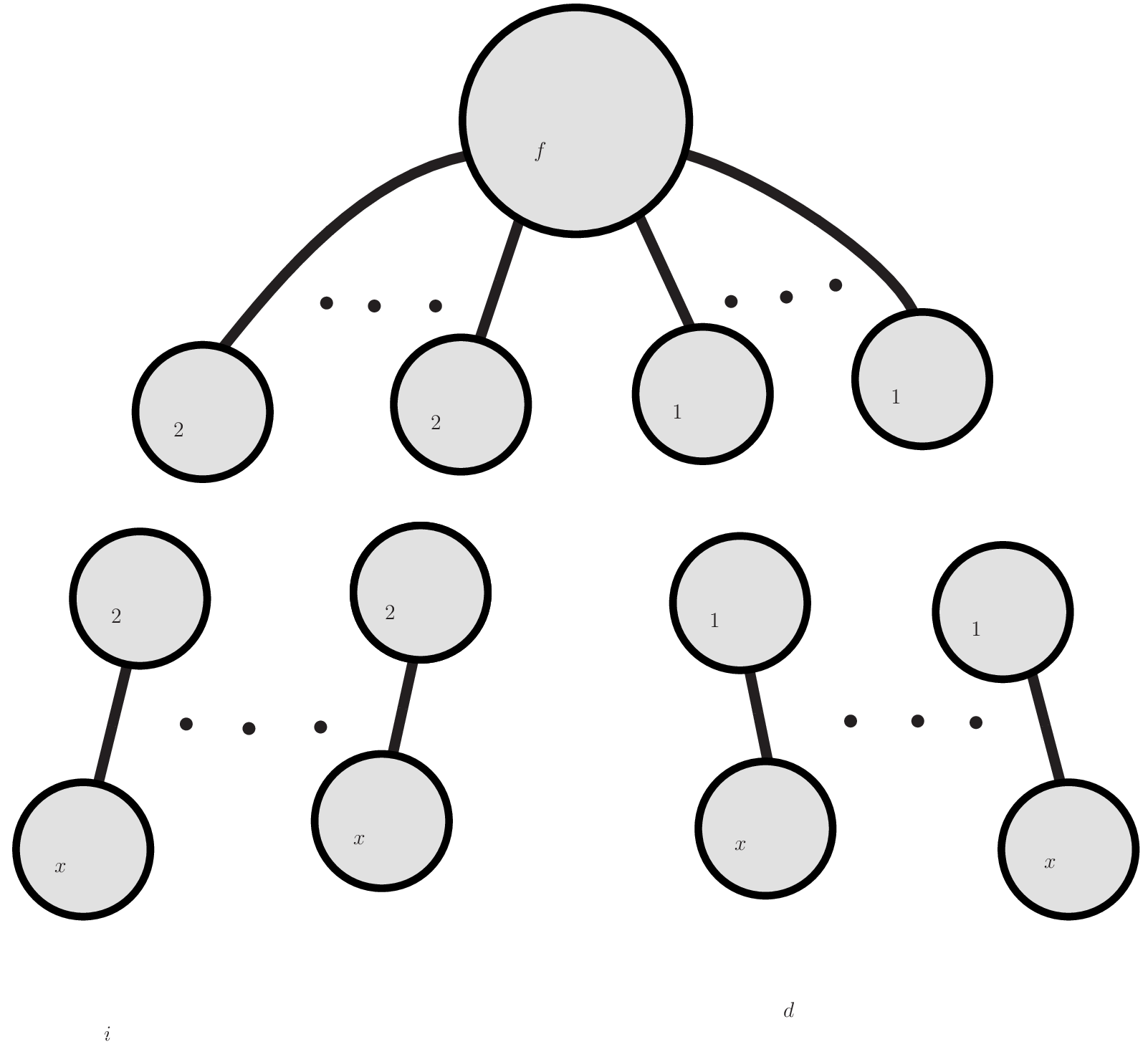}}
\qquad=
a_i\ x_1^{d-i}x_2^i, 
\]
where $F=(a_0,\ldots,a_d\cbrac x_1,x_2)$. 
If $d$ is even, then (\ref{condition.tcY}) is equivalent to
\[
(s_0 \, a_0,\ldots,s_d \, a_d\cbrac x_1,x_2)^d =F\ .
\]
Whereas, if $d$ is odd then (\ref{condition.tcY})
is equivalent to $G_{+} G_{-}=0$, 
where
\[
G_{\pm}=(s_0 \, a_0,\ldots,s_d \, a_d\cbrac x_1,x_2)^d  \pm F\ .
\]
This concludes the proof of Theorem~\ref{can.form}. \qed

\bigskip 

\noindent {\sc Acknowledgements:} {\small 
The second author has received financial support from NSERC, Canada. We are grateful to 
Daniel Grayson and Michael Stillman, the authors of Macaulay-2.}

\vspace{1cm} 

\centerline{--} 

\vspace{1cm} 

\parbox{7cm} 
{Abdelmalek Abdesselam \\
Department of Mathematics, \\
University of Virginia, \\
P. O. Box 400137, \\
Charlottesville, VA 22904-4137, \\
USA.\\
{\tt malek@virginia.edu}} 
\hfill 
\parbox{7cm} 
{Jaydeep Chipalkatti \\ 
Department of Mathematics, \\ 
Machray Hall, \\ 
University of Manitoba, \\ 
Winnipeg, MB R3T 2N2, \\ Canada. \\ 
{\tt chipalka@cc.umanitoba.ca}}

\end{document}